\documentclass[3p]{elsarticle}

\usepackage{lineno,hyperref}
\usepackage{placeins}\usepackage{tikz}\usetikzlibrary{shapes,shadings}
\usepackage{pgfplots}
\pgfplotsset{compat=newest}
\pgfplotsset{plot coordinates/math parser=false}
\newlength\figureheight
\newlength\figurewidth

\usepackage{tabularx}

\usepackage{amsmath}
\usepackage{amsfonts}

\definecolor{myred}{RGB}{128,0,0}
\definecolor{myblue}{RGB}{0,0,128}
\definecolor{mygreen}{RGB}{0,128,128}

\modulolinenumbers[1]
\usepackage[normalem]{ulem}
\def\checkmark{\tikz\fill[scale=0.4](0,.35) -- (.25,0) -- (1,.7) -- (.25,.15) -- cycle;}
\newcommand{\addtxt}[1]{{#1}}
\newcommand{\removetxt}[1]{}

\journal{Elsevier}

\begin{document}
\begin{frontmatter}

\title{A Class of Analytic Solutions for Verification and Convergence Analysis of Linear and Nonlinear Fluid-Structure Interaction Algorithms}

\author[Stuttgart]{Andreas Hessenthaler\corref{corauth}}
\ead{hessenthaler@mechbau.uni-stuttgart.de}

\author[KCL]{Maximilian Balmus\corref{corauth}}
\ead{maximilian.balmus@kcl.ac.uk}

\author[Stuttgart]{Oliver R\"ohrle}

\author[UM,KCL]{David Nordsletten}

\address[Stuttgart]{Institute for Modelling and Simulation of Biomechanical Systems,
University of Stuttgart, Pfaffenwaldring 5a, 70569~Stuttgart, Germany}

\address[KCL]{School of Biomedical Engineering and Imaging Sciences, King's College London, 4th FL Rayne Institute, St.~Thomas~Hospital, London, SE1 7EH}
\address[UM]{Department of Biomedical Engineering and Cardiac Surgery, University of Michigan, NCRC B20, 2800 Plymouth Rd, Ann~Arbor, 48109}

\cortext[corauth]{Authors are acknowledged as joint-first and joint-corresponding authors.}

\begin{abstract}
Fluid-structure interaction (FSI) problems are pervasive in the computational engineering community.
The need to address challenging FSI problems has led to the development
of a broad range of numerical methods addressing a variety of application-specific demands.
While a range of numerical and experimental benchmarks are present in the literature, few solutions are available that enable both verification and spatiotemporal convergence analysis.
In this paper, we introduce a class of analytic solutions to FSI problems involving shear in channels and pipes.
Comprised of 16 separate analytic solutions, our approach is permuted to enable progressive verification and analysis of FSI methods and implementations,
in two and three dimensions, for static and transient scenarios as well as for linear and hyperelastic solid materials.
Results are shown for a range of analytic models exhibiting progressively complex behavior.
The utility of these solutions for analysis of convergence behavior is further demonstrated using a previously published monolithic FSI technique.
The resulting class of analytic solutions addresses a core challenge in the development of novel FSI algorithms and implementations, providing a progressive testbed for verification and detailed convergence analysis.
\end{abstract}


\end{frontmatter}
\newif\iffigure
\figuretrue
\section{Introduction}\label{introduction-sec}

Within computational engineering, examples of fluid-structure interaction (FSI) are pervasive and represent an increasingly important set of problems.
Addressing the disparate requirements of different FSI applications, the scientific community has responded by generating a broad range of numerical methods.
From now classic arbitrary Lagrangian-Eulerian (ALE) boundary fitted approaches~\cite{hirt1974arbitrary,donea1982arbitrary,hughes1981lagrangian},
to space-time ALE variational multiscale (ALE VMS)~\cite{tezduyar1991stabilized,tezduyar1992new,takizawa2011multiscale},
unified continuum methods~\cite{HoffmanJanssonStoeckli2011,JanssonDegirmenciHoffman2017},
immersed boundary methods~\cite{peskin1973flow,peskin2002immersed,mittal2005immersed},
fictitious domain methods~\cite{glowinski1994fictitious,glowinski1999distributed},
immersed structural potential methods~\cite{gil2010immersed,gil2013enhanced},
and overlapping domain methods~\cite{StegerDoughertyBenek1983,StegerBenek1987,ChesshireHenshaw1990,HouzeauxCodina2003,WallGamnitzerGerstenberger2008,BalmusMassingHoffmanRazaviNordsletten2019_preprint} (to name a few) many innovative FSI techniques have been introduced that address application-specific challenges.
Trailing this methodological development enabling complex simulations was a boom in applications, further pushing the numerical envelope to accommodate bigger problems with more physical models that emulate the complex fluid-structure dynamics of real-world systems.

Part and parcel to the development of FSI techniques, which often involve bespoke or in-house codes,
comes the consistent need for verification and validation.
Extending the tradition of documented numerical results established in fluid mechanics
(e.g., lid-driven cavity~\cite{schreiber1983driven,schreiber1983spurious,kim1985application}),
FSI verification problems have been developed in both
two and three dimensions~\cite{GhattasLi1995,Wall1999,Mok2001,TurekHron2006,BatheLedezma2007}
(a more exhaustive list can be found in \cite{HessenthalerGaddumHolubSinkusRoehrleNordsletten2017},
Table 1), often involving elastic structures immersed in a steady or periodic flow.
These methods have been extensively used to compare
results~\cite{HeilHazelBoyle2008,TurekHronRazzaqWobkerSchaefer2011}
across codes and provide a measure of numerical consistency.
Similarly, experiments have been proposed~\cite{BertramTscherry2006,
GomesLienhart2006,IdelsohnMartiSoutoiglesiasOnate2008,
GomesLienhart2010, NayerKalmbachBreuerSicklingerWuechner2014,HessenthalerRoehrleNordsletten2017}
that provide experimental data for validation.

While these approaches provide important mechanisms for observing numerical performance and fidelity, the above benchmark problems and experiments typically present challenges for assessing the accuracy and convergence of FSI methods and implementations.
Given the complexity of FSI solutions, spatiotemporal convergence analysis presents a critical, yet difficult, aspect of FSI method verification.
The method of \addtxt{so-called} manufactured solutions~\cite{steinberg1985symbolic,roache2002code,salari2000code}, whereby a compatible solution is selected and appropriate forcing terms are prescribed, provides a potential avenue for addressing these needs.
However, implementation of forcing terms can be complex, code specific, and prone to spatiotemporal errors introduced through their numerical inclusion.
\addtxt{It further requires software to integrate with these libraries or tools,
and how source terms are integrated into code numerically can have detrimental impacts on observed convergence rates (or even mask problems).}
An alternative is to provide analytic solutions, as are used in other fields, for spatiotemporal convergence
analysis~\cite{ethier1994exact}.
In this case, the fluid and solid state variables are known to satisfy the FSI problem, enabling comprehensive numerical assessment extending to the limits of machine precision.
While promising, in theory, the availability of analytic solutions for FSI remain severely limited due, in part, to the complexity of the system.

In this work, we aim to address this core need for method verification and spatiotemporal convergence analysis by introducing a novel class of analytic solutions for FSI.
Building from Womersley's solution for pulsatile fluid flow~\cite{womersley1955method}, we introduce analytic solutions that characterize fluid flow and solid motion in shear.
To enable incremental testing, solutions are derived for two (channel) and three dimensions (tube) for transient and steady fluids and solids.
Additionally, we introduce analytic solutions for both linear elastic and nonlinear hyperelastic (neo-Hookean) solid models.
In total, 16 analytic solutions are presented, exhibiting a wide range of solution complexity in both space and time.
Code implementing these solutions is provided (see \emph{Supplementary Material}), enabling comparisons with numerical results and evaluation of FSI solutions for different parameter combinations.
Derived solutions are subsequently used to examine spatiotemporal convergence
and accuracy of our previously published FSI method and implementation~\cite{NordslettenKaySmith2010,LeeEtAl2016},
demonstrating the efficacy of these analytic solutions for providing meaningful analysis.

In what follows, we begin by outlining the general FSI boundary value problem in both two
and three dimensions and deriving their analytic solutions (Section~\ref{methodology-sec}).
Details of the previously published FSI method~\cite{NordslettenKaySmith2010}
are briefly reviewed (Section~\ref{fsi-implementation-sec})
along with details of the spatiotemporal convergence analysis protocol (Section~\ref{num-experiments}).
In Section~\ref{results-sec}, results illustrating the behavior of the analytic solutions for two specific cases
(transient two-dimensional fluid / linear solid, transient three-dimensional fluid / nonlinear solid)
are demonstrated, accompanied by convergence results showing expected numerical behavior with spatiotemporal refinement.
These results are discussed in-depth in Section~\ref{discussion-sec},
followed by concluding remarks in Section~\ref{conclusion-sec}.

\section{Methodology} \label{methodology-sec}
\subsection{Derivations}\label{derivations-seq}
In the interest of condensing the material presented in this work,
we limit the number of fully derived cases to two major
categories: 2D linear and 3D nonlinear solid and fluid mechanics.\footnote{For
a complete presentation of other cases
(i.e.\ 2D nonlinear and 3D linear solid and fluid mechanics),
please consult the \emph{Supplementary Materials SM\pageref{online-supplement-sec}}.}
Each category is further divided to include all possible permutations of
quasi-static and transient conditions.
In the following, we derive a general solution for the fluid and solid problems separately.
Subsequently, depending on the combination of temporal behaviors,
we derive a unique FSI solution which satisfies the kinematic
and traction coupling conditions.
\subsubsection{Linear fluid/solid in two dimensions}
\label{derivation-2D-sec}

We begin the series of derivations with the simplest case presented in this work:
the interaction between a pulsatile-flowing liquid in a two dimensional channel
and a linear elastic solid, on the top and bottom, which undergoes shear deformation,
see Figure~\ref{domain-2D-fig}.
One of the aims of this section is to help the reader get familiarized
with the notation, main assumptions and steps taken throughout the derivation, all of
which are reused and adapted in the more complex cases.
For simplicity, we assume that the physics are symmetrical with respect to the
axis running along the middle of the channel, allowing us to limit the problem domain
to one of the halves.
In the reference space,
the solid and fluid domains are represented by two quadrilaterals, $\Omega^0_f$ and $\Omega^0_s$,
which are separated by an interface $\Gamma^{\lambda}$.
Their length, the width of the channel and the width of the domain
 are denoted by $L$, $H_i$ and $H_o$. The top boundary
represents the outer wall of the solid and is denoted by $\Gamma^W_s$, while the bottom
one corresponds to the axis of symmetry of the channel and we refer to it as
$\Gamma^W_f$. The left and right boundaries are the inlet ($\Gamma_f^I$ and $\Gamma_s^I$)
and the outlet ($\Gamma_f^O$ and $\Gamma_s^O$).

The FSI problem takes the following general form:
\begin{align}
\rho_f \partial_t \boldsymbol{v}_f - \nabla \cdot \boldsymbol{\sigma}_f &= \boldsymbol{0} && \text{in } \Omega_f,
\label{eq:strong_form_linear_fluid_momentum_balance}\\
\nabla \cdot \boldsymbol{v}_f &= 0 && \text{in } \Omega_f,
\label{eq:strong_form_linear_fluid_mass_balance}\\
[\boldsymbol{v}_f]_y & = 0 && \text{on } \Gamma_f^I \cup \Gamma_f^O \cup \Gamma_f^W, 
\label{eq:strong_form_linear_Dirichlet_condition_fluid} \\
(\boldsymbol{\sigma}_f \cdot \boldsymbol{n}_f) \cdot \boldsymbol{e}_x &=
[\boldsymbol{t}_f]_x && \text{on } \Gamma_f^I \cup \Gamma_f^O,
\label{eq:strong_form_linear_Neumann_condition_fluid} \\
\frac{\partial \boldsymbol{v}_f}{\partial y} & = \boldsymbol{0}
&& \text{on } \Gamma^W_f, \\
\label{eq:strong_form_smoothness_condition}
\boldsymbol{v}_f (\cdot, 0) &= \boldsymbol{v}_f^0 && \text{in } \Omega_f (0), \\
\hline
\rho_s \partial_{tt} \boldsymbol{u}_s -\nabla \cdot \boldsymbol{\sigma}_s &= \boldsymbol{0} && \text{in } \Omega_s,
\label{eq:strong_form_linear_solid_momentum_balance}\\
\nabla \cdot \partial_t \boldsymbol{u}_s &= 0 && \text{in } \Omega_s,
\label{eq:strong_form_linear_solid_mass_balance}\\
\boldsymbol{u}_s (\cdot, t) &= \boldsymbol{0} && \text{on } \Gamma^W_s,
\label{eq:strong_form_wall_condition}\\
[\boldsymbol{u}_s]_y & = 0 && \text{on } \Gamma_s^I \cup \Gamma_s^O, \\
(\boldsymbol{\sigma}_s \cdot \boldsymbol{n}_s) \cdot  \boldsymbol{e}_x
&= [\boldsymbol{t}_s ]_x
&& \text{on } \Gamma_s^I \cup \Gamma_s^O, \\
\boldsymbol{u}_s (\cdot, 0) &= \boldsymbol{u}_s^0 && \text{in } \Omega_s (0), \\
\boldsymbol{v}_s (\cdot, 0) & = \boldsymbol{v}_s^0 && \text{in } \Omega_s (0), \\
\hline
\boldsymbol{\sigma}_f \cdot \boldsymbol{n}_f + \boldsymbol{\sigma}_s \cdot \boldsymbol{n}_s &= \boldsymbol{0} && \text{on } \Gamma^\lambda, \label{eq:strong_form_linear_dynamic_constraint}\\
\boldsymbol{v}_f - \boldsymbol{v}_s &= \boldsymbol{0} && \text{on } \Gamma^\lambda,
\label{eq:strong_form_linear_kinematic_constraint}
\end{align}
with the Cauchy stress tensor for the fluid and solid defined as:
\begin{equation*}
	\boldsymbol{\sigma}_f = \mu_f [\nabla \boldsymbol{v}_f + \nabla^T \boldsymbol{v}_f] - p_f \boldsymbol{I}
	\quad
	\text{and}
	\quad
	\boldsymbol{\sigma}_s = \mu_s [\nabla \boldsymbol{u}_s + \nabla^T \boldsymbol{u}_s] - p_s \boldsymbol{I},
\end{equation*}
respectively. Naturally, individual transient and quasi-static cases may be
obtained by nullifying or assigning strictly positive values to the two density parameters,
$\rho_f$ and $\rho_s$. The solid stiffness and fluid viscosity are denoted by $\mu_s$ and $\mu_f$.
For now we will assume that
the initial value fields (i.e. $\boldsymbol{v}_f^0$, $\boldsymbol{u}_s^0$ and $\boldsymbol{v}_s^0$), and the
inlet and outlet surface tractions, (i.e. $\boldsymbol{t}_f$ and $\boldsymbol{t}_s$) are given.
As we shall see, knowing what these fields are is not prerequisite to the derivation of the analytical solutions.
Here, they are presented in anticipation of the weak form of the problem used in the numerical results section.
\removetxt{Also, note that the domains at time zero do not necessarily coincide with reference ones
(e.g., when $\boldsymbol{u}_s^0\neq\boldsymbol{0}$, $\Omega_s(0)\neq\Omega_s^0$).}

\addtxt{It should be noted that applying traction boundary conditions at the inlet of the fluid
problem, as in Eq.~\ref{eq:strong_form_linear_Neumann_condition_fluid}, will not inherently
lead to a stable formulation, see~\cite{hughes2005conservation,moghadam2011comparison}.
That being said, from our experience, stabilization was not required for this class of problems 
in order to get accurate solutions. 
Alternatively, one may achieve a more robust problem formulation by changing
the treatment of the boundary to Dirichlet or other stabilized forms. }

 In a similar fashion to the classical pulsatile channel flow problem,
 our  strategy is based on the assumptions that the fluid velocity and solid displacement
 are constant in the $x$-direction.
 Hence, the solid displacement and fluid velocity take the form:
\begin{equation*}
	\boldsymbol{v}_f = v_f(y,t) \boldsymbol{e}_x \quad \text{and} \quad
	\boldsymbol{u}_s = u_s(y,t) \boldsymbol{e}_x.
\end{equation*}

To obtain a general form of the fluid solution,
we first re-write the momentum balance equation
in~\eqref{eq:strong_form_linear_fluid_momentum_balance}
in component form and apply our knowledge of the flow behavior:
\begin{align}
	\rho_f\partial_t v_f = & - \frac{\partial p_f}{\partial x} + \mu_f  \frac{\partial^2 v_f}{\partial y^2} ,
	\label{eq:balance_fluid_lin_1} \\
	0 = & -\frac{\partial p_f }{\partial y} .
	 \label{eq:balance_fluid_lin_2}
\end{align}
From~\eqref{eq:balance_fluid_lin_2}, we see that $p_f$ is constant in $y$.
Furthermore, applying $\partial(\cdot)/\partial x$ to~(\ref{eq:balance_fluid_lin_1}) yields:
\begin{equation}
	-\frac{\partial^2 }{\partial x^2} p_f(x,t) = 0. \nonumber
\end{equation}
Thus, $p_f$ is a constant or linear function in $x$. We then suppose the velocity
and pressure are periodic and separable, where the pulsatile behavior
is defined by a single harmonic frequency, $\omega$, i.e.:
\begin{equation}
	v_f(y,t) = \Re\{  \text{v}_f(y)  e^{i\omega t}  \}
	\quad
	\text{and}
	\quad
	p_f(x,t) = \Re\{ P_f (L - x) e^{i\omega t} \},
	\label{eq:fluid_lin_pulse}
\end{equation}
where $P_f \in \mathbb{C}$ and $\text{v}_f:[0,H_i]\rightarrow \mathbb{C}$
denote the fluid pressure over the domain length $L$ and velocity amplitudes,
respectively.\footnote{The purpose of the complex pressure parameter is to allow for a temporal phase shift.
This applies to all subsequent $P$ type parameters.}
The first momentum balance equation reduces to:
\begin{equation}
	i\omega\rho_f \text{v}_f = P_f + \mu_f  \text{v}_f''.
\end{equation}
In this case, the general solution takes the form:
\begin{equation}
	\begin{aligned}
		\text{v}_f(y) &= -\frac{P_f}{2\mu_f} y^2 + c_2 y + c_1,  & (\rho_f = 0), \\
		\text{v}_f(y) &= -\frac{iP_f}{\rho_f  \omega} + c_1 e^{k_fy} + c_2 e^{-k_fy}, & (\rho_f > 0).
	\end{aligned}
	\label{eq:fluid_lin_general}
\end{equation}
where $k_f = \sqrt{i\rho_f\omega / \mu_f}$ (see Table \ref{nomenclature-tab} for a list of short-form constants).
Note that the integration constants are re-used here to simplify notation
and they are not related to each other.
Furthermore, based on the smoothness condition
in Equation~\eqref{eq:strong_form_smoothness_condition},
which is equivalent to $\text{v}_f' = 0$,
we can also conclude that $c_2 = 0$ (quasi-static case) or
$c_2 = c_1$ (transient case).

Moving on to the solid component of the problem,
the corresponding momentum balance equation can
be re-written to incorporate our assumptions:
\begin{align}
	\rho_s \partial_{tt} u_s = & - \frac{\partial p_s}{\partial x}  + \mu_s \frac{\partial^2 u_s}{\partial y^2} ,
	\label{eq:balance_solid_lin_1}   \\
	0 = & - \frac{\partial p_s}{\partial y} .
\end{align}
Assuming periodicity, separability and the same period length for both fluid and solid, we write:
\begin{equation}
	u_s (y,t) = \Re\{ \text{u}_s(y) e^{i\omega t} \}
	\quad
	\text{and}
	\quad
	p_s(x,t) = \Re\{P_s(L-x) e^{i\omega t}\}, \label{eq:solid_lin_pulse}
\end{equation}
where $P_s\in\mathbb{C}$ and $\text{u}_s:[R_i,R_0]\rightarrow\mathbb{C}$ are
the amplitudes of the solid solution.
Here, we skipped the derivation of $p_s$ as it is
analogous to that of $p_f$.
Consequently, Equation~\eqref{eq:balance_solid_lin_1} reduces to:
\begin{equation}
	-\rho\omega^2 \text{u}_s = P_s + \mu_s \text{u}_s ''. \nonumber
\end{equation}
resulting in the general solution for solid displacement:
\begin{equation}
	\begin{aligned}
		\text{u}_s(y) & = - \frac{P_s}{2\mu_s} y^2 + c_3 y + c_4, &(\rho_s = 0), \\
		\text{u}_s{(y)} &= - \frac{P_s}{\rho_s \omega^2} + c_3 \sin(k_s y) +c_4 \cos(k_s y), & (\rho_s > 0),
	\end{aligned}
	\label{eq:solid_lin_general}
\end{equation}
where $k_s = \sqrt{\rho_s \omega^2 / \mu_s}$.
The last unknown integration constants (i.e. $c_1$, $c_3$ and $c_4$)
can only be identified by ensuring that the coupling condition are satisfied.
In the following we identify their closed form depending on the different
temporal behavior combinations.

\textbf{Quasi-static fluid and quasi-static solid ($\rho_f = \rho_s = 0$)}

We begin the coupling process by expanding the traction
boundary condition in~\eqref{eq:strong_form_linear_dynamic_constraint}:
\begin{align}
	\mu_f \left. \frac{\partial \text{v}_f }{\partial y}  \right|_{y = H_i} -
	\mu_s \left.\frac{\partial  \text{u}_s }{\partial y} \right|_{y= H_i} = & 0, \label{eq:lin-qq-dyn-1} \\
	-\left. p_f \right|_{y=H_i} + \left. p_s \right|_{y=H_i} = & 0.\label{eq:lin-qq-dyn-2}
\end{align}
Based on equation~(\ref{eq:lin-qq-dyn-2})
it can easily be shown to shown that $P_f  = P_s = P$. Note that, as we will see,
this result is independent of the quasi-static/transient property of the FSI problem.
Furthermore, we substituting equations
\eqref{eq:fluid_lin_general} and \eqref{eq:solid_lin_general} into Equation~\eqref{eq:lin-qq-dyn-1}
we can see that
\begin{equation}
	-P H_i - \mu_s \left( -\frac{P H_i }{\mu_s}+ c_3\right) = 0
	\hspace{2mm}\Rightarrow \hspace{2mm}
	c_3 = 0.
\end{equation}
Considering the fixed wall constraint on the solid (i.e. $\text{u}_s(H_o) = 0$),
we can use Equation~\eqref{eq:solid_lin_general} to show that:
\begin{equation}
	c_4 = \frac{P H_o^2 }{2\mu_s}.
\end{equation}
Finally, by expanding the kinematic constraint in~\eqref{eq:strong_form_linear_kinematic_constraint}, we obtain the final unknown constant:
\begin{equation}
	c_1 = \frac{PH_i^2}{2\mu_f} + i\omega\frac{P}{2 \mu_s} (H_o^2 - H_i^2).
\end{equation}

\textbf{Quasi-static fluid and transient solid} ($\rho_f = 0$,~$\rho_s > 0$)

Based on our general formulations for the
fluid velocity and solid displacement in~\eqref{eq:fluid_lin_pulse}, \eqref{eq:fluid_lin_general},
\eqref{eq:solid_lin_pulse} and \eqref{eq:solid_lin_general}, the
kinematic interface condition~\eqref{eq:strong_form_linear_kinematic_constraint} reduces to
\begin{equation}
	c_1 - c_3 i \omega \sin(k_s H_i) - c_4  i \omega \cos(k_s H_i) = \frac{PH_i^2}{2\mu_f} - \frac{iP}{\rho_s \omega}.
	\label{eq:lin_qt_sys1}
\end{equation}
A second equation relating the integration constants can
be obtained from the condition~\eqref{eq:strong_form_wall_condition} for fixed outer walls:
\begin{equation}
	c_3 \sin(k_s H_o) + c_4 \cos(k_s H_o) = \frac{P}{\rho_s \omega^2}.
	\label{eq:lin_qt_sys2}
\end{equation}
The third and final equation can be obtained through the expansion of
shear component of the traction  equation~(\ref{eq:strong_form_linear_dynamic_constraint}), e.g.,
\begin{equation}
	c_3 \mu_s k_s \cos(k_s H_i) - c_4\mu_s k_s \sin(k_s H_i) = - P H_i.
	\label{eq:lin_qt_sys3}
\end{equation}
Finally, $c_1$, $c_3$ and $c_4$ can be determined by solving the system of equations formed by~\eqref{eq:lin_qt_sys1},
\eqref{eq:lin_qt_sys2} and \eqref{eq:lin_qt_sys3}. Provided there is a unique solution,
the three constants can be written as the following closed-form expressions:
\begin{align}
	c_1  & = \frac{PH_i^2}{2\mu_f} - \frac{iP}{\rho_s \omega}
	+ PH_i \frac{i\omega}{\mu_s k_s} \tan[k_s(H_o - H_i)]
	+ \frac{iP}{\rho_s \omega} \sec[k_s(H_o - H_i)], \\
	c_3 & =  \left[\frac{1}{\rho_s \omega^2} \sin(k_s H_i)
	- \frac{H_i}{\mu_s k_s} \cos(k_s H_o)\right] P \sec[k_s(H_o - H_i)], \\
	c_4 & = \left[ \frac{H_i}{\mu_s k_s} \sin(k_s H_o) +
	\frac{1}{\rho_s \omega^2} \cos(k_s H_i)\right] P \sec[k_s(H_o - H_i)].
\end{align}

This approach of finding the closed form of the constants by solving a system
of equations is repeated throughout the paper. Consequently, it should be noted
that the existence of unique solutions is conditioned by the non-singularity of
the system matrix. As we will see, identifying the appropriate set of parameters
which lead to a singular matrix is not a trivial problem and, as a consequence,
we will only perform this analysis for this specific case, which is more tractable.
Here, it can be shown that for a set of problem parameters with finite values,
the determinant of the system is null when:
\begin{equation}
\cos\left[k_s(H_o -H_i)\right] = 0. \nonumber
\end{equation}
Rearranging this, we find a series of resonance frequencies for which this is true:
\begin{equation}
\omega_n = \frac{(2n+1)\pi}{2(R_o - R_i)} \sqrt{\frac{\mu_s}{\rho_s}},
\text{  } n\in \mathbb{Z}.
\label{eq:resonance-condition}
\end{equation}

\textbf{Transient fluid and quasi-static solid} ($\rho_f > 0,~\rho_s =0$)

As in the previous problem permutation, the constants (i.e. $c_1$, $c_3$ and $c_4$)
can be identified by solving the system of equations derived from~\eqref{eq:strong_form_linear_dynamic_constraint},
\eqref{eq:strong_form_wall_condition} and \eqref{eq:strong_form_linear_kinematic_constraint}:
\begin{align}
\left(e^{k_fH_i} + e^{-k_fH_i}\right)c_1  - i \omega H_i  c_3 - i \omega c_4 &
= \frac{i \omega  P}{\rho_f \omega} -\frac{i\omega P H_i^2}{2\mu_s} , \\
H_o c_3 + c_4 & = \frac{P H_o^2}{2\mu_s} ,\\
\mu_f k_f \left(e^{k_f H_i} - e^{-k_f H_i} \right)c_1 - \mu_s c_3 & = - P H_i.
\end{align}
The resulting closed-form expressions for the three constants are:
\begin{align}
	c_1 & = \frac{\frac{i\mu_s P}{\rho_f \omega} + \frac{i\omega P}{2}(H_o^2 - H_i^2) - i\omega P (H_o - H_i) H_i}
	{\alpha \mu_s + i \omega (H_o - H_i) \beta}
	\\
	c_3 & = \frac{\left[\frac{i P }{\rho_f \omega} + \frac{i\omega P}{2\mu_s} (H_o^2 - H_i^2)\right]\beta
	+ PH_i \alpha}
	{\alpha \mu_s + i \omega (H_o - H_i) \beta}
	\\
	c_4 &= \frac{PH_o^2}{2\mu_s} -  H_o\frac{\left[\frac{i P }{\rho_f \omega} + \frac{i\omega P}{2\mu_s} (H_o^2 - H_i^2)\right]\beta
		+ PH_i \alpha}
	{\alpha \mu_s + i \omega (H_o - H_i) \beta}
\end{align}
where $\alpha =  e^{k_f H_i} + e^{-k_f H_i}$ and $\beta = \mu_f k_f \left(e^{k_f H_i} - e^{-k_f H_i}\right)$.

\textbf{Transient fluid and transient solid} ($\rho_f,~\rho_s > 0$)

In the transient fluid / transient solid case, the three equations obtained from the expansion
of~\eqref{eq:strong_form_linear_dynamic_constraint},
\eqref{eq:strong_form_wall_condition} and \eqref{eq:strong_form_linear_kinematic_constraint} are:
\begin{align}
	\left( e^{k_f H_i} + e^{-k_f H_i}\right) c_1 -
	i\omega \sin(k_s H_i) c_3 -
	i\omega \cos(k_s H_i) c_4
	& = - \left( \frac{P}{i\rho_f\omega} + \frac{iP}{\rho_s\omega} \right), \\
	\mu_f k_f \left(e^{k_f H_i} - e^{-k_f H_i}\right) c_1
	-\mu_s k_s \cos(k_s H_i)  c_3
	+\mu_s k_s \sin(k_s H_i)   c_ 4
	& = 0, \\
	\sin(k_s H_o)  c_3 + \cos(k_s H_o) c_4
	& = \frac{P}{\rho_s \omega^2}.
\end{align}
The resulting closed-form expressions are:
\begin{align}
	c_1 & = \frac{\mu_s k_s \cos \left[k_s(H_i - H_o)\right] \left(\frac{P}{i \rho_f\omega} + \frac{i P}{\rho_s \omega}\right)
	-i\omega\mu_s k_s  \frac{P}{\rho_s \omega^2}}
	{-\mu_s k_s \alpha \cos\left[k_s(H_i - H_o)\right]
		+ i\omega \beta \sin\left[k_s(H_i - H_o)\right]}, \\
	c_3 & = \frac{ \beta \cos(k_s H_o)   \left(\frac{P}{i \rho_f\omega} + \frac{i P}{\rho_s \omega}\right)
	- \left[i\omega\beta\cos(k_s H_i) + \alpha \mu_s k_s \sin(k_s H_i) \right] \frac{P}{\rho_s \omega^2}
	}{ -\mu_s k_s \alpha \cos\left[k_s(H_i - H_o)\right]
	+ i\omega \beta \sin\left[k_s(H_i - H_o)\right]}, \\
                    c_4 & =\frac{P}{\rho_s \omega^2} \sec(k_s H_o) - \tan(k_sH_o)  c_3
    .
\end{align}

\subsubsection{Nonliear fluid/solid in three dimensions}
\label{derivation-3D-sec}
In the three dimensional setting,
the channel is replaced with a tube. Similarly to the previous case,
we now consider the interaction between the pulsatile flow and the
hyperelastic wall which undergoes a shearing deformation along
the flow direction.
The domains of the two media in the reference configuration are shown
in Figure~\ref{domain-3D-fig}, where we used an analogous notation
for domains, boundaries and interfaces as in the case of the two-dimensional
FSI problem. Similarly, $H_i$ and $H_o$ now denote the inner and outer radii of the
tube. For ease, we use a cylindrical coordinate system defined by radial ($r$),
angular ($\theta$) and axial ($z$) positions.

The FSI problem takes the following general form:
\begin{align}
\rho_f \partial_t \boldsymbol{v}_f +  \rho_f \boldsymbol{v}_f \cdot \nabla_{\boldsymbol{x}}
 \boldsymbol{v}_f
+\nabla_{\boldsymbol{x}} \cdot \boldsymbol{\sigma}_f &= \boldsymbol{0} && \text{in } \Omega_f,
\label{eq:strong_form_nonlinear_fluid_momentum_balance}
\\
\nabla_{\boldsymbol{x}} \cdot \boldsymbol{v}_f &= 0 && \text{in } \Omega_f,
\\
[\boldsymbol{v}_f]_k & = 0 && \text{on } \Gamma_f^I \cup \Gamma_f^O,~\text{where } k\in\{r,\theta\},
\\
(\boldsymbol{\sigma}_f \cdot \boldsymbol{n}_f ) \cdot \boldsymbol{e}_z
& = [\boldsymbol{t}_f]_z && \text{on } \Gamma_f^I \cup \Gamma_f^O,
\label{eq:strong_form_nonlinear_fluid_Neumann_condition}
\\
\boldsymbol{v}_f (\cdot, 0) &= \boldsymbol{v}_f^0 && \text{in } \Omega_f (0),
\\
\hline
\rho_s \partial_{tt} \boldsymbol{u}_s -\nabla_{\boldsymbol{X}}
\cdot \boldsymbol{P}_s &= \boldsymbol{0} && \text{in } \Omega_s(0)\times[0,T],
\label{eq:strong_form_nonlinear_solid_momentum_balance}
\\
J-1 &= 0 && \text{in } \Omega_s(0) \times [0,T],
\\
\boldsymbol{u}_s (\cdot, t) &= \boldsymbol{0} && \text{on } \Gamma_s^W,
\\
[\boldsymbol{u}_s]_k & = 0 && \text{on }
\Gamma_s^I \cup \Gamma_s^O,~\text{where } k\in \{r, \theta \}
\\
(\boldsymbol{P}_s \cdot \boldsymbol{N}_s) \cdot \boldsymbol{e}_z
&= \boldsymbol{t}_s && \text{on }
\Gamma_s^I \cup \Gamma_s^O,
\\
\boldsymbol{u}_s (\cdot, 0) &= \boldsymbol{u}_s^0 && \text{in } \Omega_s (0),
\\
\boldsymbol{v}_s (\cdot, 0) & = \boldsymbol{v}_s^0 && \text{in } \Omega_s (0),
\\
\hline
\boldsymbol{\sigma}_f \cdot \boldsymbol{n}_f + \boldsymbol{P}_s \cdot \boldsymbol{N}_s &= \boldsymbol{0} && \text{on } \Gamma^\lambda, \label{eq:strong_form_nonlinear_dynamic_constraint}\\
\boldsymbol{v}_f - \boldsymbol{v}_s &= \boldsymbol{0} && \text{on } \Gamma^\lambda,
\label{eq:strong_form_nonlinear_kinematic_constraint}
\end{align}
where the fluid Cauchy stress tensor is the same as in the linear case.
For the solid, the first Piola-Kirchhoff stress tensor is defined as:
\begin{equation*}
	\boldsymbol{P}_s = \frac{\mu_s}{J^{2/3}}
	 \left[\boldsymbol{F} - \frac{\boldsymbol{F} : \boldsymbol{F}}{3} \boldsymbol{F}^{-T}\right] - p_s \boldsymbol{F}^{-T}.
\end{equation*}
The domains at time zero, $\Omega_s(0)$ and $\Omega_f(0)$, do not necessarily coincide with the reference
domains, $\Omega_s^0$ and $\Omega_f^0$, e.g., when $\boldsymbol{u}^0_s\neq\boldsymbol{0}$.
For now, we will consider
$\boldsymbol{u}_s^0$, $\boldsymbol{v}_s^0$, $\boldsymbol{t}_f$ and $\boldsymbol{t}_s$ to be known, since
their definitions are not used to arrive to the analytical solutions. 
\addtxt{Note that applying traction boundary conditions at the inlet of the fluid
	problem, as in Eq.~\ref{eq:strong_form_nonlinear_fluid_Neumann_condition}, will not necessarily
	lead to a stable formulation, see~\cite{hughes2005conservation,moghadam2011comparison}.
	That being said, from our experience, stabilization was not required for this class of problems 
	in order to get accurate solutions. 
	Alternatively, one may achieve a more robust problem formulation by changing
	the treatment of the boundary to Dirichlet or other stabilized forms. }

Similarly to the 2D problem in Section~\ref{derivation-2D-sec}, our derivation is based
on the assumption that the flow and deformation
field are axisymmetric and axially invariant,
i.e. $\boldsymbol{v}_f = v_f(r,t) \boldsymbol{e}_z$ and
$\boldsymbol{u}_s = u_s(r,t)\boldsymbol{e}_z$.
A first consequence of these is that the advective term in
Equation~\eqref{eq:strong_form_nonlinear_fluid_momentum_balance}
is always null. Secondly, we can keep using $\rho_f$
as a switch between the quasi-static and transient modes.
Thirdly, we can use the same coordinate vectors
(i.e. $\boldsymbol{e}_r$, $\boldsymbol{e}_{\theta}$ and $\boldsymbol{e}_z$),
to refer to both the reference and current configurations.
The same notation simplification is also valid for
 $r$ and $\theta$. Conversely, we use $Z$ and $z$
 to distinguish the axial coordinates in the reference and deformed configuration.

In order to derive the general analytical solution for the fluid, we first observe that:
\begin{equation*}
	\nabla_{\boldsymbol{x}} \boldsymbol{v}_f = \frac{\partial v_f(r,t)}{\partial r} \boldsymbol{e}_r \otimes \boldsymbol{e}_z.
\end{equation*}
Based on this, we can expand the fluid momentum equation~\eqref{eq:strong_form_nonlinear_fluid_momentum_balance} into its three vector components:
\begin{equation}
	\frac{\partial p_f}{\partial r} \boldsymbol{e}_r + \frac{1}{r} \frac{\partial p_f}{\partial \theta} \boldsymbol{e}_{\theta} +
	\left[\rho_f\partial_t {v}_f - \frac{\mu_f}{r}\frac{\partial}{\partial r}
	\left( r\frac{\partial v_f}{\partial r} \right) + \frac{\partial p_f}{\partial z}\right] \boldsymbol{e}_z
	= \boldsymbol{0}.
	\label{eq:momentum_balance_fluid_nonlin_expanded}
\end{equation}
From the radial and circumferential directions of the equation,
we can conclude that $p_f$ is constant in the plane perpendicular
to the flow direction.
Taking the partial derivative in the axial direction, we obtain that $\partial^2 p_f / \partial z^2 = 0$.
Thus, $p_f$ is either a constant or a linear function in $z$,
as observed in the 2D linear case.
If we assume the velocity and pressure fields are periodic
(with a single harmonic frequency) and that they are separable,
we arrive at the following forms:
\begin{equation}
	v_f = \Re\{\text{v}_f(r) e^{i\omega t} \} \hspace{2mm} \text{ and } \hspace{2mm}
	p_f = \Re\{P (L - z) e^{i\omega t}  \},
	\label{eq:3D_fluid_nonlinear_solution_periodic}
\end{equation}
where $\text{v}_f : [0, H_i] \rightarrow \mathbb{C}$ and $P_f \in \mathbb{C}$
are the velocity and pressure over domain length amplitudes.
 Returning to the axial component of
\eqref{eq:momentum_balance_fluid_nonlin_expanded} and simplifying the complex exponential term,
we obtain the following ODE:
\begin{equation*}
	\frac{1}{r}\frac{\partial}{\partial r} \left(r \frac{\partial \text{v}_f}{\partial r}\right) - k_f^2 \text{v}_f
	+ \frac{P}{\mu_f} = 0,
\end{equation*}
which, depending on the value of $\rho_s$, is satisfied by the following general solutions:
\begin{equation}
	\begin{aligned}
		\text{v}_f & = -\frac{P r^2}{4 \mu_f} + c_1 + c_2\ln(r), & (\rho_s = 0),\\
		\text{v}_f & = \frac{P }{\mu_f k_f^2} + c_1 J_0(i k_f r) + c_2 Y_0(i k_f r), & (\rho_s > 0).
	\end{aligned}
	\label{eq:fluid_vel_3D_nonlinear}
\end{equation}
Here, $J_0$ and $Y_0$ denote the first and second kind Bessel functions
of order 0. Finally, we can remark that $c_2 = 0$ in order to prevent
the singularity which we would otherwise have at $r = 0$.

Moving on to the solid, we begin by defining the deformation tensor as
$\boldsymbol{F} = \boldsymbol{I} + \frac{\partial u_s}{\partial r}
\boldsymbol{e}_r \otimes \boldsymbol{e}_z$. Thus, we can write the
first Piola-Kirchhoff into its vector components:
\begin{equation}
	\boldsymbol{P} = \mu_s\left[2\frac{\partial u_s}{\partial r}
	\text{sym}(\boldsymbol{e}_r \otimes \boldsymbol{e}_z) -
	\left(\frac{1}{3} \left[\frac{\partial u_s}{\partial r}\right]^2
	+ \frac{p_s}{\mu_s}  \right)
	\left(\boldsymbol{I} - \frac{\partial u_s}{\partial r}
	\boldsymbol{e}_z \otimes \boldsymbol{e}_r\right)
	\right].
	\label{eq:3D-nonlinear-solid-piola-kirchhoff}
\end{equation}
Substituting this into~\eqref{eq:strong_form_nonlinear_solid_momentum_balance},
we can expand the momentum balance equation as:
\begin{equation}
	\left[ - \frac{\partial}{\partial r}
	\left(\frac{\mu_s}{3}\left[\frac{\partial {u}_s}{\partial r}
	\right]^2 + {p}_s \right)
	+ \frac{\partial {u}_s}{\partial r} \frac{\partial {p}_s }{\partial z}
	\right] \boldsymbol{e}_r         	-\frac{1}{r}\frac{\partial {p}_s}{\partial \theta} \boldsymbol{e}_{\theta}   	+ \left[\frac{\mu_s}{r}\frac{\partial}{\partial r}
	\left(r \frac{\partial {u}_s}{\partial r}\right) - \frac{\partial {p}_s}{\partial z}
	-\rho_s \partial_{tt} {u}_s
	\right]  \boldsymbol{e}_z        	\label{eq:3D_nonlinear_solid_momentum_balance}
	= \boldsymbol{0}.
\end{equation}
Applying the axial partial derivative to the $\boldsymbol{e}_z$-component, we obtain that $\partial^2 {p}_s
/\partial z^2 = 0$. From this we conclude that ${p}_s$ is either a constant or linear with
respect to $z$. Furthermore, by applying the axial derivative to the $\boldsymbol{e}_r$-component, we find that
$\partial^2 {p}_s/\partial r\partial z = 0$, i.e. $\partial {p}_s/ \partial z$ is radially constant. Based on these results,
we can integrate the radial component to obtain the following formula for the pressure:
\begin{equation}
	{p}_s = c_I(t) + (Z + {u}_s) c_{II}(t)  - \frac{\mu_s}{3} \left[\frac{\partial u_s}{\partial r}\right]^2.
	\label{eq:3D_solid_nonlinear_pressure_general}
\end{equation}
Similarly to the 2D linear case, we assume that the solid displacement
and the $c_{II}$ constant are separable
periodic (with a single harmonic frequency):
\begin{align*}
u_s  & = \Re\{\text{u}_s e^{i\omega t} \} \quad \text{and} \quad c_{II} = \Re\{c_3 e^{i\omega t} \}
\end{align*}
with $\text{u}_s:[H_i, H_o]\rightarrow \mathbb{C}$ and $c_3 \in \mathbb{C}$. Consequently, we can
simplify the temporal complex exponentials in Equation~\eqref{eq:3D_nonlinear_solid_momentum_balance}, which yields:
\begin{equation}
\frac{1}{r}\frac{d}{d r}
\left(r \frac{d \text{u}_s}{d r}\right)
+ k_s^2 \text{u}_s - \frac{c_3}{\mu_s} = 0.
\label{eq:3D-nonlin-bessel-pde}
\end{equation}
For the transient case, i.e.  $\rho_s > 0$ and $k_s \neq  0$, we can recognize this as a
non-homogeneous Bessel's ODE, with the order of the Bessel function
equal to zero. Hence, the general solid displacement solutions are:
\begin{equation}
	\begin{aligned}
		\text{u}_s & = \frac{c_3 r^2}{4\mu_s} + c_4 \ln(r) + c_5 , &&(\rho_s = 0) \\
		\text{u}_s & = \frac{c_3}{\mu_s k_s^2} + c_4 J_0(-k_s r) + c_5Y_0 (-k_s r), &&(\rho_s > 0).
	\end{aligned}
	\nonumber
\end{equation}
Furthermore, by including the fixed wall boundary condition, i.e.  $\text{u}_s(H_o) = 0$,
we arrive at the following simplified forms:
\begin{equation}
\begin{aligned}
\text{u}_s & = \frac{c_3}{4\mu_s} (r^2 - H_o^2) + c_4 \ln\left(\frac{r}{H_o}\right) , &&(\rho_s = 0) \\
\text{u}_s & = \frac{c_3}{\mu_s k_s^2} \left[ 1 - \frac{Y_0(-k_s r)}{Y_{0,s}^r} \right]
+ c_4 \left[J_0(-k_s r) -\gamma Y_{0}(-k_s r) \right], &&(\rho_s > 0).
\end{aligned}
\label{eq:3D_nonlinear_disp_sol}
\end{equation}

The remaining unknown constants (i.e. $c_1$, $c_3$, $c_4$ and $c_I$)
can only be identified
by verifying that the kinematic and traction boundary conditions
are satisfied. In the following, we look at the different
transient / quasi-static permutations and derive for each of them the
appropriate closed formulations for the set of constants.

\textbf{Quasi-static fluid and quasi-static solid} ($\rho_f = \rho_s = 0$)

For the fluid and solid problems to be coupled, the traction  and
kinematic conditions in equations~\eqref{eq:strong_form_nonlinear_dynamic_constraint}
and~\eqref{eq:strong_form_nonlinear_kinematic_constraint}, respectively, need to hold.
The former can now be expanded as follows:
\begin{equation}
	 \left[ \frac{\mu_s}{3} \left.\left(\frac{\partial u_s}{\partial r}\right)^2\right|_{r = H_i} + p_s(H_i,z,t) -p_f(z,t)  \right] \boldsymbol{e}_r
	 +\left[\mu_f\left.\frac{\partial v_f}{\partial r}\right|_{r = H_i} - \mu_s \left.\frac{\partial u_s}{\partial r}\right|_{r = H_i}\right]
	 \boldsymbol{e}_z
	= \boldsymbol{0},
\end{equation}
for all $z\in[0,L]$ and $t\in[0,T]$.
Based on the solid and pressure formulations in~\eqref{eq:3D_fluid_nonlinear_solution_periodic}
and~\eqref{eq:3D_solid_nonlinear_pressure_general}, the radial component can be further reduced to:
\begin{equation}
	c_I(t) + [Z+ u_s(H_i,t)]c_3 e^{i\omega t} - P (L - Z) e^{i\omega t} =0.
	\label{eq:3D_traction_balance_shear}
\end{equation}
From this, we can deduce that $c_3 = -P$ and $c_I(t) = [L  + u_s(H_i,t)]\Re\{Pe^{i\omega t}\}$, two
results which are valid for all transient/quasi-static permutation.
Consequently, \eqref{eq:3D_solid_nonlinear_pressure_general}
can be expanded using \eqref{eq:3D_nonlinear_disp_sol} in quasi-static form to yield the
following pressure formula:
\begin{align}
	p_s(r,Z,t) =&~\left[L - Z + u_s(H_i,t) - u_s(r,t)\right]\Re\{Pe^{i\omega t}\}  - \frac{\mu_s}{3} \left[\frac{\partial u_s(r,t)}{\partial r}\right]^2 \nonumber \\
	    =&~\left[L - Z + \Re\left\{c_4 e^{i\omega t}  \ln \left(\frac{H_i}{r}\right)\right\}  +
	        \Re\left\{\frac{P}{4\mu_s}(r^2 - H_i^2)e^{i\omega t} \right\} \right]\Re\{P e^{i\omega t}\} \nonumber \\
	    &- \frac{\mu_s}{3} \left[\Re\left\{\frac{c_4 e^{i\omega t}}{r} - \frac{P r e^{i\omega t}}{2\mu_s}  \right\}\right]^2.
	\label{eq:3D_solid_nonlinear_pressure_general2}
\end{align}
In the case of shear component of the traction coupling
condition we obtain that:
\begin{equation}
	- \mu_s \frac{P H_i}{2\mu_s} + \mu_s \frac{P H_i}{2 \mu_s} + \frac{c_4}{H_i} = 0, \nonumber
\end{equation}
which is equivalent to $c_4 = 0$. Using the kinematic coupling condition, we can identify the
last integration constant:
\begin{equation}
	c_1 = \frac{PH_i^2}{4\mu_f} + i\omega \frac{P(H_o^2 - H_i^2)}{4\mu_s}.
\end{equation}

\textbf{Transient fluid and quasi-static solid} ($\rho_f > 0,~\rho_s = 0$)

Similarly to the previous example, the two integration constants can be found by solving the system of two equations formed by the shear components of the traction and the  kinematic coupling conditions:
\begin{align}
	\left[\begin{array}{c c}
	i \mu_f J^*_{1,f} &  \mu_s / H_i \\
	J^*_{0,f}  & i\omega \ln\left(H_o / H_i\right)
	\end{array}
	\right]
	\left[\begin{array}{c}
	c_1  \\ c_4
	\end{array}
	\right]
	= \left[\begin{array}{c}
	P H_i /2 \\
	i \omega P(H_o^2 - H_i^2) / (4\mu_s) - P/(\mu_f k_f^2)
	\end{array}\right],
\end{align}
where $J^*_{0,f} = J_0(ik_f H_i)$ and $J^*_{1,f} = k_f J_1(ik_f H_i)$.
The resulting closed-formulation of the two constants are:
\begin{align}
					c_1 & = -iP \frac{\omega H_i^2 \ln (H_i / H_o) / 2 + \omega (H_o^2 -H_i^2) / 4  + i \mu_s/(\mu_f k_f^2)
	}
	{\mu_f \omega H_i \ln(H_i / H_o) J^*_{1,f}- \mu_s J^*_{0,f}},  \\
					c_4 & = - P \frac{\mu_s J^*_{0,f} H_i /2  + \mu_f \omega J^*_{1,f} (R^2_o - R^2_i)/4 + i\mu_s J^*_{1,f}/ k_f^2
	}{\mu_f\mu_s\omega \ln(H_i / H_o ) J^*_{1,f} - \mu_s^2 J^*_{0,f} / H_i}.
	\label{eq:c4_3D_nonlin_tf_qs}
\end{align}
The general pressure solution in Equation~\eqref{eq:3D_solid_nonlinear_pressure_general2} continues to apply in this case as
well (with $c_3 = -P$), provided one uses the definition of $c_4$ found in \eqref{eq:c4_3D_nonlin_tf_qs} rather than $c_4 = 0$.

\textbf{Quasi-static fluid and transient solid} ($\rho_f = 0,~ \rho_s > 0$)

The closed-formulation for the $c_4$ constant is obtained by rearranging the kinematic condition
in~\eqref{eq:strong_form_nonlinear_kinematic_constraint}:
\begin{equation}
    c_4 =- P \frac{2 Y_{1,s}^* + i H_i k_s^2 Y_{0,s}^r
    }{2 \mu_s k_s^2 Y_{0,s}^r \Delta_1},
\end{equation}
where $J_{0,s}^* = J_0(- k_s H_i)$, $J_{1,s}^* = i k_s J_1(- k_s H_i)$,
$Y_{0,s}^* =  Y_0(- k_s H_i)$, $Y_{1,s}^* = i k_s Y_1(- k_s H_i)$,
$\Delta_0 = J_{0,s}^* - \gamma Y_{0,s}^*$ and
$\Delta_1 = J_{1,s}^* - \gamma Y_{1,s}^*$.
Using this result and expanding
the axial component of the traction coupling condition
in~\eqref{eq:strong_form_nonlinear_solid_momentum_balance},
we can also derive the closed-formulation of the last constant:
\begin{equation}
    c_1 = \frac{P H_i^2}{4 \mu_f}
    - \frac{{i\omega P}}{\mu_s k_s^2} \left(1 - \nu_0\right)
    - P \omega \frac{2 i Y_{1,s}^* - H_i k_s^2 Y_{0,s}^r}{2\mu_s k_s^2 Y_{0,s}^r}
    \frac{\Delta_0}{\Delta_1},
\end{equation}
where $\nu_k = Y_{k,s}^* / Y_{0,s}^r$ and $k\in\{0,1\}$.

To find the solution for the solid pressure, we expand the more
general result in Equation~\eqref{eq:3D_solid_nonlinear_pressure_general}
using the results of Equation~\eqref{eq:3D_traction_balance_shear} and the corollaries,
$c_3 = -P$ and $c'(t) = [L  + u_s(H_i,t)]\Re\{Pe^{i\omega t}\}$,
plus the transient version of the solid displacement formula in Equation~\ref{eq:3D_nonlinear_disp_sol}:
\begin{align}
 	p_s (r,Z,t) =&~\left[L - Z + u_s(H_i,t) - u_s(r,t)\right]\Re\{Pe^{i\omega t}\}  - \frac{\mu_s}{3} \left[\frac{\partial u_s(r,t)}{\partial r}\right]^2 \nonumber \\
 	=&~\left[ L - Z - \Re \left\{\frac{P}{\mu_s k_s^2} \left(\frac{Y_0(-k_s r)}{Y_{0,s}^r} - \nu_0 \right) e^{i\omega t} \right\} \right] \Re\{Pe^{i\omega t}\} \nonumber \\
 	&+ \Re \{c_4 \left[\Delta_0 - J_0(-k_s r) + \gamma Y_0(-k_s r)\right] e^{i\omega t}\}
 	    \Re\{P e^{i\omega t}\} \nonumber \\
 	&- \frac{\mu_s}{3} \left[ \Re \left\{\frac{P}{\mu_s k_s} \frac{Y_1(-k_s r)}{Y_{0,s}^r} e^{i\omega t} +
 	    c_4 k_s (J_1(-k_s r) - \gamma Y_1(-k_s r)) e^{i \omega t} \right\} \right]^2.
 	\label{eq:3D_nonlinear_solid_pressure_transient}
\end{align}

\textbf{Transient fluid and transient solid} ($\rho_f, \rho_s > 0$)

In this case, the system of equations resulting from the axial components of the traction and kinematic coupling conditions can be written as:
\begin{equation}
    \left[
    \begin{array}{c c}
        -\mu_fJ_{1,f}^* & \mu_s \Delta_1 \\
        -J_{0,f}^* & -i\omega\Delta_0
    \end{array}
    \right]
    \left[
    \begin{array}{c}
    c_1 \\ c_4
    \end{array}
    \right]
    =
    \left[
    \begin{array}{c}
    - P\nu_1/k_s^2 \\
    - i \omega P(1-\nu_0) / (\mu_sk_s^2) + P/(\mu_sk_f^2)
    \end{array}
    \right].
\end{equation}
Solving the system gives the following closed forms:
\begin{align}
	c_1 & = - \frac{i\omega\Delta_0\nu_1 \frac{P}{k_s^2} +
		\mu_s \Delta_1 \left[(1-\nu_0)\frac{i\omega P}{\mu_s k_s^2} + \frac{P}{\mu_f k_f^2}	\right]
	}{\mu_s J_{0,f}^* \Delta_1 - i\omega\mu_f J_{1,f}^*\Delta_0},\\
	c_4 & = - \frac{J_{0,f}^*\nu_1 \frac{P}{k_s^2} +
		\mu_f J_{1,f}^* \left[(1-\nu_0) \frac{i\omega P}{\mu_s k_s^2} + \frac{P}{\mu_f k_f^2}	\right]
	}{\mu_s J_{0,f}^* \Delta_1 - i\omega\mu_f J_{1,f}^*\Delta_0}.
	\label{eq:c4_3D_nonlin_tf_ts}
\end{align}
The analytical solution for the solid pressure retains the form presented in~\eqref{eq:3D_nonlinear_solid_pressure_transient}
with $c_3 = -P$, but using the version of $c_4$ in \eqref{eq:c4_3D_nonlin_tf_ts}.
\section{Numerical fluid-structure interaction implementation}
\label{fsi-implementation-sec}
As a first use of the benchmark and as guide to future applications and results reporting,
we also carried out a validation test of our finite element implementation
\cite{LeeEtAl2016,HessenthalerRoehrleNordsletten2017}.
This section presents a brief review of the structure of the weak form and discretization scheme
employed to solve different model permutations (presented later
in Section \ref{num-experiments}).
\subsection{Finite element weak form}
\label{cheart-implementation-sec}
The fluid-structure coupling strategy used in this work is based on the technique
described in~\cite{NordslettenKaySmith2010,
HessenthalerRoehrleNordsletten2017,
HessenthalerFriedhoffRoehrleNordsletten2016}.
 In this approach, a Lagrange multiplier
variable, $\boldsymbol{\lambda} = \boldsymbol{t}_f = -\boldsymbol{t}_s$,
is introduced such that it weakly enforces both the kinematic and dynamic
constraints, resulting in an FSI system which is solved monolithically.
This strategy allows us to choose the domain discretization and polynomial
interpolation schemes which are more appropriate for the representation of either
fluid and solid solutions. In the following, we outline how this applies to our case.

The general discrete weak form including the FSI coupling conditions
can be written as follows:
Find
$\boldsymbol{s}^n :=\left(\boldsymbol{v}^{n}_f,
\boldsymbol{v}^n_s,
\boldsymbol{\lambda}^n, p^n_f, p^n_s\right)
\in
\boldsymbol{\mathcal{S}}^h_D :=
\boldsymbol{\mathcal{V}}^h_0 \times
\boldsymbol{\mathcal{U}}^h_0 \times
\boldsymbol{\mathcal{M}}^h_D \times
\mathcal{W}^h_f \times
\mathcal{W}^h_s$,
such that for every
$\boldsymbol{d} :=
\left(\boldsymbol{w}_f, \boldsymbol{w}_s,
\boldsymbol{q}_\lambda, q_f, q_s\right)
\in
\boldsymbol{\mathcal{S}}^h_0 :=
\boldsymbol{\mathcal{V}}^h_0 \times
\boldsymbol{\mathcal{U}}^h_0 \times
\boldsymbol{\mathcal{M}}^h_0 \times
\mathcal{W}^h_f \times
\mathcal{W}^h_s$:
\begin{align}
	R(\boldsymbol{s}^{n}, \boldsymbol{s}^{n-1}, \boldsymbol{d}) :=
	 & \int_{\Omega^n_{f,h}} \rho_f \left[ \frac{\mathbf{v}^n_f - \mathbf{v}^{n-1}_f}{\Delta^n_t}
	+\delta_{nl}\left(\boldsymbol{v}^n_f - \hat{\boldsymbol{v}}^n_f\right)
	\cdot \nabla_{\boldsymbol{x}} \boldsymbol{v}^n_f\right] \cdot \boldsymbol{w}_f~d\Omega
	\nonumber \\
	+ & \int_{\Omega^n_{f,h}}  \boldsymbol{\sigma}_f^n : \nabla_{\boldsymbol{x}} \boldsymbol{w}_f
	+ q_f \nabla_{\boldsymbol{x}} \cdot \boldsymbol{v}_f^n~d\Omega
	\nonumber \\
	- & \int_{\Gamma^{I,n}_{f,h} \cup \Gamma^{O,n}_{f,h}}
	\boldsymbol{t}^n_{f} \cdot \boldsymbol{w}_f~d\Omega
	\nonumber \\
	+ & \int_{\Omega^0_{s,h}} \delta_{nl}
	\left[ J^n_s \rho_s \frac{\boldsymbol{v}^n_s - \boldsymbol{v}^{n-1}_s}{\Delta^n_t}
	+ \boldsymbol{P}^n_s:\nabla_{\boldsymbol{X}} \boldsymbol{w}_s
	+ q_s (J^n_s - 1) \right] ~d\Omega
	\nonumber\\
	+ & \int_{\Omega^0_{s,h}}(1 - \delta_{nl} )
	\left[ \rho_s \frac{\boldsymbol{v}^n_s-\boldsymbol{v}^{n-1}_s}{\Delta^n_t}
	+ \boldsymbol{\sigma}_s^n : \nabla_{\boldsymbol{X}} \boldsymbol{w}_s
	+ q_s \nabla_{\boldsymbol{X}} \cdot \boldsymbol{v}^n_s \right]~d\Omega
	\nonumber \\
	- & \int_{\Gamma^{I,n}_{s,h} \cup \Gamma^{O,n}_{s,h}} \boldsymbol{t}^n_s \cdot \boldsymbol{w}_s~d\Omega
	\nonumber \\
	+ & \int_{\Omega^0_{\lambda,h}} \boldsymbol{\lambda}^n
	\cdot \left(\boldsymbol{w}_f - \boldsymbol{w}_s \right)
	+ \boldsymbol{q}_{\lambda}
	\cdot \left(\boldsymbol{v}^n_f - \boldsymbol{v}^n_s \right)~d\Omega= 0,
\end{align}
where the $\delta_{nl}$ term acts as a switch between the linear and nonlinear cases and
the current displacement is defined as
$\boldsymbol{u}^n_s = \boldsymbol{u}^{n-1}_s + \Delta_t^n \boldsymbol{v}^n_s$
with the length of the current time step $\Delta_t^n$.
Further, $\boldsymbol{t}^n_f$ and $\boldsymbol{t}^n_s$
are the discrete representations
of the analytical solutions of the surface tractions acting on the inflow and outflow
surfaces of the fluid and solid, respectively.

In the nonlinear cases, where the ALE term is active,
the weak form is expanded to include the fluid domain motion problem,
which takes the form: Find $\hat{\boldsymbol{v}}_f^n \in \boldsymbol{\mathcal{W}}^h_D$
such that for any $\boldsymbol{z}\in \boldsymbol{\mathcal{W}}^h_0$:
\begin{align}
	R_{g} (\hat{\boldsymbol{v}}^n_f, \hat{\boldsymbol{v}}^{n-1}_f, \boldsymbol{z}) :=
	\int_{\Omega^0_{f,h}}	\delta_{nl} \left[
	\frac{\hat{\boldsymbol{v}}^n_f - \hat{\boldsymbol{v}}^{n-1}_f}{\Delta^n_t}
	\cdot \boldsymbol{z} - \boldsymbol{\Phi} \nabla_{\boldsymbol{X}}\hat{\boldsymbol{v}}^n_f
	: \nabla_{\boldsymbol{X}} \boldsymbol{z}
	\right] ~d\Omega	= 0.
\end{align}
The newly introduced field, $\hat{\boldsymbol{v}}_f^n$, represents the
artificial grid velocity and is used to account for the motion of the
domain in the ALE formulation. The domain deformation is preferentially
weighted in the axial direction by means of the diffusion field $\boldsymbol{\Phi}$
in order to avoid radial variations.

The definitions of the function spaces are:
\begin{equation}
	S^k(\Omega^0_{i,h})  = \{f:\Omega^0_{i,h} \rightarrow \mathbb{R}~|~
	f\in\mathcal{C}^0(\bar{\Omega}^0_{i,h}),~
	f|_{\tau_e} \in \mathbb{P}^k(\tau_e),~\forall~\tau_e\in\mathcal{T}^h_i \},
\end{equation}
which represent the general continuous $k^{th}$-order piecewise polynomial spaces
defined on $\Omega^0_{i,h}$. Consequently, we can define:
\begin{align*}
	& \boldsymbol{\mathcal{V}}^h = \left[S^2(\Omega^0_{f,h})\right]^d, &&
	\boldsymbol{\mathcal{U}}^h = \left[S^2(\Omega^0_{s,h})\right]^d, &&
	\boldsymbol{\mathcal{M}}^h = \left[S^2(\Gamma^{\lambda}_{f,h}) \right]^d,
	& \mathcal{W}^h_f = S^1(\Omega^0_{f,h}),  &&
	\mathcal{W}^h_s = S^1(\Omega^0_{s,h}).
\end{align*}
In the case of solid and fluid domain velocity, further restriction are applied on their
respective spaces in order to take incorporate the Dirichlet and homogeneous
boundary conditions:
\begin{align}
	\boldsymbol{\mathcal{V}}^h_0 & = \{\boldsymbol{v} \in \boldsymbol{\mathcal{V}}^h
	~|~ \boldsymbol{v}_{\bot} = \boldsymbol{0} \text{ on } \Gamma^{I,n}_{f,h}
	\cup \Gamma^{O,n}_{f,h} \}, \\
	\boldsymbol{\mathcal{U}}^h_0 & = \{\boldsymbol{v} \in \boldsymbol{\mathcal{U}}^h
	~|~ \boldsymbol{v}_{\bot} = \boldsymbol{0} \text{ on }
	 \Gamma^{I}_{s,h} \cup \Gamma^{O}_{s,h}
	 \text{ and } \boldsymbol{v} = \boldsymbol{0} \text{ on } \Gamma^{W}_{s,h}
	\}, \\
	\boldsymbol{\mathcal{W}}^h_D & = \{\boldsymbol{v} \in \boldsymbol{\mathcal{V}}^h
	~|~ \boldsymbol{v} = \boldsymbol{v}_f \text{ on } \Gamma^{\lambda}_h \}, \\
	\boldsymbol{\mathcal{W}}^h_0 & = \{\boldsymbol{v} \in \boldsymbol{\mathcal{V}}^h
	~|~ \boldsymbol{v} = \boldsymbol{0} \text{ on } \Gamma^{\lambda}_h \},
	\\
	\boldsymbol{\mathcal{M}}^h_D & = \{\boldsymbol{\lambda} \in
	\boldsymbol{\mathcal{M}}^h~|~
	\boldsymbol{\lambda}_{\bot} = \boldsymbol{t}^n_{h,\bot}
	\text{ on } \left( \Gamma^{\lambda}_h\cap\Gamma^I_{f,h} \right)
	\cup \left( \Gamma^{\lambda}_h\cap\Gamma^O_{f,h} \right) \} ,
	\\
	\boldsymbol{\mathcal{M}}^h_0 & = \{\boldsymbol{\lambda} \in
	\boldsymbol{\mathcal{M}}^h~|~
	\boldsymbol{\lambda}_{\bot} = \boldsymbol{0}
	\text{ on } \left( \Gamma^{\lambda}_h\cap\Gamma^I_{f,h} \right)
	\cup \left( \Gamma^{\lambda}_h\cap\Gamma^O_{f,h} \right) \} .
\end{align}
Here, we used the $\bot$ symbol to indicate more generally the components of
$\boldsymbol{v}$ and $\boldsymbol{\lambda}$ which are perpendicular to the flow direction, such as
$\hat{\boldsymbol{e}}_y$ in two dimensions, as well as $\hat{\boldsymbol{e}}_x$ and $\hat{\boldsymbol{e}}_y$
in the three-dimensional case.

\section{Numerical experiments}
\label{num-experiments}
\subsection{Numerical validation of analytic solutions}
\label{validation-of-derivations-seq}
The 16 analytic solutions presented in Section~\ref{derivations-seq}
and in the \emph{Supplementary Materials SM\pageref{online-supplement-sec}}
were validated numerically.
In the linear two dimensional case (Section~\ref{derivation-2D-sec}),
numerical differentiation was employed to verify that the analytic solutions
satisfy the (strong form) momentum and mass balance
(i.e.\
Equation~\ref{eq:strong_form_linear_fluid_momentum_balance}~/~\ref{eq:strong_form_linear_fluid_mass_balance}
and
Equation~\ref{eq:strong_form_linear_solid_momentum_balance}~/~\ref{eq:strong_form_linear_solid_mass_balance}),
the boundary conditions,
and the dynamic and kinematic coupling constraints
(i.e.\ Equation~\ref{eq:strong_form_linear_dynamic_constraint}~/~\ref{eq:strong_form_linear_kinematic_constraint}).
All other analytic solutions were validated similarly.
\subsection{Numerical experiments: Space-time discretization}
In Section \ref{results-sec}, transient and linear fluid / solid models are considered
in two dimensions, and transient and nonlinear fluid / solid models are considered in three dimensions.
In what follows, details about the space-time discretization of the computational domains are discussed.
We further define the domain dimensions that are complemented with material parameters
in Section~\ref{results-2D-analytic-sec} and \ref{results-3D-analytic-sec}.
Parameters were selected to highlight key solution features
and for conducting a space-time convergence analysis.
But they are arbitrary in the sense that they were chosen without a particular application area in mind.
Thus, they are given in a unitless format.
In contrast, Section~\ref{results-physiological-parameters-sec} considers physiologically relevant
parameters, which are given in appropriate units.
\subsubsection{Linear FSI in two dimensions}
For the linear FSI case in two dimensions, the channel length and width
are set to $L = H_i = 1$ and the wall thickness is set to $H_o - H_i = 0.2$,
see Figure \ref{domain-2D-fig}.
The temporal cycle length is selected as $T = 1.024$ and the length of the time domain is set as $10 T$.

The fluid and solid domains are discretized using quadrilateral elements
with quadratic (and linear) interpolation for the velocity / displacement (and pressure) variables
and the Lagrange multiplier domain is discretized using line elements with quadratic interpolation.
Three different spatial refinement levels are considered, and referred to as \textit{coarse},
\textit{medium} and \textit{fine}, see Table~\ref{space-time-discretization-tab}.
The temporal domain is discretized using equidistant time points and seven different time step sizes,
such that we may consider temporal refinement for a given spatial refinement level
as well as three space-time refinement levels with $\delta_t / \delta_x^3 = \text{const}$.
\subsubsection{Nonlinear FSI in three dimensions}
In the nonlinear FSI case in three dimensions, the fluid domain length and radius
are set to $L = 1$ and $H_i = 0.7$ and the wall thickness is set to $H_o - H_i = 0.3$,
see Figure \ref{domain-3D-fig}.
To reduce the computational cost,
we only consider one quarter of the domain in numerical experiments (i.e.\ for $x, y \geq 0$)
and set the length of the time domain to $7 T$ with cycle length $T = 1.024$.

The fluid domain is discretized using tetrahedral elements (quadratic-linear interpolation
for velocity and pressure),
the solid domain using hexahedral elements (quadratic-linear interpolation for velocity / displacement
and pressure) and the Lagrange multiplier domain using triangle elements (quadratic interpolation).

Similar to the linear FSI case in two dimensions, three different refinement levels are considered,
see Figure~\ref{mesh-refinements-fig}.
Since the analytic solution for the solid pressure variable
(see Equation \eqref{eq:3D_solid_nonlinear_pressure_general}) is now a function of the displacement
(in contrast to the linear case), the mesh for the solid domain is more refined in the radial direction
than in the axial and circumferential directions.

The coarse solid mesh has $5$ elements in the circumferential direction, $8$ elements in the radial
direction and $8$ elements in the axial direction.
The coarse fluid mesh is selected, such that four tetrahedral element faces conform with one
hexahedral element face at the interface boundary.
Tetrahedral element sizes are approximately constant throughout the fluid domain.
The corresponding triangular mesh for the interface domain is embedded in the interface boundary
of the fluid mesh and each triangle element conforms with a tetrahedral mesh face of the fluid mesh.

The medium and fine refinement levels are achieved by using a uniform refinement of the solid domain
and an approximately uniform refinement of the fluid domain.
For all refinement levels, we note that all solid, interface and fluid interface nodes
(including midpoint nodes for quadratic elements) lay exactly on circles with radius $H_i$
around the $z$-axis.
We select time step sizes (see Table~\ref{space-time-discretization-tab}), such that we may consider
temporal refinement for a given spatial mesh, as well as three space-time refinement levels
with $\delta_t / \delta_x^3 = \text{const}$.
\subsection{Numerical solution}
The linear and nonlinear transient FSI cases are used to validate the {CH}eart \cite{LeeEtAl2016}
implementation of the FSI method detailed in Section \ref{cheart-implementation-sec}.
{CH}eart is based on the matrix solver {MUMPS} \cite{AmestoyDuffLexcellent2000}.
A Newton-Raphson-Shimanskii solver \cite{Shamanskii1967}
is employed to reduce the computational cost by re-using the Jacobian matrix (and its inverse)
as long as the residual norm could be reduced by a factor of $3/4$
\cite{HessenthalerRoehrleNordsletten2017}.

Complex expression evaluators were implemented to set initial and boundary conditions,
according to the settings
in Equation~\eqref{eq:strong_form_linear_fluid_momentum_balance}
- \eqref{eq:strong_form_linear_kinematic_constraint} and
Equation~\ref{eq:strong_form_nonlinear_fluid_momentum_balance}
- \eqref{eq:strong_form_nonlinear_kinematic_constraint}, respectively.
We further evaluated the analytic FSI solutions in {CH}eart after each time step
to compute the space-time error on-the-fly.
\section{Results}\label{results-sec}
This section considers the transient fluid and transient linear solid model in two spatial dimensions
and the transient fluid and transient nonlinear solid model in three spatial dimensions.
In the following sections,
we present the velocity, displacement and pressure solutions
for a unitless parameter set that is suitable for space-time convergence
analysis (Section~\ref{results-2D-analytic-sec} and \ref{results-3D-analytic-sec})
and for a physiologically relevant parameter set (Section~\ref{results-physiological-parameters-sec}).
In Section~\ref{numerical-solution-t-conv-cheart-sec}, we demonstrate the use of the analytic solutions for validating a non-conforming monolithic
FSI method (see Section~\ref{cheart-implementation-sec}).
\subsection{Transient fluid and transient linear solid in 2D}
\label{results-2D-analytic-sec}

To illustrate the linear FSI case in two spatial dimensions,
we consider the following material parameters:
the fluid density and viscosity are selected as $\rho_f = 1$ and $\mu_f = 0.01$,
and the solid density and stiffness are chosen as $\rho_s = 1$ and $\mu_s = 0.1$.
We further set the pressure amplitude as $P = 1$.
Here, the Womersley number is
$W = H_i \sqrt{\omega \rho_f / \mu_f}
\approx 24.77$,
and the Reynolds number is
$Re = 2 \rho_f V_f H_i / \mu_f
\approx 88.32$
(with $V_f = \max_{y,t} | v_f (y,t) | \approx 0.4416$).

Figure~\ref{2D-linear-tf-ts-solution-along-y-fig} illustrates the temporal variation of the
velocity, displacement and pressure solution profiles along the $y$-axis,
see Figure~\ref{domain-2D-fig}.
It further highlights the temporal variation of the fluid and solid pressure solutions
at the inlet ($x = 0$), the mid-way point in the $x$-direction ($x = L/2$)
and the outlet ($x = L$).
In the linear case, the fluid and solid pressure solution is linear along the $x$-axis,
continuous at the coupling boundary and constant along the $y$-axis.
Figure~\ref{2D-linear-tf-ts-solution-along-y-fig} further highlights that
$p_f (x,t) = p_s (x,t) = P (L - x) cos(\omega t)$.

The velocity and displacement solutions are constant along the $x$-axis, but varying in $y$.
The fluid velocity solution resembles a typical Womersley profile over time
and the maximum fluid velocity occurs at the interface, where it matches the deformation rate
of the solid due to the kinematic coupling condition.
The solid deformation rate at the wall (i.e.\ at $y = H_o$) is zero,
similar to the displacement, as enforced by the Dirichlet boundary condition.
The maximum absolute displacement $U_s$ is related to the maximum absolute
deformation rate $V_s$ by approximately $U_s \approx V_s / (2 \pi)$.
\subsection{Transient fluid and transient nonlinear solid in 3D}
\label{results-3D-analytic-sec}
To illustrate the transient nonlinear FSI solution in three dimensions,
the fluid density and viscosity were selected as $\rho_f = 2.1$ and $\mu_f = 0.03$,
and the solid density and stiffness were chosen as $\rho_s = 1$ and $\mu_s = 0.1$.
We further set the pressure amplitude as $P = 1$.
Compared to Section~\ref{results-2D-analytic-sec}, the Womersley number, $W \approx 14.51$,
and the Reynolds number, $Re \approx 8.036$
(with $V_f = \max_{r,t} | v_f (r,t) | \approx 0.082$) are smaller.

Figure~\ref{3D-nonlinear-tf-ts-solution-along-y-fig} illustrates the temporal variation
of the velocity, displacement and pressure profiles along the $y$-axis,
see Figure~\ref{domain-3D-fig}.
Compared to the transient linear FSI case in the previous section,
the peak solid deformation rate is of similar magnitude,
however, the peak fluid velocity is smaller.
Furthermore, the peak values are now present in the interior of the fluid and solid domains.
While the fluid velocity gradient near the wall is smaller compared to the linear FSI case,
the solid deformation rate and displacement is varying along the $y$-axis in a similar fashion.

Figure~\ref{3D-nonlinear-tf-ts-solution-along-t-fig} further highlights the temporal variation
of the fluid and solid pressure solutions at the inlet ($x = 0$),
the mid-way point in the $x$-direction ($x = L/2$) and the outlet ($x = L$).
The fluid and solid pressure solutions are no longer continuous at the interface,
which becomes more clear towards the outlet,
see Figure~\ref{3D-nonlinear-tf-ts-solution-along-t-fig} (bottom)
and Figure~\ref{3D-nonlinear-tf-ts-solution-along-y-fig} (bottom-right).
The fact that the solid pressure solution is now a function of
$\left[ \partial u_s / \partial r \right]^2$
(see Equation~\eqref{eq:3D_solid_nonlinear_pressure_general}) results in a higher temporal frequency
toward the outlet at $z = L$ as compared to the inlet at $z = 0$,
and in large pressure gradients in the radial direction,
see Figure~\ref{3D-nonlinear-tf-ts-solution-along-t-fig}.

Lastly, we emphasize that the response of the fluid and / or solid can be changed significantly
by modifying even just one parameter. As an exemplifying case, we modify the solid density
and set $\rho_s = 5$.
This does not significantly affect the flow (fluid velocity and pressure); however, changes
the solid response, see Figure~\ref{3D-nonlinear-tf-ts-solution-along-y-rhos5-fig}.
Similar to the mode shapes of an Euler-Bernoulli beam \cite{BauchauCraig2009},
the (transient) nonlinear solid model exhibits deformation modes
corresponding to the natural frequencies of the solid material.
\subsection{Transient fluid and transient nonlinear solid in 3D: Physiologically relevant parameters}
\label{results-physiological-parameters-sec}

While the parameter choices in Section~\ref{results-2D-analytic-sec}
and Section~\ref{results-3D-analytic-sec}
were motivated from the viewpoint of space-time convergence analysis
(see Section~\ref{numerical-solution-t-conv-cheart-sec}),
one may select a more physiologically relevant parameter set.
For example, consider parameters that are similar to those found in models of the ascending
aorta (parameters taken from \cite{BalmusMassingHoffmanRazaviNordsletten2019_preprint,
LiKamelRandoAndersonKumbasarLimaBluemke2003,ReymondCrosettoDeparisQuarteroniStergiopulos2012}),
with a fluid density and viscosity of $\rho_f = 1.03~g/cm^3$
and $\mu_f = 0.03~g/(cm \cdot s^2)$,
and a solid density and stiffness of $\rho_s = 1.03~g/cm^3$
and $\mu_s = 2 \cdot 10^5~g/(cm \cdot s^2)$.
The pressure amplitude is selected as $P = 583$ to yield a maximum inflow velocity of
$V_f \approx 100~cm / s$.
We further modify the dimensions of the fluid and solid domains:
the fluid domain radius is $H_i = 0.7~cm$, the solid wall thickness is $0.223~cm$
and the fluid domain length is $5.53~cm$.
Here, the Womersley number is typical for aortic blood flow, $W \approx 10.16$,
resulting in the typical Womersley flow profile
(see Figure~\ref{3D-nonlinear-tf-ts-solution-along-y-physiological-fig}),
while the Reynolds number is, $Re \approx 4813.40$.

Due to the significantly increased stiffness, the deformation of the solid is small
(see Figure~\ref{3D-nonlinear-tf-ts-solution-along-y-physiological-fig})
and the deformation rate is several orders of magnitudes smaller than the peak velocity of the fluid.
This is in line with observations that the shear deformations in the ascending aorta
are much smaller compared to the radial deformations (which are constrained here).

Although the deformation of the solid is small, it cannot be neglected because it affects
the pressure solution. For example, the solid pressure solution at / near the outlet
has a steep gradient in the radial direction.
Furthermore, similar to the previous nonlinear FSI case in this section,
the fluid and solid pressure solution is discontinuous at the interface,
see Figure~\ref{3D-nonlinear-tf-ts-solution-along-y-physiological-fig}
and Figure~\ref{3D-nonlinear-tf-ts-solution-along-t-physiological-fig}.
This indicates that the presented analytic FSI solutions are suited for benchmarking
methods that are able to capture such discontinuities as well as methods that assume continuity
at the interface.
\subsection{Numerical solution and space-time convergence}
\label{numerical-solution-t-conv-cheart-sec}
In this section, results for convergence tests for the FSI method detailed
in Section~\ref{cheart-implementation-sec} are presented.
Figure~\ref{2D-linear-tf-ts-t-conv-vf-us-lambda-fig} illustrates observed
temporal convergence of the algorithm applied to the linear FSI case
for the fluid velocity, solid displacement and Lagrange multiplier
for three different spatial resolutions.

For the \emph{coarse} refinement level, error reduction under temporal refinement is of first-order
initially, then degrades as the limit of spatial discretization accuracy is approached.
Such a degradation is not observed for the \emph{medium} and \emph{fine} spatial resolutions
for all considered time step sizes.
Furthermore, temporal convergence under space-time refinement
(with $\delta_t / \delta_x^3 = \text{const}$) is of (optimal) rate,
with $O (\delta_t)$ and $O (\delta_x^3)$ in the $L^2 (\Omega_t; L^2 (\Omega_k^0))$ norm.

Space-time convergence of the FSI algorithm for the nonlinear FSI case in three dimensions is illustrated
in Figure~\ref{3D-nonlinear-tf-ts-t-conv-vf-us-lambda-pf-ps-fig}.
Here, observed temporal convergence for the fluid velocity, solid displacement
and solid pressure are presented for three different spatial resolutions.
Similar to error reduction for the linear FSI case, temporal convergence of the non-conforming FSI algorithm
is initially first-order for the \emph{coarse} refinement level,
then degrades as discretization accuracy is approached.
For the \emph{medium} refinement level, the error can be reduced further under temporal refinement.
Error reduction in the fluid velocity and solid pressure solution, however, still degrades,
whereas the error in the solid displacement can still be decreased.
For the \emph{fine} spatial resolution, the error in the fluid velocity, solid displacement
and solid pressure variable can be reduced with optimal first-order rate
for all considered time step sizes.\footnote{The number of Newton iterations taken in the nonlinear solve
of the combination largest time step size / finest spatial mesh increased significantly.
To avoid excessive runtimes, this case was omitted in the convergence study.}
Similar to the case of two spatial dimensions,
temporal convergence under space-time refinement
(with $\delta_t / \delta_x^3 = \text{const}$) is of (optimal) rate,
with $O (\delta_t)$ and $O (\delta_x^3)$ in the $L^2 (\Omega_t; L^2 (\Omega_k^0))$ norm.
\section{Discussion}
\label{discussion-sec}
\subsection{Class of analytic solutions for FSI}
The class of analytic FSI solutions presented in Section~\ref{methodology-sec} provides
a rich and comprehensive test bed for the validation of FSI implementations.
It is comprised of permutations of linear or nonlinear, and quasi-static or transient behavior
in two and three dimensions.
Depending on the considered case, the velocity / displacement (pressure) solutions
may exhibit simplistic quadratic (linear) dependence in one spatial dimension
or more complex dependence on Bessel functions, with the pressure solution
varying in multiple spatial dimensions.

The analytic solutions are functions of the material and geometric parameters
and as such, the behavior and properties of the solutions can vary widely
(even for the same FSI solution; see Section~\ref{results-3D-analytic-sec} and \ref{results-physiological-parameters-sec}).
For example, the flow profile can be close to a parabolic profile or resemble a typical Womersley
profile; similarly, the Reynolds number can be small or (by changing, e.g., material parameters) high.
On the other hand, the solid can exhibit small
(e.g., as demonstrated in Section~\ref{results-physiological-parameters-sec})
or large nonlinear deformations.
It is even possible to select parameters
to observe higher-order deformation modes (similar to the Euler-Bernoulli mode shapes;
see Figure~\ref{3D-nonlinear-tf-ts-solution-along-y-rhos5-fig})
or to construct cases without unique solution (resonance-frequency-type solutions;
e.g., see Equation~\eqref{eq:resonance-condition}).
This property of the analytic solutions accentuates one of the avenues
to make the numerical solution process more challenging.

A commonly used benchmark problem is that of flow in an elastic tube (e.g., \cite{Womersley1957});
however, it is restricted to validating pulse wave propagation but not spatiotemporal behavior of the solutions
in the entire domain, which is possible with the class of analytic FSI solutions.
While the flow in the nonlinear FSI case is still Stokes-like (the advective term in the Navier-Stokes equations
drops out because of the assumption of no radial motion),
the stress response of the neo-Hookean solid material is nonlinear and does not simplify.
In fact, to the knowledge of the authors, the possibility to validate a nonlinear solid model
in a fluid-structure interaction framework using an analytic solution is the first of its kind.
Furthermore, the pressure solution varies along the axial direction and, in the nonlinear case,
also in the radial direction.
Thus, the analytic FSI solution cannot be reduced to (or solved in)
a lower-dimensional manifold in space.
\subsection{Analytic solutions utility for verification of FSI algorithms}
\label{analytic-solutions-utility-for-verification-of-fsi-algorithms-sec}
The class of analytic FSI solutions presented here enable the testing
of FSI implementations at various stages:
Starting from the quasi-static linear case in two dimensions, complexity can be added
and the respective implementation validated step-by-step.
For example, transient behavior can be added in the fluid and / or solid models, respectively,
or the FSI domain changed from two dimensions to three dimensions.
While the FSI geometry is relatively simple compared to other benchmarks
(e.g., \cite{TurekHron2006,HessenthalerGaddumHolubSinkusRoehrleNordsletten2017}),
it does not pose any challenges in setting up the computational model.
In particular, the analytic solutions are smooth and discretization errors stemming from discretizing
the spatial domains can be minimized
(e.g., by selecting appropriate shape functions in finite element discretizations).

Another key aspect is the possibility to change Neumann-type to Dirichlet-type boundary conditions
or vice-versa, complementing the validation of the implementation at the equation level
to a more infrastractural aspect.
The validation of the coupling constraints can further be isolated
by first validating the fluid- and solid-only cases (e.g., by setting Dirichlet boundary conditions
at the coupling boundary) and then adding the coupling constraint.
Coupling approaches can vary widely \cite{NordslettenKaySmith2010,
LandajuelaVidrascuChapelleFernandez2017,JanssonDegirmenciHoffman2017}
and affect numerical performance of an FSI implementation.
Furthermore, some FSI methods assume continuity of the pressure across the fluid / solid interface,
whereas other methods do not make such an assumption (e.g., the non-conforming monolithic FSI method
assessed in Section~\ref{numerical-solution-t-conv-cheart-sec}).
It is an open research question how such assumptions influence the numerical solution.

A benefit of analytic solutions is their usefulness for assessing
numerical tricks such as implicit-explicit splits
\cite{LandajuelaVidrascuChapelleFernandez2017} and
novel methods or time integration algorithms, such as
multigrid-reduction-in-time \cite{FriedhoffFalgoutKolevMaclachlanSchroder2012,
FalgoutFriedhoffKolevMaclachlanSchroder2014} (MGRIT).
For example, combined with convergence theory
\cite{DobrevKolevPeterssonSchroder2017,Southworth2019,
HessenthalerSouthworthNordslettenRoehrleFalgoutSchroder2019_preprint},
the class of analytic FSI solutions can be an invaluable tool for understanding the performance
of MGRIT for FSI problems and help to identify reasons of potentially slow(er) convergence,
as well as support the development of remedies thereof.
\subsection{Analytic solutions utility for spatiotemporal convergence analysis}
The availability of analytic solutions further enables the assessment of the accuracy and convergence
of FSI methods and implementations.
Even if neglecting the complex and code-specific task of implementing forcing terms for the validation
through \addtxt{so-called} manufactured FSI solutions, spatiotemporal errors of such forcing terms
(due to their numerical inclusion) may affect the validation process.
Analytic solutions make it easier to isolate potential implementation errors, as described
in Section~\ref{analytic-solutions-utility-for-verification-of-fsi-algorithms-sec}.
Conducting stability analyses under refinement is also simplified,
e.g., when investigating the stability of an explicit method and identifying time step size limits.
Performing space-time convergence studies can be difficult with numerical benchmarks
because the smoothness of solutions may not be sufficient for testing optimal convergence
and because one needs to rely on reference numerical solutions.
These difficulties can be avoided with the class of analytic solutions,
as demonstrated in Section~\ref{numerical-solution-t-conv-cheart-sec}.
By conducting a thorough convergence study, it can be shown that a given FSI method
was not only implemented correctly but that the theoretical best-case convergence rate can be achieved;
for example, the backward Euler time-discretization with quadratic / linear finite elements in space
(see Section~\ref{cheart-implementation-sec}), was shown to converge at the optimal first-order rate in time
under space-time refinement with $\delta_t / \delta_x^3 = \text{const}$
(measuring errors in $L^2 (\Omega_t; L^2 (\Omega_k^0))$).
\subsection{Study limitations}
Although the derivation of the nonlinear FSI solutions was based on the incompressible
Navier-Stokes equations as a fluid model, the advective (nonlinear) term vanishes.
Thus, any potential nonlinearities in the behavior of the fluid in numerical experiments stem from
sources like approximating the spatial domains, entirely.
In the linear case, the velocity and displacement solutions are further constant in the axial
component. While more complicated and radially varying analytic solutions would be desirable,
the assumption of no radial motion (and thus, the restriction to shear deformations)
was an essential ingredient to enable the derivation of the class of analytic solutions
presented in this work. Naturally, it would be beneficial to have more complexity in space and time
but the class of analytic solutions already represents an important first step
and provides key testing for the dynamics and nonlinear mechanical coupling.

\addtxt{These restrictions can potentially limit the benchmark's effectiveness for code testing.
Examples of scenarios where this may happen are setting either or both the advection and the 
$x$- and $y$-components of the solution to zero. It should still be noted that numerical solutions
are prone to secondary flow and displacements that are dependent on the discretization, which
should expose underlying issues in the code. In the nonlinear case, we have the benefit of additional
spatiotemporal complexity to the pressure, requiring a precise implementation in order for
convergence to be achieved. Consequently, the benchmark represents a necessary
(but not necessarily sufficient) approach to verify the code.}

\subsection{Matlab implementation}
All analytic solutions were implemented\footnote{\label{implementation-fn}The code
is available in the \emph{Supplementary Materials} and as an open-source implementation.
Bitbucket repository: \url{https://bitbucket.org/hessenthaler/fsi-analytical-solutions-matlab}.}
in Matlab R2018a to assist with the validation process.
For example, the code can be used to evaluate and export the analytic solutions for a given FSI case
to set appropriate initial and boundary conditions in a computational model,
as well as for the comparison with the numerical solution.
The material, temporal and geometric parameters can be user-defined and a user-mesh can be imported
to evaluate the analytic solutions at relevant coordinates.
We have further implemented functionality to visualize the analytic solutions similar
to Figure~\ref{3D-nonlinear-tf-ts-solution-along-y-fig} and included movies of the animated solutions
for all linear / nonlinear, 2D / 3D and quasi-static / transient cases in the \emph{Supplementary Materials}
(for parameters in Section~\ref{results-2D-analytic-sec} and Section~\ref{results-3D-analytic-sec}).
\addtxt{The provided Matlab code further validates
that the derived analytic solutions satisfy the momentum balance equations
and coupling constraints
(as detailed in Section~\ref{validation-of-derivations-seq}) using numerical differentiation.}
\section{Conclusion}\label{conclusion-sec}
In this work, we have presented a class of analytic FSI solutions that provide
a rich test bed for validation and convergence testing of linear and nonlinear FSI implementations.
It enables a step-by-step validation process (two and three dimensions, linear and nonlinear,
static and transient behavior, etc.) and further allows for isolated testing of certain parts
of an FSI implementation, such as the coupling conditions.
To the knowledge of the authors, it is the first analytic FSI solution
that includes a hyperelastic solid model.
Along with a description and derivation of the class of analytic FSI solutions,
we have demonstrated their usefulness for validating numerical FSI methods
and highlighted their value in spatiotemporal convergence testing.
An implementation of all analytic solutions is provided in the \emph{Supplementary Materials}
and as an open-source implementation\footnote{See footnote \ref{implementation-fn}.}
to assist with the validation process.
\section{Acknowledgements}
M.B. acknowledges funding from the King's College London and Imperial
College London EPSRC Centre for Doctoral Training in Medical
Imaging (EP/L015226/1), the support of Wellcome
EPSRC Centre for Medical Engineering at King's College London
(WT 203148/Z/16/Z) and of the National Institute for Health Research
(NIHR) Biomedical Research Centre award to Guy and St Thomas'
NHS Foundation Trust in partnership with King's College London.
The views expressed are those of the authors and not necessarily
those of the NHS, the NIHR or the Department of Health.
D.N. acknowledges funding form the Engineering and Physical
Sciences (EP/N011554/1 and EP/R003866/1).
\section*{Additional Information}
Declarations of interest: none.
\section*{References}
\bibliography{main}

\begin{thebibliography}{10}
\expandafter\ifx\csname url\endcsname\relax
  \def\url#1{\texttt{#1}}\fi
\expandafter\ifx\csname urlprefix\endcsname\relax\def\urlprefix{URL }\fi
\expandafter\ifx\csname href\endcsname\relax
  \def\href#1#2{#2} \def\path#1{#1}\fi

\bibitem{hirt1974arbitrary}
C.~W. Hirt, A.~A. Amsden, J.~L. Cook, {An arbitrary Lagrangian-Eulerian
  computing method for all flow speeds}, Journal of computational physics
  14~(3) (1974) 227--253.

\bibitem{donea1982arbitrary}
J.~Donea, S.~Giuliani, J.-P. Halleux, {An arbitrary Lagrangian-Eulerian finite
  element method for transient dynamic fluid-structure interactions}, Computer
  Methods in Applied Mechanics and Engineering 33~(1-3) (1982) 689--723.

\bibitem{hughes1981lagrangian}
T.~J. Hughes, W.~K. Liu, T.~K. Zimmermann, {Lagrangian-Eulerian finite element
  formulation for incompressible viscous flows}, Computer Methods in Applied
  Mechanics and Engineering 29~(3) (1981) 329--349.

\bibitem{tezduyar1991stabilized}
T.~E. Tezduyar, Stabilized finite element formulations for incompressible flow
  computations, in: Advances in Applied Mechanics, Vol.~28, Elsevier, 1991, pp.
  1--44.

\bibitem{tezduyar1992new}
T.~E. Tezduyar, M.~Behr, S.~Mittal, J.~Liou, A new strategy for finite element
  computations involving moving boundaries and interfaces -- the
  deforming-spatial-domain/space-time procedure: Ii. computation of
  free-surface flows, two-liquid flows, and flows with drifting cylinders,
  Computer Methods in Applied Mechanics and Engineering 94~(3) (1992) 353--371.

\bibitem{takizawa2011multiscale}
K.~Takizawa, T.~E. Tezduyar, Multiscale space-time fluid-structure interaction
  techniques, Computational Mechanics 48~(3) (2011) 247--267.

\bibitem{HoffmanJanssonStoeckli2011}
J.~Hoffman, J.~Jansson, M.~St{\"o}ckli, Unified continuum modeling of
  fluid-structure interaction, Mathematical Models and Methods in Applied
  Sciences 21~(03) (2011) 491--513.

\bibitem{JanssonDegirmenciHoffman2017}
J.~Jansson, N.~C. Degirmenci, J.~Hoffman, Adaptive unified continuum fem
  modeling of a 3d fsi benchmark problem, International Journal for Numerical
  Methods in Biomedical Engineering 33~(9) (2017) e2851.

\bibitem{peskin1973flow}
C.~S. Peskin, Flow patterns around heart valves: a digital computer method for
  solving the equations of motion, IEEE Transactions on Biomedical
  Engineering~(4) (1973) 316--317.

\bibitem{peskin2002immersed}
C.~S. Peskin, The immersed boundary method, Acta Numerica 11 (2002) 479--517.

\bibitem{mittal2005immersed}
R.~Mittal, G.~Iaccarino, Immersed boundary methods, Annual Review of Fluid
  Mechanics 37 (2005) 239--261.

\bibitem{glowinski1994fictitious}
R.~Glowinski, T.-W. Pan, J.~Periaux, {A fictitious domain method for Dirichlet
  problem and applications}, Computer Methods in Applied Mechanics and
  Engineering 111~(3-4) (1994) 283--303.

\bibitem{glowinski1999distributed}
R.~Glowinski, T.-W. Pan, T.~I. Hesla, D.~D. Joseph, {A distributed Lagrange
  multiplier/fictitious domain method for particulate flows}, International
  Journal of Multiphase Flow 25~(5) (1999) 755--794.

\bibitem{gil2010immersed}
A.~J. Gil, A.~A. Carre{\~n}o, J.~Bonet, O.~Hassan, {The immersed structural
  potential method for haemodynamic applications}, Journal of Computational
  Physics 229~(22) (2010) 8613--8641.

\bibitem{gil2013enhanced}
A.~J. Gil, A.~A. Carreno, J.~Bonet, O.~Hassan, An enhanced immersed structural
  potential method for fluid-structure interaction, Journal of Computational
  Physics 250 (2013) 178--205.

\bibitem{StegerDoughertyBenek1983}
J.~L. Steger, F.~C. Dougherty, J.~A. Benek, A chimera grid scheme, in: Advances
  in Grid Generation, Vol. ASME FED-5, 1983, pp. 59--69.

\bibitem{StegerBenek1987}
J.~L. Steger, J.~A. Benek, On the use of composite grid schemes in
  computational aerodynamics, Computer Methods in Applied Mechanics and
  Engineering 64~(1-3) (1987) 301--320.

\bibitem{ChesshireHenshaw1990}
G.~Chesshire, W.~D. Henshaw, Composite overlapping meshes for the solution of
  partial differential equations, Journal of Computational Physics 90~(1)
  (1990) 1--64.

\bibitem{HouzeauxCodina2003}
G.~Houzeaux, R.~Codina, {A Chimera method based on a Dirichlet/Neumann(Robin)
  coupling for the Navier-Stokes equations}, Computer Methods in Applied
  Mechanics and Engineering 192~(31-32) (2003) 3343--3377.

\bibitem{WallGamnitzerGerstenberger2008}
W.~A. Wall, P.~Gamnitzer, A.~Gerstenberger, Fluid-structure interaction
  approaches on fixed grids based on two different domain decomposition ideas,
  International Journal of Computational Fluid Dynamics 22~(6) (2008) 411--427.

\bibitem{BalmusMassingHoffmanRazaviNordsletten2019_preprint}
M.~Balmus, A.~Massing, J.~Hoffman, R.~Razavi, D.~Nordsletten,
  \href{https://arxiv.org/abs/1902.06168}{A partition of unity approach to
  fluid mechanics and fluid-structure interaction}, arXiv preprint
  arXiv:1902.06168.
\newline\urlprefix\url{https://arxiv.org/abs/1902.06168}

\bibitem{schreiber1983driven}
R.~Schreiber, H.~B. Keller, Driven cavity flows by efficient numerical
  techniques, Journal of Computational Physics 49~(2) (1983) 310--333.

\bibitem{schreiber1983spurious}
R.~Schreiber, H.~B. Keller, Spurious solutions in driven cavity calculations,
  Journal of Computational Physics 49 (1983) 165--172.

\bibitem{kim1985application}
J.~Kim, P.~Moin, {Application of a fractional-step method to incompressible
  Navier-Stokes equations}, Journal of Computational Physics 59~(2) (1985)
  308--323.

\bibitem{GhattasLi1995}
O.~Ghattas, X.~Li, A variational finite element method for stationary nonlinear
  fluid-solid interaction, Journal of Computational Physics 121~(2) (1995)
  347--356.

\bibitem{Wall1999}
W.~A. Wall, {Fluid-Struktur-Interaktion mit stabilisierten Finiten Elementen},
  Ph.D. thesis, Institut f\"ur Baustatik, Universit\"at Stuttgart (1999).

\bibitem{Mok2001}
D.~P. Mok, {Partitionierte L{\"o}sungsans{\"a}tze in der Strukturdynamik und
  der Fluid-Struktur-Interaktion}, Ph.D. thesis, Institut f\"ur Baustatik,
  Universit\"at Stuttgart (2001).

\bibitem{TurekHron2006}
S.~Turek, J.~Hron, Proposal for numerical benchmarkingof fluid-structure
  interaction between an elastic object and laminar incompressible flow, in:
  Fluid-structure interaction, Springer, 2006, pp. 371--385.

\bibitem{BatheLedezma2007}
K.-J. Bathe, G.~A. Ledezma, Benchmark problems for incompressible fluid flows
  with structural interactions, Computers \& Structures 85~(11-14) (2007)
  628--644.

\bibitem{HessenthalerGaddumHolubSinkusRoehrleNordsletten2017}
A.~Hessenthaler, N.~Gaddum, O.~Holub, R.~Sinkus, O.~R{\"o}hrle, D.~Nordsletten,
  Experiment for validation of fluid-structure interaction models and
  algorithms, International Journal for Numerical Methods in Biomedical
  Engineering.

\bibitem{HeilHazelBoyle2008}
M.~Heil, A.~L. Hazel, J.~Boyle, Solvers for large-displacement fluid-structure
  interaction problems: segregated versus monolithic approaches, Computational
  Mechanics 43~(1) (2008) 91--101.

\bibitem{TurekHronRazzaqWobkerSchaefer2011}
S.~Turek, J.~Hron, M.~Razzaq, H.~Wobker, M.~Sch{\"a}fer, Numerical benchmarking
  of fluid-structure interaction: A comparison of different discretization and
  solution approaches, in: Fluid Structure Interaction II, Springer, 2011, pp.
  413--424.

\bibitem{BertramTscherry2006}
C.~D. Bertram, J.~Tscherry, The onset of flow-rate limitation and flow-induced
  oscillations in collapsible tubes, Journal of Fluids and Structures 22~(8)
  (2006) 1029--1045.

\bibitem{GomesLienhart2006}
J.~P. Gomes, H.~Lienhart, Experimental study on a fluid-structure interaction
  reference test case, in: Fluid-Structure Interaction I -- Modelling,
  Simulation, Optimization, Vol.~53, Lecture Notes in Computational Science and
  Engineering, Springer, 2006, pp. 356--370.

\bibitem{IdelsohnMartiSoutoiglesiasOnate2008}
S.~Idelsohn, J.~Marti, A.~Souto-Iglesias, E.~Onate, Interaction between an
  elastic structure and free-surface flows: experimental versus numerical
  comparisons using the pfem, Computational Mechanics 43~(1) (2008) 125--132.

\bibitem{GomesLienhart2010}
J.~P. Gomes, H.~Lienhart, Experimental benchmark: Self-excited fluid-structure
  interaction test cases, in: Fluid-Structure Interaction II -- Modelling,
  Simulation, Optimization, Vol.~73, Lecture Notes in Computational Science and
  Engineering, Springer, 2010, pp. 383--411.

\bibitem{NayerKalmbachBreuerSicklingerWuechner2014}
G.~De~Nayer, A.~Kalmbach, M.~Breuer, S.~Sicklinger, R.~W{\"u}chner, {Flow past
  a cylinder with a flexible splitter plate: A complementary
  experimental-numerical investigation and a new FSI test case (FSI-PfS-1a)},
  Computers \& Fluids 99 (2014) 18--43.

\bibitem{HessenthalerRoehrleNordsletten2017}
A.~Hessenthaler, O.~R{\"o}hrle, D.~Nordsletten, Validation of a non-conforming
  monolithic fluid-structure interaction method using phase-contrast mri,
  International Journal for Numerical Methods in Biomedical Engineering 33~(8).

\bibitem{steinberg1985symbolic}
S.~Steinberg, P.~J. Roache, Symbolic manipulation and computational fluid
  dynamics, Journal of Computational Physics 57~(2) (1985) 251--284.

\bibitem{roache2002code}
P.~J. Roache, Code verification by the method of manufactured solutions,
  Journal of Fluids Engineering 124~(1) (2002) 4--10.

\bibitem{salari2000code}
K.~Salari, P.~Knupp, Code verification by the method of manufactured solutions,
  Tech. rep., Sandia National Laboratories, Albuquerque, NM 87185 (US); Sandia
  National Laboratories, Livermore, CA 94550 (US) (2000).

\bibitem{ethier1994exact}
C.~R. Ethier, D.~Steinman, {Exact fully 3D Navier-Stokes solutions for
  benchmarking}, International Journal for Numerical Methods in Fluids 19~(5)
  (1994) 369--375.

\bibitem{womersley1955method}
J.~R. Womersley, Method for the calculation of velocity, rate of flow and
  viscous drag in arteries when the pressure gradient is known, The Journal of
  Physiology 127~(3) (1955) 553--563.

\bibitem{NordslettenKaySmith2010}
D.~Nordsletten, D.~Kay, N.~Smith, A non-conforming monolithic finite element
  method for problems of coupled mechanics, Journal of Computational Physics
  229~(20) (2010) 7571--7593.

\bibitem{LeeEtAl2016}
J.~Lee, A.~Cookson, I.~Roy, E.~Kerfoot, L.~Asner, G.~Vigueras, T.~Sochi,
  S.~Deparis, C.~Michler, N.~P. Smith, D.~A. Nordsletten, Multiphysics
  computational modeling in cheart, SIAM Journal on Scientific Computing 38~(3)
  (2016) C150--C178.

\bibitem{hughes2005conservation}
T.~J. Hughes, G.~N. Wells, Conservation properties for the {G}alerkin and
  stabilised forms of the advection-diffusion and incompressible
  {N}avier-{S}tokes equations, Computer Methods in Applied Mechanics and
  Engineering 194~(9-11) (2005) 1141--1159.

\bibitem{moghadam2011comparison}
M.~E. Moghadam, Y.~Bazilevs, T.-Y. Hsia, I.~E. Vignon-Clementel, A.~L. Marsden,
  et~al., A comparison of outlet boundary treatments for prevention of backflow
  divergence with relevance to blood flow simulations, Computational Mechanics
  48~(3) (2011) 277--291.

\bibitem{HessenthalerFriedhoffRoehrleNordsletten2016}
A.~Hessenthaler, S.~Friedhoff, O.~R{\"o}hrle, D.~A. Nordsletten, 3d
  fluid-structure interaction experiment and benchmark results, PAMM 16~(1)
  (2016) 451--452.

\bibitem{AmestoyDuffLexcellent2000}
P.~R. Amestoy, I.~S. Duff, J.-Y. L'excellent, Multifrontal parallel distributed
  symmetric and unsymmetric solvers, Computer Methods in Applied Mechanics and
  Engineering 184~(2-4) (2000) 501--520.

\bibitem{Shamanskii1967}
V.~Shamanskii, A modification of newton's method, Ukrainian Mathematical
  Journal 19~(1) (1967) 118--122.

\bibitem{BauchauCraig2009}
O.~A. Bauchau, J.~I. Craig, Euler-bernoulli beam theory, in: Structural
  Analysis, Springer, 2009, pp. 173--221.

\bibitem{LiKamelRandoAndersonKumbasarLimaBluemke2003}
A.~E. Li, I.~Kamel, F.~Rando, M.~Anderson, B.~Kumbasar, J.~A. Lima, D.~A.
  Bluemke, Using mri to assess aortic wall thickness in the multiethnic study
  of atherosclerosis: distribution by race, sex, and age, American Journal of
  Roentgenology 182~(3) (2004) 593--597.

\bibitem{ReymondCrosettoDeparisQuarteroniStergiopulos2012}
P.~Reymond, P.~Crosetto, S.~Deparis, A.~Quarteroni, N.~Stergiopulos,
  Physiological simulation of blood flow in the aorta: comparison of
  hemodynamic indices as predicted by 3-d fsi, 3-d rigid wall and 1-d models,
  Medical Engineering \& Physics 35~(6) (2013) 784--791.

\bibitem{Womersley1957}
J.~R. Womersley, Oscillatory {F}low in {A}rteries: the {C}onstrained {E}lastic
  {T}ube as a {M}odel of {A}rterial {F}low and {P}ulse {T}ransmission, Physics
  in Medicine and Biology 2~(2) (1957) 178--187.

\bibitem{LandajuelaVidrascuChapelleFernandez2017}
M.~Landajuela, M.~Vidrascu, D.~Chapelle, M.~A. Fern{\'a}ndez, Coupling schemes
  for the fsi forward prediction challenge: comparative study and validation,
  International Journal for Numerical Methods in Biomedical Engineering 33~(4)
  (2017) e2813.

\bibitem{FriedhoffFalgoutKolevMaclachlanSchroder2012}
S.~Friedhoff, R.~D. Falgout, T.~V. Kolev, S.~MacLachlan, J.~B. Schroder, A
  multigrid-in-time algorithm for solving evolution equations in parallel,
  Tech. rep., Lawrence Livermore National Lab.(LLNL), Livermore, CA (United
  States) (2012).

\bibitem{FalgoutFriedhoffKolevMaclachlanSchroder2014}
R.~D. Falgout, S.~Friedhoff, T.~V. Kolev, S.~P. MacLachlan, J.~B. Schroder,
  Parallel time integration with multigrid, SIAM Journal on Scientific
  Computing 36~(6) (2014) C635--C661.

\bibitem{DobrevKolevPeterssonSchroder2017}
V.~A. Dobrev, T.~V. Kolev, N.~A. Petersson, J.~B. Schroder, Two-level
  convergence theory for multigrid reduction in time ({MGRIT}), SIAM Journal on
  Scientific Computing 39~(5) (2017) S501--S527.

\bibitem{Southworth2019}
B.~S. Southworth, Necessary conditions and tight two-level convergence bounds
  for parareal and multigrid reduction in time, SIAM J. on Matrix Analysis and
  Applications.

\bibitem{HessenthalerSouthworthNordslettenRoehrleFalgoutSchroder2019_preprint}
A.~Hessenthaler, B.~S. Southworth, D.~Nordsletten, O.~R{\"o}hrle, R.~D.
  Falgout, J.~B. Schroder, \href{https://arxiv.org/abs/1812.11508}{Multilevel
  convergence analysis of multigrid-reduction-in-time}, arXiv preprint
  arXiv:1812.11508.
\newline\urlprefix\url{https://arxiv.org/abs/1812.11508}

\end{thebibliography}
\newpage
\section*{Tables}
\FloatBarrier
\begin{table}[ht!]
	\centering	\begin{tabularx}{0.7\textwidth}{X X}
		\hline\hline
		\multicolumn{2}{l}{\textbf{Physical and numerical parameters}}
		\\
		\hline \hline
		Material properties    \\
		\hline
		$\rho_f$ & Fluid density \\
		$\mu_f$ & Fluid viscosity \\
		$\rho_s$ & Solid density \\
		$\mu_s$ & Solid stiffness
		\\
		\hline
		Discretization \\
		\hline
		$\delta_t$ & Time step size \\
		$\delta_x, \delta_y, \delta_z$ & Spatial step size \\
		\hline
		\multicolumn{2}{l}{Domain dimensions and forces}
		\\
		\hline
		$H_i$ & Fluid domain height \\
		$H_o$ & Fluid-solid domain height \\
		$L$ & Domain length \\
		$T$ & Temporal cycle length \\
		$P$ & Pressure over Domain length \\
		\hline\hline
		\textbf{Derived constants} \\
		\hline\hline
		$k_f$ & $\sqrt{\rho_f i \omega / \mu_f}$ \\
		 $k_s$ & $\omega \sqrt{\rho_s / \mu_s}$ \\
		 $\alpha$ & $e^{k_f H_i} + e^{-k_f H_i}$ \\
		 $\beta$ & $\mu_f k_f \left(e^{k_f H_i} - e^{-k_f H_i}\right)$ \\
		 $\gamma$ & $J_{0,s}^r/ Y_{0,s}^r$ \\
		 $J_{0,s}^r$ & $J_0(- k_s H_o)$ \\
		 $Y_{0,s}^r$ & $Y_0(- k_s H_o)$ \\
		 $J^*_{0,f}$ & $J_0(i k_f H_i)$ \\
		 $J^*_{1,f}$ & $k_f J_1(i k_f H_i)$ \\
		 $J_{0,s}^*$ & $J_0(- k_s H_i)$ \\
		 $J_{1,s}^*$ & $i k_s J_1(- k_s H_i)$ \\
		 $Y_{0,s}^*$ & $Y_0(- k_s H_i)$ \\
		 $Y_{1,s}^*$ & $i k_s Y_1(- k_s H_i)$ \\
		 $\Delta_0$ & $J_{0,s}^* - \gamma Y_{0,s}^*$ \\
		 $\Delta_1$ & $J_{1,s}^* - \gamma Y_{1,s}^*$ \\
		 $\nu_0$ & $Y_{0,s}^* / Y_{0,s}^r$ \\
		 $\nu_1$ & $Y_{1,s}^* / Y_{0,s}^r$ \\
		 $\xi_1$ & $\sin(k_sH_i)+\cot(k_sH_o)\cos(k_sH_i)$ \\
		 $\xi_2$ & $\cot(k_sH_o)\sin(k_sH_i) - \cos(k_sH_i)$ \\
		 $\zeta_1$ & $\csc(k_sH_o)\cos(k_sH_i)$ \\
		 $\zeta_2$ & $1 - \sin(k_sH_i)\csc(k_sH_o)$
		 \\
		 \hline\hline
	\end{tabularx}
	\caption{Material parameters and constants.}
	\label{nomenclature-tab}
\end{table}
\FloatBarrier
\FloatBarrier
\begin{table}[ht!]
    \centering
    \begin{tabular}{ l | c | c }
        & 2D & 3D \\
        \hline
        Time step size & $\delta_t \in \{ 0.064, 0.032, \ldots, 0.001 \}$
        & $\delta_t \in \{ 0.04096, 0.02048, \ldots, 0.00064 \}$ \\
                Number of time steps & $N_t \in \{ 40, 80, \ldots, 2560 \}$
        & $N_t \in \{ 175, 350, \ldots, 11200 \}$ \\
                Spatial step size & $\delta_x, \delta_y \in \{ 0.1, 0.05, 0.025 \}$
        & $\delta_x, \delta_y, \delta_z \in \{ 0.314, 0.157, 0.079 \}$ \\
                Number of spatial DOFs & $N_x \in \{ 1288, 4728, 18088 \}$ & $N_x \in \{ 32793, 245009, 1679243 \}$
    \end{tabular}
    \caption{The spatial and temporal step sizes, $\delta_x, \delta_y, \delta_z$ and $\delta_t$,
        and numbers of spatial and temporal degrees-of-freedom (DOFs), $N_x$ and $N_t$,
        for the numerical experiments.
        For the 3D case, the spatial step size is the approximate step size at the fluid / solid interface.}
    \label{space-time-discretization-tab}
\end{table}
\FloatBarrier
\newpage
\section*{Figures}
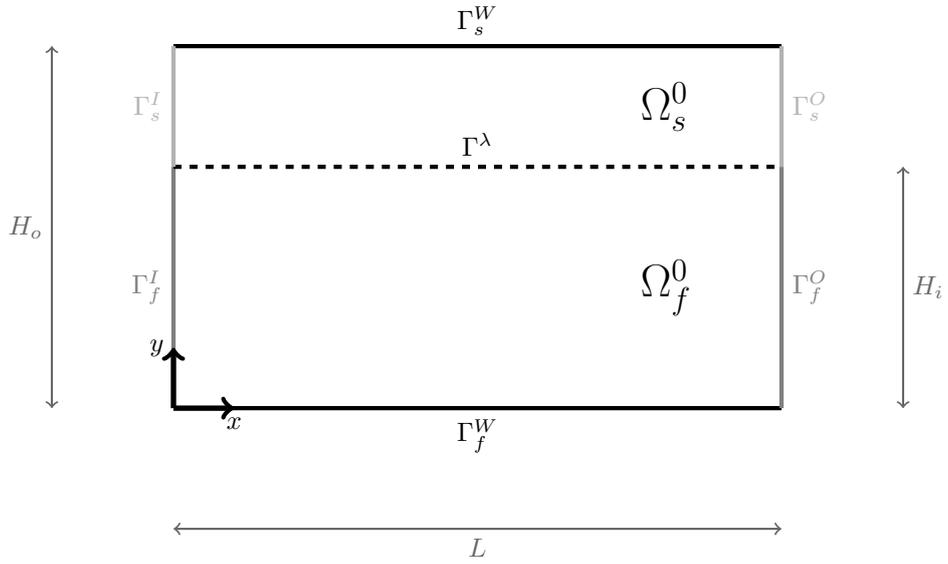
\begin{figure}[ht!]
    \centering
    \begin{tikzpicture}[scale=0.8]
                \def\L{10.0};
        \def\axl{0.1*\L};
        \def\ri{4.0};
        \def\ro{6.0};
                \draw[ultra thick, black] (0,0) -- (\L,0);
        \draw[ultra thick, black, below] (0.5*\L,0) node {$\Gamma_f^W$};
                \draw[ultra thick, black!50] (0,0) -- (0,\ri);
        \draw[ultra thick, black!50, left] (0,0.5*\ri) node {$\Gamma_f^I$};
                \draw[ultra thick, black!50] (\L,0) -- (\L,\ri);
        \draw[ultra thick, black!50, right] (\L,0.5*\ri) node {$\Gamma_f^O$};
                \draw[ultra thick, black, dashed] (0,\ri) -- (\L,\ri);
                        \draw[ultra thick, black, above] (0.5*\L,\ri) node {$\Gamma^\lambda$};
                \draw[ultra thick, black] (0,\ro) -- (\L,\ro);
        \draw[ultra thick, black, above] (0.5*\L,\ro) node {$\Gamma_s^W$};
                \draw[ultra thick, black!30] (0,\ri) -- (0,\ro);
        \draw[ultra thick, black!30, left] (0,0.5*\ro+0.5*\ri) node {$\Gamma_s^I$};
                \draw[ultra thick, black!30] (\L,\ri) -- (\L,\ro);
        \draw[ultra thick, black!30, right] (\L,0.5*\ro+0.5*\ri) node {$\Gamma_s^O$};
                \draw[ultra thick, black, right] (0.75*\L,0.5*\ri+0.5*\ro) node { \LARGE $\Omega_s^0$};
        \draw[ultra thick, black, right] (0.75*\L,0.5*\ri) node { \LARGE $\Omega_f^0$};
                \draw[->, line width=2pt,black] (0,0) -- (\axl,0);
        \draw[black, below] (\axl,0) node {$x$};
        \draw[->, line width=2pt,black] (0,0) -- (0,\axl);
        \draw[black, left] (0,\axl) node {$y$};
                \draw[<->, thick, black!60] (0.0,-2*\axl) -- (\L,-2*\axl);
        \draw[black!60, below] (0.5*\L,-2*\axl) node {$L$};
        \draw[<->, thick, black!60] (-2*\axl,0) -- (-2*\axl,\ro);
        \draw[black!60, left] (-2*\axl,0.5*\ro) node {$H_o$};
        \draw[<->, thick, black!60] (\L+2*\axl,0) -- (\L+2*\axl,\ri);
        \draw[black!60, right] (\L+2*\axl,0.5*\ri) node {$H_i$};
    \end{tikzpicture}
        \caption{Fluid and solid reference domains, $\Omega_f^0$ and $\Omega_s^0$, in two dimensions
            with respective boundaries at the inlet ($\Gamma_f^I$ and $\Gamma_s^I$),
            the outlet ($\Gamma_f^O$ and $\Gamma_s^O$) and the wall ($\Gamma_f^W$ and $\Gamma_s^W$).
            The common interface boundary is denoted as $\Gamma^\lambda$.
            Further, the domain length is given as $L$, the fluid domain height as $H_i$
            and the fluid / solid domain height as $H_o$.}
    \label{domain-2D-fig}
\end{figure}
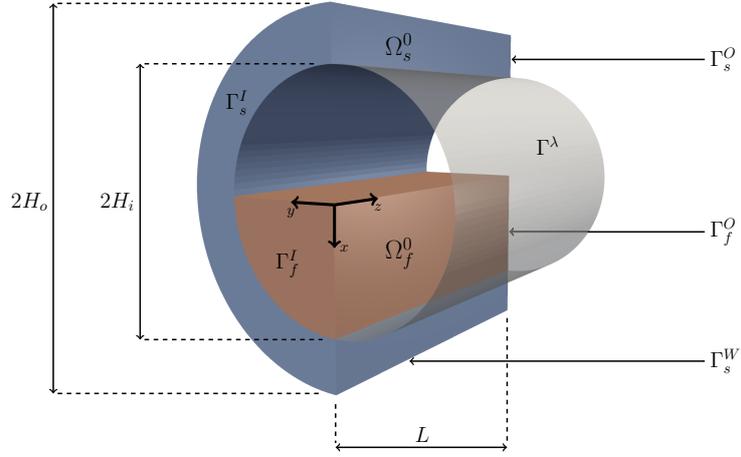
\begin{figure}[ht!]
    \centering
\resizebox{0.6\textwidth}{!}{
    \begin{tikzpicture}
                \node [anchor=south west] at (0,0) { \includegraphics[width=0.6\linewidth]{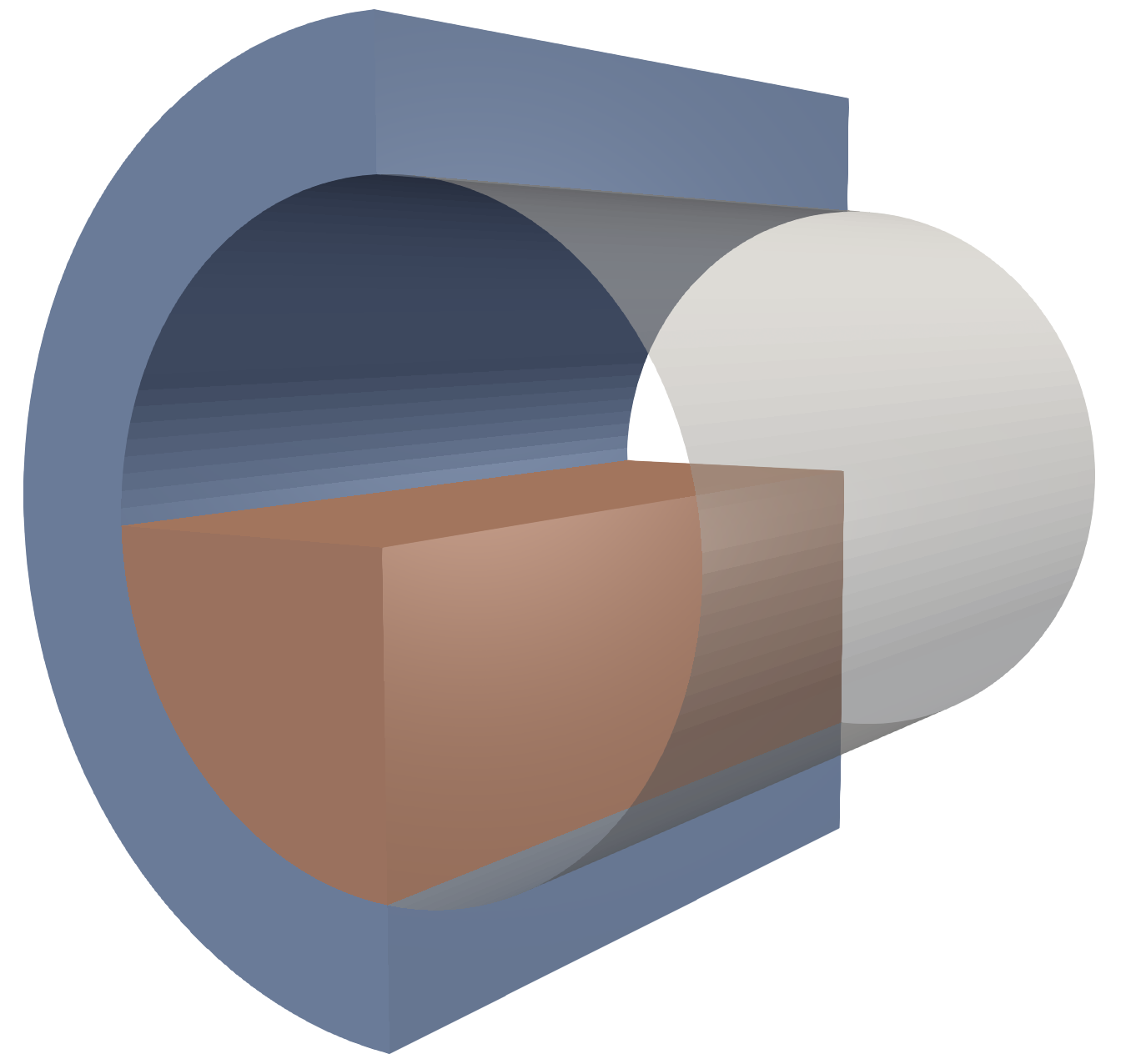} };
                \draw[->, line width=2pt,black] (3.475,4.625) -- (4.475,4.775);
        \draw[black, below] (4.475,4.775) node {$z$};
        \draw[->, line width=2pt,black] (3.475,4.625) -- (2.475,4.69);
       \draw[black, below] (2.475,4.675) node {$y$};
        \draw[->, line width=2pt,black] (3.475,4.625) -- (3.475,3.625);
        \draw[black, right] (3.4755,3.625) node {$x$};
                \draw[ultra thick, black, right] (4.5,3.5) node { \LARGE $\Omega_f^0$};
        \draw[ultra thick, black, right] (4.5,8.275) node { \LARGE $\Omega_s^0$};
                        \draw[black, right] (8,6) node {\Large $\Gamma^\lambda$};
                \draw[black, right] (2,3.25) node {\Large $\Gamma_f^I$};
                \draw[<-, line width=1pt,black] (7.5,4) -- (12,4);
        \draw[<-, line width=1pt,black!60] (7.5,4) -- (9.31,4);
        \draw[black, right] (12,4) node {\Large $\Gamma_f^O$};
                \draw[black, right] (0.85,7) node {\Large $\Gamma_s^I$};
                \draw[<-, line width=1pt,black] (7.55,8) -- (12,8);
        \draw[black, right] (12,8) node {\Large $\Gamma_s^O$};
                \draw[<-, line width=1pt,black] (5.2,1) -- (12,1);
        \draw[black, right] (12,1) node {\Large $\Gamma_s^W$};
                \draw[dashed, line width=1pt, black] (-0.9,1.5) -- (3.2,1.5);
        \draw[dashed, line width=1pt, black] (-0.9,7.9) -- (3,7.9);
        \draw[<->, line width=1pt, black] (-1,1.5) -- (-1,7.9);
        \draw[black, left] (-1,4.7) node {\Large $2 H_i$};
                \draw[dashed, line width=1pt, black] (-2.9,0.25) -- (3.2,0.25);
        \draw[dashed, line width=1pt, black] (-2.9,9.3) -- (3,9.3);
        \draw[<->, line width=1pt, black] (-3,0.25) -- (-3,9.3);
        \draw[black, left] (-3,4.7) node {\Large $2 H_o$};
                \draw[dashed, line width=1pt, black] (3.5,0) -- (3.5,-0.9);
        \draw[dashed, line width=1pt, black] (7.45,2) -- (7.45,-0.9);
        \draw[<->, line width=1pt, black] (3.5,-1) -- (7.45,-1);
        \draw[black, above] (5.5,-1) node {\Large $L$};
    \end{tikzpicture}
}
        \caption{Fluid and solid reference domains in three dimensions:
            The fluid reference domain $\Omega_f^0$ is shown in red for $x, y > 0$,
            the solid reference domain $\Omega_s^0$ is shown in blue for $y > 0$,
            and the the common interface boundary is indicated in opaque gray.
            The respective boundaries are denoted as $\Gamma_f^I$ and $\Gamma_s^I$ at the inlet,
            $\Gamma_f^O$ and $\Gamma_s^O$ at the outlet, and $\Gamma_f^W$ and $\Gamma_s^W$ at the wall.
            Furthermore, the domain length is given as $L$, the fluid domain radius as $H_i$
            and the fluid / solid domain radius as $H_o$.}
    \label{domain-3D-fig}
\end{figure}
\iffigure
\FloatBarrier
\begin{figure}[ht!]
    \centering
    \includegraphics[width=0.3\linewidth]{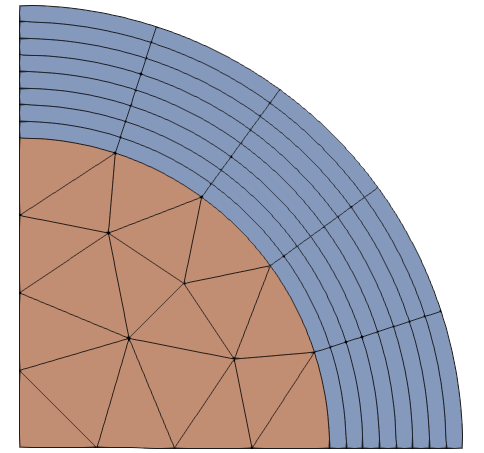}\quad
    \includegraphics[width=0.3\linewidth]{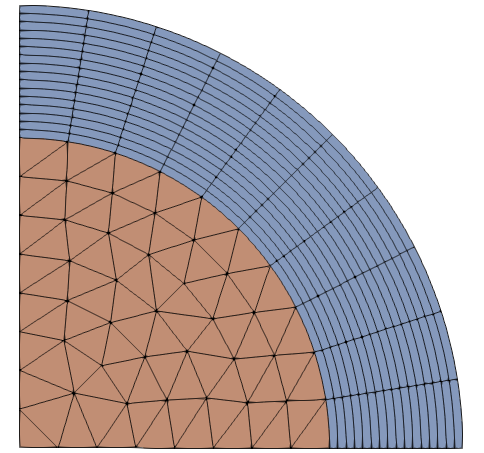}\quad
    \includegraphics[width=0.3\linewidth]{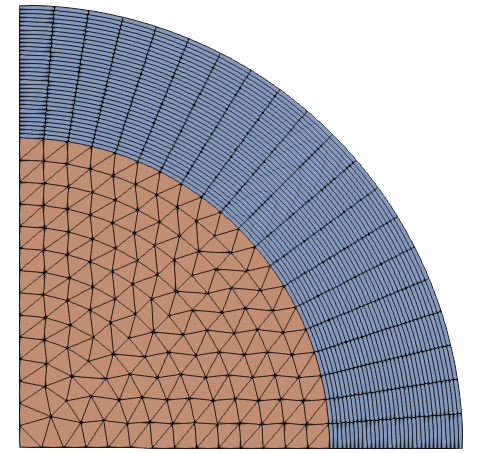}
    \caption{\emph{Coarse}, \emph{medium} and \emph{fine} mesh refinement levels
        for simulating the nonlinear transient FSI case in three dimensions.}
    \label{mesh-refinements-fig}
\end{figure}
\FloatBarrier
\fi
\iffigure
\begin{figure}[ht!]
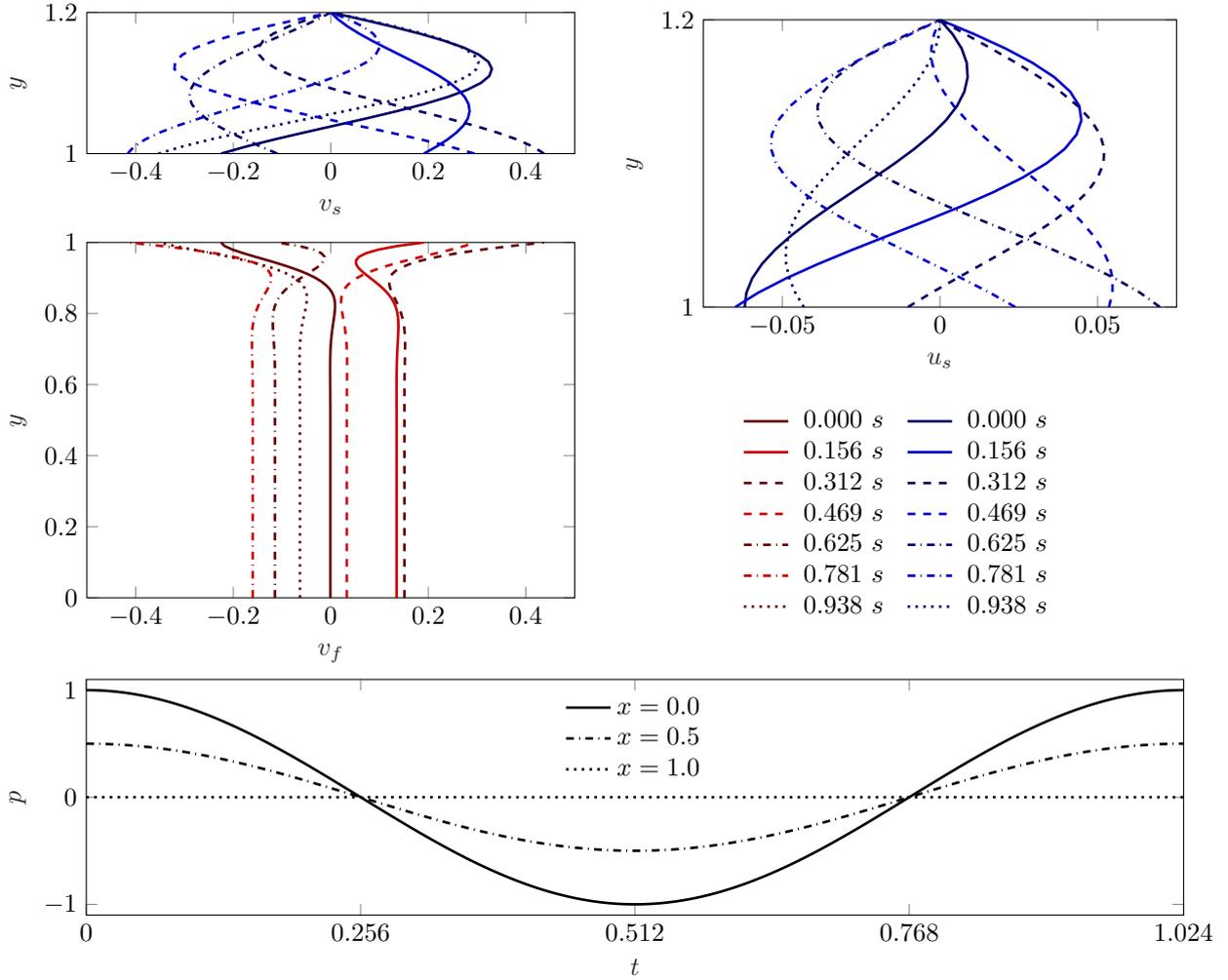

    \centering
    \begin{minipage}{\textwidth}
        \hspace{-10pt}
        \begin{minipage}{0.5\textwidth}
            \centering
            ~\\[3ex]
            \setlength{\figurewidth}{0.8\textwidth}
            \setlength{\figureheight}{0.35\textwidth}
    \caption{The analytic solution for the transient linear FSI case in two dimensions
        with density $\rho_f = \rho_s = 1$, fluid viscosity $\mu_f = 0.01$, solid stiffness $\mu_s = 0.1$
        and cycle length $T = 1.024$:
        Fluid and solid velocity, $v_f$ (left) and $v_s$ (top left), and solid displacement, $u_s$ (top right),
        along the $y$-axis, and fluid / solid pressure, $p_f$ and $p_s$ (bottom), over time $t$ at three positions
        along the $x$-axis.}
    \label{2D-linear-tf-ts-solution-along-y-fig}
\end{figure}
\FloatBarrier
\fi
\iffigure
\FloatBarrier
\begin{figure}[ht!]
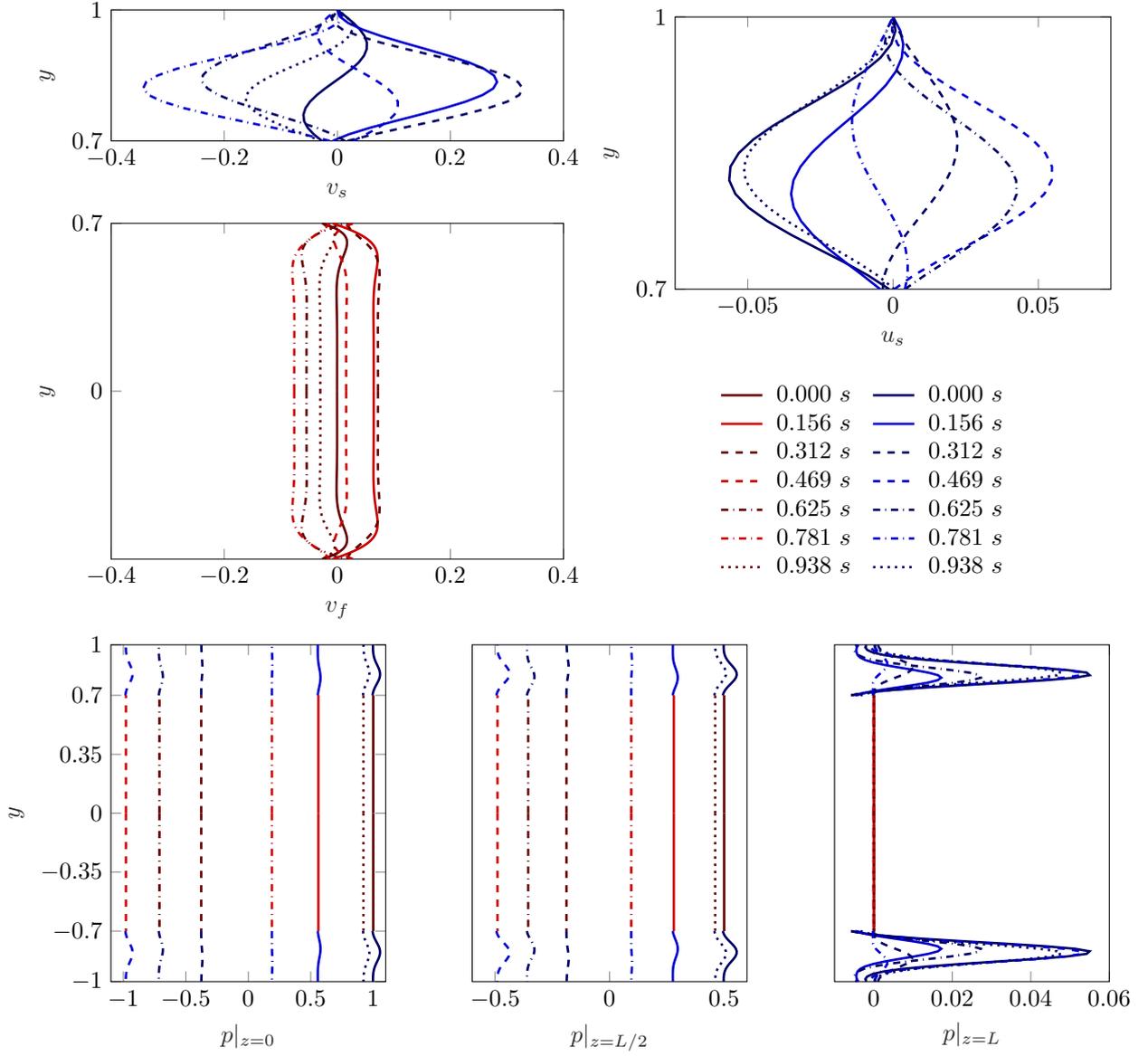

    \centering
    \begin{minipage}{\textwidth}
        \hspace{8.5pt}
        \begin{minipage}{0.5\textwidth}
            \centering
            ~\\[3ex]
            \setlength{\figurewidth}{0.8\textwidth}
            \setlength{\figureheight}{0.35\textwidth}

    \end{minipage}
    \caption{The analytic solution for the transient nonlinear FSI case in three dimensions
        with fluid and solid density $\rho_f = 2.1$ and $\rho_s = 1$, fluid viscosity $\mu_f = 0.03$,
        solid stiffness $\mu_s = 0.1$ and cycle length $T = 1.024$:
        Fluid and solid velocity, $v_f$ (left) and $v_s$ (top left), and solid displacement, $u_s$ (top right),
        along the $y$-axis.
        Further, fluid and solid pressure, $p_f$ and $p_s$ (bottom) over time $t$
        at three positions along the $z$-axis.}
    \label{3D-nonlinear-tf-ts-solution-along-y-fig}
\end{figure}
\FloatBarrier
\FloatBarrier
\begin{figure}[ht!]
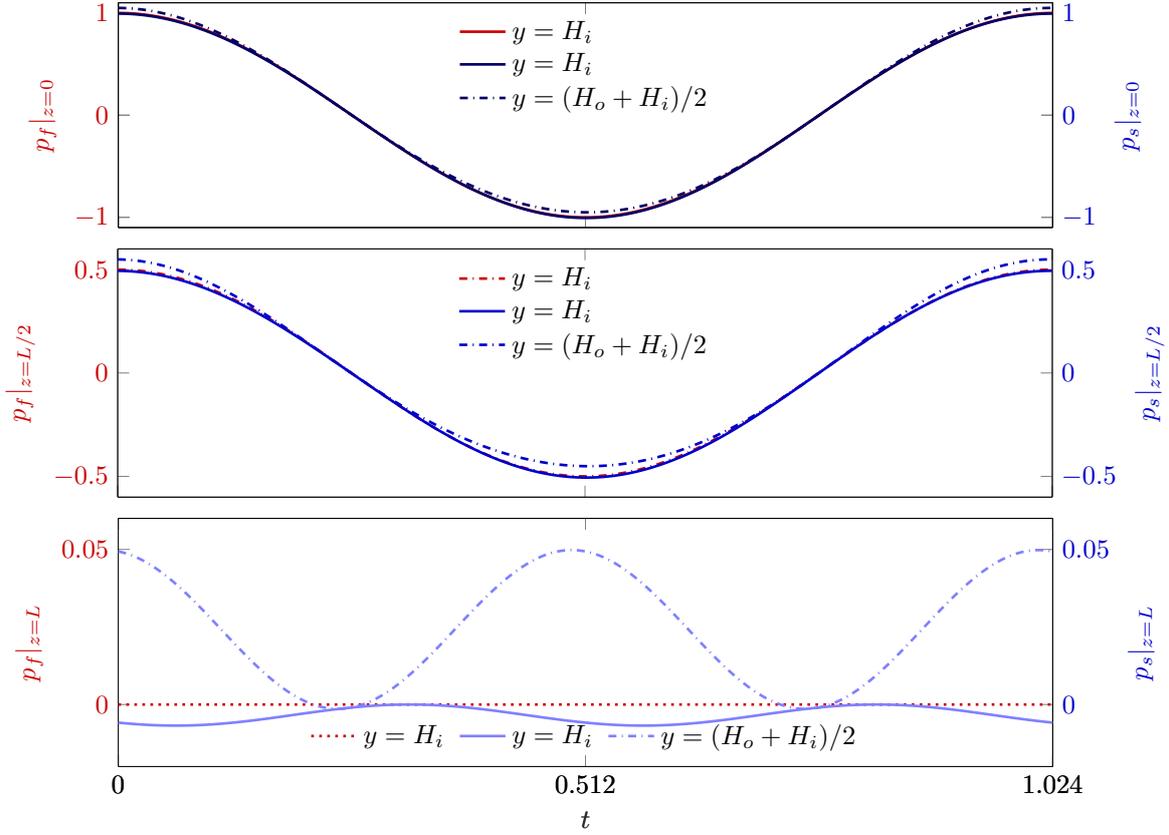

    \centering
    \setlength{\figurewidth}{0.7475\textwidth}
    \setlength{\figureheight}{0.2\textwidth}
    \caption{The analytic solution for the transient nonlinear FSI case in three dimensions
        with fluid and solid density $\rho_f = 2.1$ and $\rho_s = 1$, fluid viscosity $\mu_f = 0.03$,
        solid stiffness $\mu_s = 0.1$ and cycle length $T = 1.024$:
        The fluid and solid pressure, $p_f$ and $p_s$, over time $t$ at
        three positions along the $z$-axis: At the inlet for $z = 0$, the mid-way point at $z = L / 2$,
        and at the outlet $z = L$ (top to bottom).}
    \label{3D-nonlinear-tf-ts-solution-along-t-fig}
\end{figure}
\FloatBarrier
\fi
\iffigure
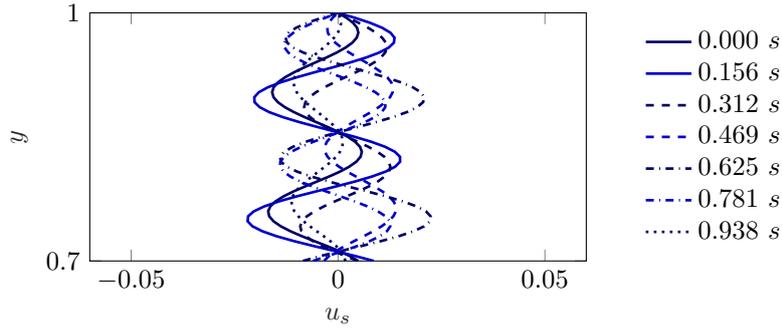
\begin{figure}[ht!]
    \centering
    \begin{minipage}{\textwidth}
        \centering
        \setlength{\figurewidth}{0.4\textwidth}
        \setlength{\figureheight}{0.2\textwidth}
        \begin{tikzpicture}

\begin{axis}[width=0.993\figurewidth,
height=\figureheight,
at={(0\figurewidth,0\figureheight)},
scale only axis,
xmin=-0.06,
xmax=0.06,
xtick={-0.05,     0,  0.05},
scaled x ticks={false},
x tick label style={
        /pgf/number format/.cd,
            fixed,
            precision=2,
        /tikz/.cd},
xlabel style={font=\color{white!15!black}},
xlabel={$u_s$},
ymin=0.7,
ymax=1,
ytick={0.7, 1},
ylabel style={font=\color{white!15!black}},
ylabel={$y$},
axis background/.style={fill=white},
legend style={at={(1.1,0.5)}, anchor=west, legend cell align=left, align=left, fill=none, draw=none}
]
\addplot [color=black!60!blue, line width=1.0pt]
  table[row sep=crcr]{0.004873017106645	0.7\\
0.00283144679977254	0.706\\
0.000260704061287028	0.712\\
-0.00266242521152174	0.718\\
-0.00573903770037263	0.724\\
-0.00876160320528981	0.73\\
-0.0115279208714771	0.736\\
-0.0138546053577002	0.742\\
-0.0155892060568327	0.748\\
-0.0166201608217608	0.754\\
-0.0168839364663006	0.76\\
-0.0163689011494832	0.766\\
-0.0151156954418737	0.772\\
-0.0132141044214244	0.778\\
-0.0107966667492037	0.784\\
-0.00802947277414564	0.79\\
-0.00510078803296106	0.796\\
-0.00220827899042884	0.802\\
0.000454294541131359	0.808\\
0.0027100266142049	0.814\\
0.00441058892444398	0.82\\
0.00544595611363335	0.826\\
0.00575149669958148	0.832\\
0.0053119743518635	0.838\\
0.00416221397669016	0.844\\
0.00238441174380331	0.85\\
0.000102292682715442	0.856\\
-0.00252747119089452	0.862\\
-0.00532598895778141	0.868\\
-0.00810429240101521	0.874\\
-0.0106760689343567	0.88\\
-0.0128701296904681	0.886\\
-0.0145417800956246	0.892\\
-0.0155823392484025	0.898\\
-0.0159261832697489	0.904\\
-0.0155548575660573	0.91\\
-0.0144980020354819	0.916\\
-0.0128310481575666	0.922\\
-0.0106698631735353	0.928\\
-0.00816271972644767	0.934\\
-0.00548014589056511	0.94\\
-0.00280334884308123	0.946\\
-0.000311996511300749	0.952\\
0.00182782043062438	0.958\\
0.00347464105427542	0.964\\
0.00452096697366287	0.97\\
0.00490026092441308	0.976\\
0.00459115665612841	0.982\\
0.0036186151607354	0.988\\
0.00205196834002717	0.994\\
1.94086425377319e-19	1\\
};
\addlegendentry{$0.000~s$}

\addplot [color=black!20!blue, line width=1.0pt]
  table[row sep=crcr]{0.00864956638213382	0.7\\
0.00415004602122081	0.706\\
-0.000789727851181441	0.712\\
-0.00583491385348512	0.718\\
-0.0106464158920077	0.724\\
-0.0149036043588178	0.73\\
-0.0183256162598766	0.736\\
-0.020689823443815	0.742\\
-0.0218462583459419	0.748\\
-0.021727066032825	0.754\\
-0.0203503906067497	0.76\\
-0.0178184801416983	0.766\\
-0.0143101819297411	0.772\\
-0.0100683729999164	0.778\\
-0.00538320486534863	0.784\\
-0.000572314286431753	0.79\\
0.00404065422785672	0.796\\
0.00814776577227186	0.802\\
0.011477298585117	0.808\\
0.0138116719224365	0.814\\
0.0150016184722564	0.82\\
0.0149756810263868	0.826\\
0.013744432225684	0.832\\
0.0113991731650235	0.838\\
0.00810523739747514	0.844\\
0.00409038657150821	0.85\\
-0.000370891362924525	0.856\\
-0.00497610216952864	0.862\\
-0.00941521003143745	0.868\\
-0.0133914655222012	0.874\\
-0.0166412005691114	0.88\\
-0.0189512727030087	0.886\\
-0.0201730069997607	0.892\\
-0.0202317260975801	0.898\\
-0.0191312598564163	0.904\\
-0.0169531661011136	0.91\\
-0.0138507494977907	0.916\\
-0.01003831308845	0.922\\
-0.00577639318294701	0.928\\
-0.00135399224844817	0.934\\
0.00293098116620824	0.94\\
0.00679174413651696	0.946\\
0.00997179815934679	0.952\\
0.0122619373862481	0.958\\
0.0135139805821259	0.964\\
0.0136503274179653	0.97\\
0.0126687248119924	0.976\\
0.0106419535007936	0.982\\
0.00771248691136147	0.988\\
0.00408251074406713	0.994\\
9.43933461912473e-20	1\\
};
\addlegendentry{$0.156~s$}

\addplot [color=black!60!blue, dashed, line width=1.0pt]
  table[row sep=crcr]{0.00473786611423621	0.7\\
0.0017798372703309	0.706\\
-0.00113820263389291	0.712\\
-0.00382098368693833	0.718\\
-0.00609062581552001	0.724\\
-0.00779839468763093	0.73\\
-0.00883441291997777	0.736\\
-0.00913469470592924	0.742\\
-0.00868505562289004	0.748\\
-0.00752166145561698	0.754\\
-0.00572820603656305	0.76\\
-0.00342993317927377	0.766\\
-0.000784926776644598	0.772\\
0.00202672775404067	0.778\\
0.00481516998633569	0.784\\
0.00739355121119897	0.79\\
0.00959052245480506	0.796\\
0.0112615912478093	0.802\\
0.0122985963576066	0.808\\
0.0126366809624718	0.814\\
0.0122583164161613	0.82\\
0.0111941290812806	0.826\\
0.00952049812910925	0.832\\
0.00735410823118227	0.838\\
0.00484384328249875	0.844\\
0.00216058229754282	0.85\\
-0.000514405084565314	0.856\\
-0.00300167729281422	0.862\\
-0.00513563190440698	0.868\\
-0.00677550684027142	0.874\\
-0.00781464242145764	0.88\\
-0.00818739629278903	0.886\\
-0.00787326430350166	0.892\\
-0.00689795031644021	0.898\\
-0.00533133372302665	0.904\\
-0.00328249131637419	0.91\\
-0.00089212621648549	0.916\\
0.00167707227421881	0.922\\
0.00425147896020991	0.928\\
0.0066582441484935	0.934\\
0.00873687766953786	0.94\\
0.0103499305881497	0.946\\
0.0113920649653262	0.952\\
0.0117969143912749	0.958\\
0.0115412896277928	0.964\\
0.0106464641951228	0.97\\
0.0091764718673066	0.976\\
0.0072335485163109	0.982\\
0.00495104114027605	0.988\\
0.00248427455074564	0.994\\
-8.9202158699377e-20	1\\
};
\addlegendentry{$0.312~s$}

\addplot [color=black!20!blue, dashed, line width=1.0pt]
  table[row sep=crcr]{-0.00338513161993004	0.7\\
-0.00217239680719139	0.706\\
-0.000474975153889381	0.712\\
0.00158926425885227	0.718\\
0.00387887508488039	0.724\\
0.00623849245124593	0.73\\
0.00850934257078974	0.736\\
0.0105398945111429	0.742\\
0.0121959415935475	0.748\\
0.0133694436176416	0.754\\
0.0139855490817144	0.76\\
0.0140073425903967	0.766\\
0.0134380180253335	0.772\\
0.0123203522210758	0.778\\
0.0107335350880237	0.784\\
0.0087875882249281	0.79\\
0.00661576336213483	0.796\\
0.00436544397516196	0.802\\
0.00218816950330211	0.808\\
0.000229455651393141	0.814\\
-0.00138090705674687	0.82\\
-0.00253743122210971	0.826\\
-0.00316582149758013	0.832\\
-0.00322772591772485	0.838\\
-0.0027230471151386	0.844\\
-0.00168967615050042	0.85\\
-0.000200684942472314	0.856\\
0.00164081706349174	0.862\\
0.00370880160376886	0.868\\
0.00586292569405024	0.874\\
0.00795803514700325	0.88\\
0.00985392537059146	0.886\\
0.0114247044323181	0.892\\
0.0125671343682555	0.898\\
0.0132073992188019	0.904\\
0.013305857170068	0.91\\
0.0128594719578392	0.916\\
0.011901775956807	0.922\\
0.0105003834961505	0.928\\
0.00875223675460804	0.934\\
0.00677691715924958	0.94\\
0.00470848256067314	0.946\\
0.00268638621537468	0.952\\
0.000846091568297669	0.958\\
-0.000689986646406859	0.964\\
-0.00182061023052682	0.97\\
-0.00247237558475771	0.976\\
-0.00260446503126268	0.982\\
-0.00221118475137587	0.988\\
-0.00132213276198229	0.994\\
-1.93509474380176e-19	1\\
};
\addlegendentry{$0.469~s$}

\addplot [color=black!60!blue, dashdotted, line width=1.0pt]
  table[row sep=crcr]{-0.00849922284000933	0.7\\
-0.00419367527109562	0.706\\
0.000610438520043223	0.712\\
0.00558687951617889	0.718\\
0.0104006008850419	0.724\\
0.0147302360972904	0.73\\
0.0182894877897723	0.736\\
0.0208459980050446	0.742\\
0.0222364598489313	0.748\\
0.0223769912676082	0.754\\
0.0212681155610334	0.76\\
0.018994058353138	0.766\\
0.0157164523879569	0.772\\
0.0116629141546526	0.778\\
0.00711129519361908	0.784\\
0.00237069366440827	0.79\\
-0.00223948006939747	0.796\\
-0.00641096979477977	0.802\\
-0.00986723267593474	0.808\\
-0.0123817235030475	0.814\\
-0.0137926981267519	0.82\\
-0.014013571591958	0.826\\
-0.0130381705033274	0.832\\
-0.0109405651116499	0.838\\
-0.00786953112306057	0.844\\
-0.00403804984286519	0.85\\
0.000291415924059576	0.856\\
0.00482485552942752	0.862\\
0.00925663144686567	0.868\\
0.0132900408283116	0.874\\
0.0166571373034552	0.88\\
0.0191364915213836	0.886\\
0.0205677157107877	0.892\\
0.0208618018551609	0.898\\
0.020006609446172	0.904\\
0.0180671676533746	0.91\\
0.015180805880737	0.916\\
0.0115474726091219	0.922\\
0.00741592205129316	0.928\\
0.00306672027790714	0.934\\
-0.0012067707789002	0.94\\
-0.0051181450813459	0.946\\
-0.0084071125340135	0.952\\
-0.0108567878117648	0.958\\
-0.0123079617116421	0.964\\
-0.0126694178951461	0.97\\
-0.0119236284267782	0.976\\
-0.0101274750049309	0.982\\
-0.00740797799541861	0.988\\
-0.00395334976406034	0.994\\
-1.25814048846413e-19	1\\
};
\addlegendentry{$0.625~s$}

\addplot [color=black!20!blue, dashdotted, line width=1.0pt]
  table[row sep=crcr]{-0.00605869880748897	0.7\\
-0.00248736548795088	0.706\\
0.00115325809553849	0.712\\
0.00461854365045963	0.718\\
0.00767765342961279	0.724\\
0.0101288689507648	0.73\\
0.0118128474155562	0.736\\
0.0126229374272347	0.742\\
0.0125118887660561	0.748\\
0.0114945368880037	0.754\\
0.00964631475454788	0.76\\
0.00709772426008487	0.766\\
0.00402516820550391	0.772\\
0.000638783648100173	0.778\\
-0.00283188723244341	0.784\\
-0.00615341456182131	0.79\\
-0.00910414029013065	0.796\\
-0.0114889319406968	0.802\\
-0.0131520510173575	0.808\\
-0.0139872896749379	0.814\\
-0.0139447179677503	0.82\\
-0.0130336152474523	0.826\\
-0.0113214173517857	0.832\\
-0.00892877869916607	0.838\\
-0.00602110736444895	0.844\\
-0.00279716443379035	0.85\\
0.000524488968143117	0.856\\
0.00372027515804818	0.862\\
0.00657661617605461	0.868\\
0.00890417646560001	0.874\\
0.010550384159237	0.88\\
0.011409404736834	0.886\\
0.0114289167879285	0.892\\
0.0106132578675056	0.898\\
0.00902275412508204	0.904\\
0.00676930391631114	0.91\\
0.00400853576333364	0.916\\
0.000929088139667677	0.922\\
-0.00226025241196623	0.928\\
-0.00534467975580243	0.934\\
-0.00811780900491925	0.94\\
-0.0103954606716161	0.946\\
-0.0120278691544625	0.952\\
-0.0129095078371507	0.958\\
-0.0129858876658598	0.964\\
-0.0122568926739313	0.97\\
-0.010776450462251	0.976\\
-0.00864858226551666	0.982\\
-0.00602011937086175	0.988\\
-0.0030705941372717	0.994\\
5.3712393510693e-20	1\\
};
\addlegendentry{$0.781~s$}

\addplot [color=black!60!blue, dotted, line width=1.0pt]
  table[row sep=crcr]{0.00176715742346486	0.7\\
0.00142986282360404	0.706\\
0.000670993217696901	0.712\\
-0.000455028771984766	0.718\\
-0.00186964947517441	0.724\\
-0.00347563993088752	0.73\\
-0.00516375500720117	0.736\\
-0.00682014142936335	0.742\\
-0.00833399393438502	0.748\\
-0.00960494619296689	0.754\\
-0.0105497048891042	0.76\\
-0.0111074897109695	0.766\\
-0.0112439251122071	0.772\\
-0.0109531357942044	0.778\\
-0.0102579196928467	0.784\\
-0.00920800158836282	0.79\\
-0.0078764986154646	0.796\\
-0.0063548473960988	0.802\\
-0.00474654342086327	0.808\\
-0.00316012006498836	0.814\\
-0.00170184229471943	0.82\\
-0.000468605728271786	0.826\\
0.000458485550839918	0.832\\
0.00101943778669758	0.838\\
0.00117923508005468	0.844\\
0.0009300072503151	0.85\\
0.000291364992435389	0.856\\
-0.000691107256519786	0.862\\
-0.00194908708404337	0.868\\
-0.00339625004062953	0.874\\
-0.00493417853186798	0.88\\
-0.00645904022486797	0.886\\
-0.00786858378472559	0.892\\
-0.00906898156206651	0.898\\
-0.00998106222267123	0.904\\
-0.0105455201450436	0.91\\
-0.0107267595845317	0.916\\
-0.01051512518062	0.922\\
-0.00992737996969155	0.928\\
-0.00900541023259975	0.934\\
-0.00781325530204302	0.94\\
-0.00643267193400582	0.946\\
-0.00495753960373452	0.952\\
-0.00348748874274366	0.958\\
-0.00212118356133415	0.964\\
-0.000949711542758395	0.97\\
-5.05217620957064e-05	0.976\\
0.00051768527584634	0.982\\
0.000718779752067652	0.988\\
0.000541488363355016	0.994\\
1.85496062803966e-19	1\\
};
\addlegendentry{$0.938~s$}

\end{axis}
\end{tikzpicture}    \end{minipage}
    \caption{The solid displacement for the transient nonlinear FSI case in three dimensions
        with fluid and solid density $\rho_f = 2.1$ and $\rho_s = 5$, fluid viscosity $\mu_f = 0.03$,
        solid stiffness $\mu_s = 0.1$ and cycle length $T = 1.024$:
        By increasing the solid density compared to Figure \ref{3D-nonlinear-tf-ts-solution-along-y-fig},
        the solid model exhibits higher-order deformation modes, similar to the mode shapes of an Euler-Bernoulli beam.}
    \label{3D-nonlinear-tf-ts-solution-along-y-rhos5-fig}
\end{figure}
\FloatBarrier
\fi
\iffigure
\FloatBarrier
\begin{figure}[ht!]
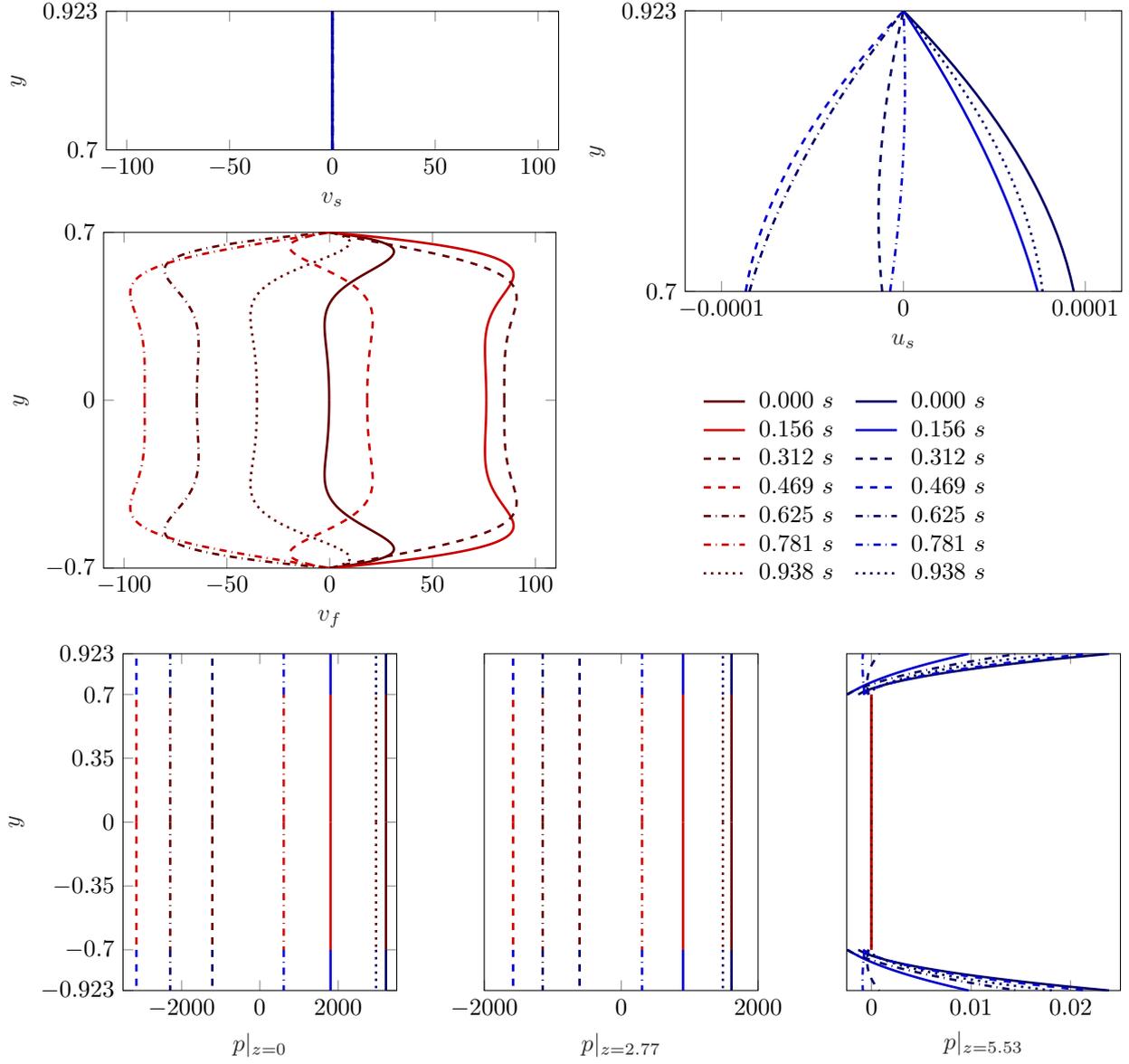

    \centering
    \begin{minipage}{\textwidth}
        \begin{minipage}{0.5\textwidth}
            \centering
            ~\\[4ex]
            \setlength{\figurewidth}{0.8\textwidth}
            \setlength{\figureheight}{0.35\textwidth}
    \end{minipage}
    \caption{The analytic solution for the transient nonlinear FSI case in three dimensions
        with fluid and solid density $\rho_f = \rho_s = 1.03~g/cm^3$, fluid viscosity $\mu_f = 0.03~g/(cm \cdot s^2)$,
        solid stiffness $\mu_s = 2 \cdot 10^5~g/(cm \cdot s^2)$, pressure amplitude $P = 583$ and cycle length $T = 1.024~s$:
        Fluid and solid velocity, $v_f$ (left) and $v_s$ (top left), and solid displacement, $u_s$ (top right),
        along the $y$-axis.
        Further, fluid and solid pressure, $p_f$ and $p_s$ (bottom) over time $t$
        at three positions along the $z$-axis.}
    \label{3D-nonlinear-tf-ts-solution-along-y-physiological-fig}
\end{figure}
\FloatBarrier
\FloatBarrier
\begin{figure}[ht!]
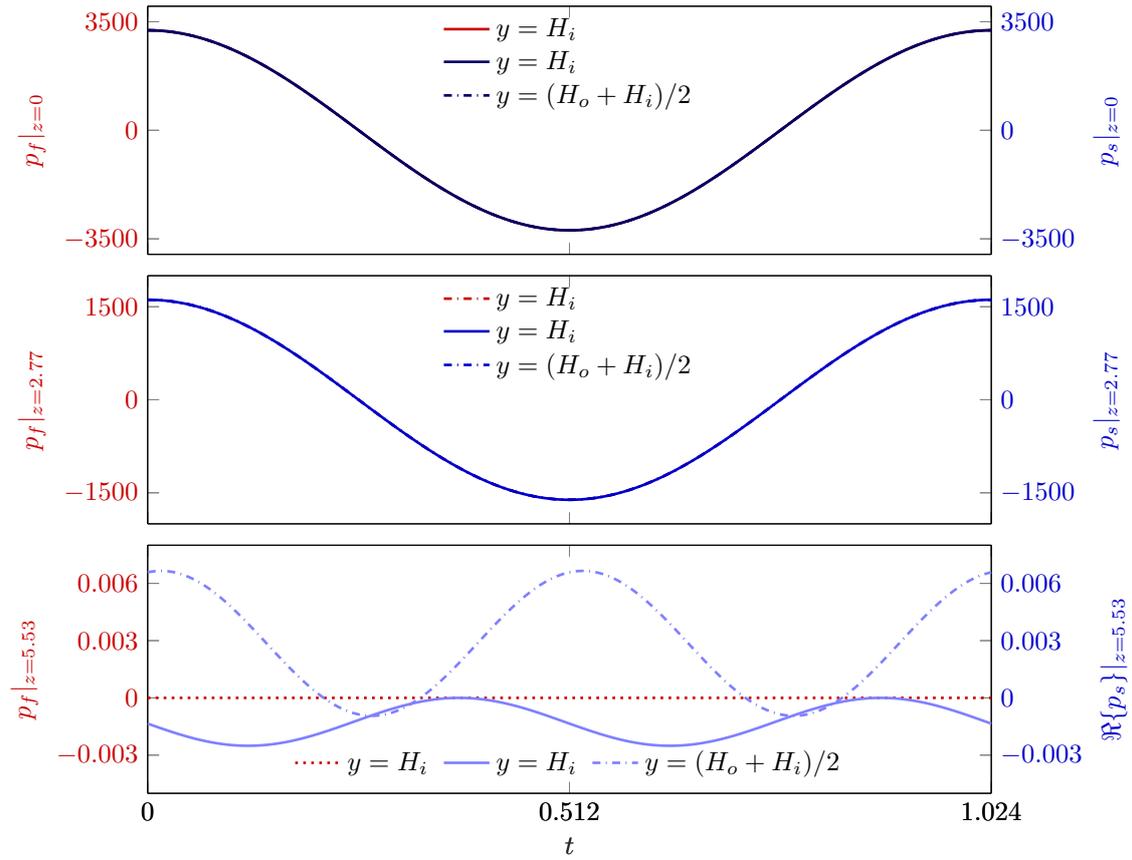

    \centering
    \setlength{\figurewidth}{0.825\textwidth}
    \setlength{\figureheight}{0.2\textwidth}
    \caption{The analytic solution for the transient nonlinear FSI case in three dimensions
        with fluid and solid density $\rho_f = \rho_s = 1.03~g/cm^3$, fluid viscosity $\mu_f = 0.03~g/(cm \cdot s^2)$,
        solid stiffness $\mu_s = 2 \cdot 10^5~g/(cm \cdot s^2)$, pressure amplitude $P = 583$ and cycle length $T = 1.024~s$:
        The fluid and solid pressure, $p_f$ and $p_s$, over time $t$ at
        three positions along the $z$-axis: At the inlet, at the mid-way point and at the outlet (top to bottom).}
    \label{3D-nonlinear-tf-ts-solution-along-t-physiological-fig}
\end{figure}
\FloatBarrier
\fi
\iffigure
\setlength{\figurewidth}{0.4\textwidth}
\setlength{\figureheight}{0.2\textwidth}
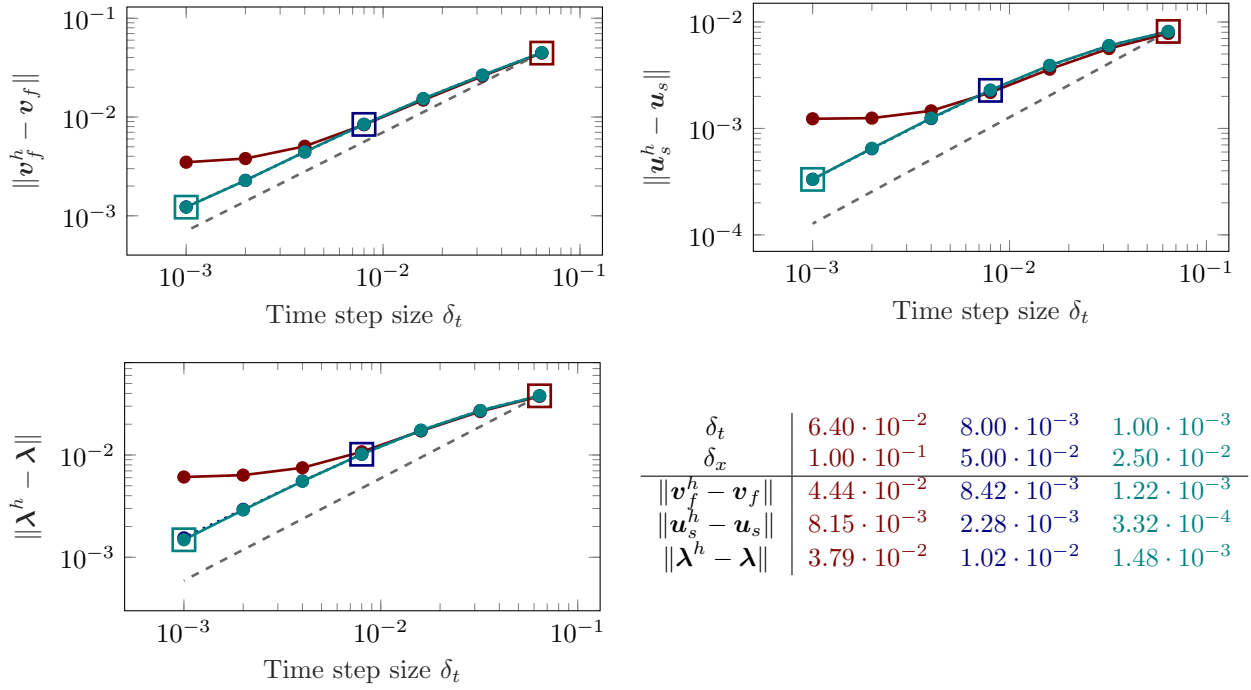
\begin{figure}[ht!]
    \centering
    \begin{minipage}{\textwidth}
        \begin{minipage}{0.5\textwidth}
            \begin{tikzpicture}

\definecolor{myred}{RGB}{128,0,0}
\definecolor{myblue}{RGB}{0,0,128}
\definecolor{mygreen}{RGB}{0,128,128}

\begin{axis}[width=0.951\figurewidth,
height=\figureheight,
at={(0\figurewidth,0\figureheight)},
scale only axis,
clip=false,
xmode=log,
xmin=0.0005,
xmax=0.13,
xtick={0.001, 0.01, 0.1},
xminorticks=true,
xlabel style={font=\color{white!15!black}},
xlabel={Time step size $\delta_t$},
ymode=log,
ymin=0.0004,
ymax=0.13,
ytick={0.001,  0.01, 0.1},
yminorticks=true,
ylabel style={font=\color{white!15!black}},
ylabel={$\| \boldsymbol{v}_f^h - \boldsymbol{v}_f \|$},
axis background/.style={fill=white}
]

\addplot [color=myred, line width=1.0pt, mark size=2pt, mark=*, mark options={solid, myred},  forget plot]
  table[row sep=crcr]{0.064                0.04437813\\
0.032                0.02594225\\
0.016                0.01484448\\
0.008                0.00834862\\
0.004                0.00505502\\
0.002                0.00379728\\
0.001                0.00348791\\
};
\addplot [color=myblue, dotted, line width=1.0pt, mark size=2pt, mark=*, mark options={solid, myblue},   forget plot]
  table[row sep=crcr]{0.064                0.04481856\\
0.032                0.02650512\\
0.016                0.01530853\\
0.008                0.00841845\\
0.004                0.00442856\\
0.002                0.00228144\\
0.001                0.00122914\\
};
\addplot [color=mygreen, line width=1.0pt,  mark size=2pt, mark=*, mark options={solid, mygreen},  forget plot]
  table[row sep=crcr]{0.064                0.04485248\\
0.032                 0.0265458\\
0.016                0.01534467\\
0.008                0.00844193\\
0.004                0.00443887\\
0.002                0.00228176\\
0.001                0.00122298\\
};
\addplot [color=myred, line width=1.0pt, draw=none, mark size=4.2pt, mark=square, mark options={solid, myred}, forget plot]
  table[row sep=crcr]{0.064                0.04437813\\
};
\addplot [color=myblue, line width=1.0pt, draw=none, mark size=4.2pt, mark=square, mark options={solid, myblue}, forget plot]
  table[row sep=crcr]{0.008                0.00841845\\
};
\addplot [color=mygreen, line width=1.0pt, draw=none, mark size=4.2pt, mark=square, mark options={solid, mygreen}, forget plot]
  table[row sep=crcr]{0.001                0.00122298\\
};
\addplot [color=white!40!black, line width=1.0pt, dashed, forget plot]
  table[row sep=crcr]{0.064   0.04485248\\
0.001   0.00070082\\
};

\end{axis}
\end{tikzpicture}        \end{minipage}
        \begin{minipage}{0.5\textwidth}
            \begin{tikzpicture}

\definecolor{myred}{RGB}{128,0,0}
\definecolor{myblue}{RGB}{0,0,128}
\definecolor{mygreen}{RGB}{0,128,128}

\begin{axis}[width=0.951\figurewidth,
height=\figureheight,
at={(0\figurewidth,0\figureheight)},
scale only axis,
clip=false,
xmode=log,
xmin=0.0005,
xmax=0.13,
xtick={0.001, 0.01, 0.1},
xminorticks=true,
xlabel style={font=\color{white!15!black}},
xlabel={Time step size $\delta_t$},
ymode=log,
ymin=0.00007,
ymax=0.015,
ytick={0.0001, 0.001,  0.01},
yminorticks=true,
ylabel style={font=\color{white!15!black}},
ylabel={$\| \boldsymbol{u}_s^h - \boldsymbol{u}_s \|$},
axis background/.style={fill=white}
]

\addplot [color=myred, line width=1.0pt, mark size=2pt, mark=*, mark options={solid, myred},  forget plot]
  table[row sep=crcr]{0.064                0.00779916\\
0.032                0.00563051\\
0.016                0.00360224\\
0.008                0.00217709\\
0.004                 0.0014629\\
0.002                0.00124972\\
0.001                0.00123215\\
};
\addplot [color=myblue, dotted, line width=1.0pt, mark size=2pt, mark=*, mark options={solid, myblue},   forget plot]
  table[row sep=crcr]{0.064                0.00813278\\
0.032                  0.005985\\
0.016                0.00389259\\
0.008                0.00227774\\
0.004                0.00123853\\
0.002                0.00064504\\
0.001                0.00033285\\
};
\addplot [color=mygreen, line width=1.0pt,  mark size=2pt, mark=*, mark options={solid, mygreen},  forget plot]
  table[row sep=crcr]{0.064                0.00815417\\
0.032                0.00600733\\
0.016                0.00391282\\
0.008                0.00229327\\
0.004                0.00124885\\
0.002                0.00065028\\
0.001                0.00033199\\
};
\addplot [color=myred, line width=1.0pt, draw=none, mark size=4.2pt, mark=square, mark options={solid, myred}, forget plot]
  table[row sep=crcr]{0.064                0.00815417\\
};
\addplot [color=myblue, line width=1.0pt, draw=none, mark size=4.2pt, mark=square, mark options={solid, myblue}, forget plot]
  table[row sep=crcr]{0.008                0.00227774\\
};
\addplot [color=mygreen, line width=1.0pt, draw=none, mark size=4.2pt, mark=square, mark options={solid, mygreen}, forget plot]
  table[row sep=crcr]{0.001                0.00033199\\
};
\addplot [color=white!40!black, line width=1.0pt, dashed, forget plot]
  table[row sep=crcr]{0.064                0.00815417\\
0.001                0.00012740890625\\
};

\end{axis}
\end{tikzpicture}        \end{minipage}
    \end{minipage}\\[2ex]
    \begin{minipage}{\textwidth}
        \begin{minipage}{0.5\textwidth}
            \begin{tikzpicture}

\definecolor{myred}{RGB}{128,0,0}
\definecolor{myblue}{RGB}{0,0,128}
\definecolor{mygreen}{RGB}{0,128,128}

\begin{axis}[width=0.951\figurewidth,
height=\figureheight,
at={(0\figurewidth,0\figureheight)},
scale only axis,
clip=false,
xmode=log,
xmin=0.0005,
xmax=0.13,
xtick={0.001, 0.01, 0.1},
xminorticks=true,
xlabel style={font=\color{white!15!black}},
xlabel={Time step size $\delta_t$},
ymode=log,
ymin=0.0003,
ymax=0.08,
ytick={0.001,  0.01},
yminorticks=true,
ylabel style={font=\color{white!15!black}},
ylabel={$\| \boldsymbol{\lambda}^h - \boldsymbol{\lambda} \|$},
axis background/.style={fill=white}
]
\addplot [color=myred, line width=1.0pt, mark size=2pt, mark=*, mark options={solid, myred},  forget plot]
  table[row sep=crcr]{0.064                0.03744063\\
0.032                0.02661045\\
0.016                0.01717316\\
0.008                0.0107552\\
0.004                0.00749908\\
0.002                0.00635519\\
0.001                0.00609329\\
};
\addplot [color=myblue, dotted, line width=1.0pt, mark size=2pt, mark=*, mark options={solid, myblue},   forget plot]
  table[row sep=crcr]{0.064                0.03787704\\
0.032                0.02715847\\
0.016                0.01743042\\
0.008                0.01016763\\
0.004                0.00556045\\
0.002                0.00293973\\
0.001                0.00154865\\
};
\addplot [color=mygreen, line width=1.0pt,  mark size=2pt, mark=*, mark options={solid, mygreen},  forget plot]
  table[row sep=crcr]{0.064                0.03790284\\
0.032                0.02719056\\
0.016                0.01745377\\
0.008                0.01017091\\
0.004                0.00554147\\
0.002                0.00289914\\
0.001                0.0014832\\
};
\addplot [color=myred, line width=1.0pt, draw=none, mark size=4.2pt, mark=square, mark options={solid, myred}, forget plot]
  table[row sep=crcr]{0.064                0.03790284\\
};
\addplot [color=myblue, line width=1.0pt, draw=none, mark size=4.2pt, mark=square, mark options={solid, myblue}, forget plot]
  table[row sep=crcr]{0.008                0.01016763\\
};
\addplot [color=mygreen, line width=1.0pt, draw=none, mark size=4.2pt, mark=square, mark options={solid, mygreen}, forget plot]
  table[row sep=crcr]{0.001                0.0014832\\
};
\addplot [color=white!40!black, line width=1.0pt, dashed, forget plot]
  table[row sep=crcr]{0.064                0.03790284\\
0.001                0.000592231875\\
};

\end{axis}
\end{tikzpicture}        \end{minipage}
        \begin{minipage}{0.5\textwidth}
            \centering
            \vspace{-25pt}
                        \resizebox{\columnwidth}{!}{
            \begin{tabular}{ c | c  c  c }                                 $\delta_t$ & \color{myred}$6.40 \cdot 10^{-2}$ & \color{myblue}$8.00 \cdot 10^{-3}$ & \color{mygreen}$1.00 \cdot 10^{-3}$ \\
                $\delta_x$ & \color{myred}$1.00 \cdot 10^{-1}$ & \color{myblue}$5.00 \cdot 10^{-2}$ & \color{mygreen}$2.50 \cdot 10^{-2}$ \\ \hline
                                                                $\| \boldsymbol{v}_f^h - \boldsymbol{v}_f \|$           & \color{myred}$4.44 \cdot 10^{-2}$ & \color{myblue}$8.42 \cdot 10^{-3}$ & \color{mygreen}$1.22 \cdot 10^{-3}$ \\
                $\| \boldsymbol{u}_s^h - \boldsymbol{u}_s \|$           & \color{myred}$8.15 \cdot 10^{-3}$ & \color{myblue}$2.28 \cdot 10^{-3}$ & \color{mygreen}$3.32 \cdot 10^{-4}$ \\
                $\| \boldsymbol{\lambda}^h - \boldsymbol{\lambda} \|$   & \color{myred}$3.79 \cdot 10^{-2}$ & \color{myblue}$1.02 \cdot 10^{-2}$ & \color{mygreen}$1.48 \cdot 10^{-3}$
            \end{tabular}
            }
        \end{minipage}
    \end{minipage}
    \caption{Spatiotemporal convergence for the transient linear FSI case in two dimensions:
        Temporal convergence of error for cycle $10$ (i.e., $\Omega_t = [9 T, 10 T]$)
        for fluid velocity $\boldsymbol{v}_f$, solid displacement $\boldsymbol{u}_s$
        and Lagrange multiplier $\boldsymbol{\lambda}$
        (and corresponding numerical approximations $\boldsymbol{v}_f^h$, $\boldsymbol{u}_s^h$
        and $\boldsymbol{\lambda}^h$)
        for \emph{coarse} (red), \emph{medium} (blue) and \emph{fine} (green) spatial resolutions.
        Square markers highlight optimal convergence of error
        for $\delta_t / \delta_x^3 = \text{const}$ (see Table).
        Dashed line illustrates optimal rate.
        All norms represent the $L^2 (\Omega_t; L^2 (\Omega_k^0))$ norm for the $k$-domain (as appropriate for the variable).}
    \label{2D-linear-tf-ts-t-conv-vf-us-lambda-fig}
\end{figure}
\FloatBarrier
\fi
\FloatBarrier
\iffigure
\setlength{\figurewidth}{0.4\textwidth}
\setlength{\figureheight}{0.2\textwidth}
\begin{figure}[ht!]
    \centering
    \begin{minipage}{\textwidth}
        \begin{minipage}{0.5\textwidth}
            \begin{tikzpicture}

\definecolor{myred}{RGB}{128,0,0}
\definecolor{myblue}{RGB}{0,0,128}
\definecolor{mygreen}{RGB}{0,128,128}

\begin{axis}[width=0.951\figurewidth,
height=\figureheight,
at={(0\figurewidth,0\figureheight)},
scale only axis,
clip=false,
xmode=log,
xmin=0.0005,
xmax=0.05,
xminorticks=true,
xlabel style={font=\color{white!15!black}},
xlabel={Time step size $\delta_t$},
ymode=log,
ymin=0.00005,
ymax=0.01,
ytick={0.00001, 0.0001, 0.001, 0.01},
yminorticks=true,
ylabel style={font=\color{white!15!black}},
ylabel={$\| \boldsymbol{v}_f^h - \boldsymbol{v}_f \|$},
axis background/.style={fill=white}
]

\addplot [color=myred, line width=1.0pt, mark size=2pt, mark=*, mark options={solid, myred},  forget plot]
  table[row sep=crcr]{0.04096                0.00621589\\
0.02048                0.00310841\\
0.01024                0.00155029\\
0.00512                0.00090497\\
0.00256                0.00071823\\
0.00128                 0.0006887\\
0.00064                0.00069171\\
};
\addplot [color=myblue, dotted, line width=1.0pt, mark size=2pt, mark=*, mark options={solid, myblue},   forget plot]
  table[row sep=crcr]{0.04096                0.00630349\\
0.02048                0.00319687\\
0.01024                 0.0015585\\
0.00512                0.00075433\\
0.00256                0.00037575\\
0.00128                0.00020699\\
0.00064                0.00014219\\
};
\addplot [color=mygreen, line width=1.0pt,  mark size=2pt, mark=*, mark options={solid, mygreen},  forget plot]
  table[row sep=crcr]{0.02048             0.00320530903\\
0.01024             0.00156470695\\
0.00512            0.000755246904\\
0.00256            0.000366883748\\
0.00128            0.000180832959\\
0.00064            9.17299047e-05\\
};
\addplot [color=myred, line width=1.0pt, draw=none, mark size=4.2pt, mark=square, mark options={solid, myred}, forget plot]
  table[row sep=crcr]{0.04096                0.00621589\\
};
\addplot [color=myblue, line width=1.0pt, draw=none, mark size=4.2pt, mark=square, mark options={solid, myblue}, forget plot]
  table[row sep=crcr]{0.00512                0.00075433\\
};
\addplot [color=mygreen, line width=1.0pt, draw=none, mark size=4.2pt, mark=square, mark options={solid, mygreen}, forget plot]
  table[row sep=crcr]{0.00064            9.17299047e-05\\
};
\addplot [color=white!40!black, line width=1.0pt, dashed, forget plot]
  table[row sep=crcr]{0.04096                0.00621589\\
0.00064            9.712328125e-05\\
};

\end{axis}
\end{tikzpicture}        \end{minipage}
        \begin{minipage}{0.5\textwidth}
            \begin{tikzpicture}

\definecolor{myred}{RGB}{128,0,0}
\definecolor{myblue}{RGB}{0,0,128}
\definecolor{mygreen}{RGB}{0,128,128}

\begin{axis}[width=0.951\figurewidth,
height=\figureheight,
at={(0\figurewidth,0\figureheight)},
scale only axis,
clip=false,
xmode=log,
xmin=0.0005,
xmax=0.05,
xminorticks=true,
xlabel style={font=\color{white!15!black}},
xlabel={Time step size $\delta_t$},
ymode=log,
ymin=0.00005,
ymax=0.01,
ytick={0.00001, 0.0001, 0.001,  0.01},
yminorticks=true,
ylabel style={font=\color{white!15!black}},
ylabel={$\| \boldsymbol{u}_s^h - \boldsymbol{u}_s \|$},
axis background/.style={fill=white}
]

\addplot [color=myred, line width=1.0pt, mark size=2pt, mark=*, mark options={solid, myred},  forget plot]
  table[row sep=crcr]{0.04096                0.00348602\\
0.02048                0.00207027\\
0.01024                0.00110255\\
0.00512                0.00055159\\
0.00256                0.00028789\\
0.00128                 0.0001967\\
0.00064                0.00018354\\
};
\addplot [color=myblue, dotted, line width=1.0pt, mark size=2pt, mark=*, mark options={solid, myblue},   forget plot]
  table[row sep=crcr]{0.04096             0.00350804122\\
0.02048             0.00211717567\\
0.01024             0.00115927728\\
0.00512            0.000598779816\\
0.00256            0.000302468795\\
0.00128            0.000152653344\\
0.00064            7.88621412e-05\\
};
\addplot [color=mygreen, line width=1.0pt,  mark size=2pt, mark=*, mark options={solid, mygreen},  forget plot]
  table[row sep=crcr]{0.02048             0.00211925421\\
0.01024             0.00116091537\\
0.00512            0.000599723633\\
0.00256            0.000302607246\\
0.00128            0.000151660833\\
0.00064            7.59693034e-05\\
};
\addplot [color=myred, line width=1.0pt, draw=none, mark size=4.2pt, mark=square, mark options={solid, myred}, forget plot]
  table[row sep=crcr]{0.04096                0.00348602\\
};
\addplot [color=myblue, line width=1.0pt, draw=none, mark size=4.2pt, mark=square, mark options={solid, myblue}, forget plot]
  table[row sep=crcr]{0.00512            0.000598779816\\
};
\addplot [color=mygreen, line width=1.0pt, draw=none, mark size=4.2pt, mark=square, mark options={solid, mygreen}, forget plot]
  table[row sep=crcr]{0.00064            7.59693034e-05\\
};
\addplot [color=white!40!black, line width=1.0pt, dashed, forget plot]
  table[row sep=crcr]{0.04096                0.00348602\\
0.00064            5.44690625e-05\\
};

\end{axis}
\end{tikzpicture}        \end{minipage}
    \end{minipage}\\[2ex]
    \begin{minipage}{\textwidth}
        \begin{minipage}{0.5\textwidth}
            \begin{tikzpicture}

\definecolor{myred}{RGB}{128,0,0}
\definecolor{myblue}{RGB}{0,0,128}
\definecolor{mygreen}{RGB}{0,128,128}

\begin{axis}[width=0.951\figurewidth,
height=\figureheight,
at={(0\figurewidth,0\figureheight)},
scale only axis,
clip=false,
xmode=log,
xmin=0.0005,
xmax=0.05,
xminorticks=true,
xlabel style={font=\color{white!15!black}},
xlabel={Time step size $\delta_t$},
ymode=log,
ymin=0.00001,
ymax=0.01,
ytick={0.00001, 0.0001, 0.001,  0.01},
yminorticks=true,
ylabel style={font=\color{white!15!black}},
ylabel={$\| p_s^h - p_s \|$},
axis background/.style={fill=white}
]
\addplot [color=myred, line width=1.0pt, mark size=2pt, mark=*, mark options={solid, myred},  forget plot]
  table[row sep=crcr]{0.04096                0.00227436\\
0.02048                0.00111384\\
0.01024                0.00055349\\
0.00512                0.00040372\\
0.00256                0.00040993\\
0.00128                 0.0004321\\
0.00064                0.00044664\\
};
\addplot [color=myblue, dotted, line width=1.0pt, mark size=2pt, mark=*, mark options={solid, myblue},   forget plot]
  table[row sep=crcr]{0.04096                0.00250079\\
0.02048                0.00129136\\
0.01024                0.00061965\\
0.00512                0.00029081\\
0.00256                0.00015282\\
0.00128                0.00011348\\
0.00064                0.00011053\\
};
\addplot [color=mygreen, line width=1.0pt,  mark size=2pt, mark=*, mark options={solid, mygreen},  forget plot]
  table[row sep=crcr]{0.02048             0.00134373505\\
0.01024            0.000659804529\\
0.00512             0.00031436234\\
0.00256            0.000149209852\\
0.00128            7.20275225e-05\\
0.00064            3.92226303e-05\\
};
\addplot [color=myred, line width=1.0pt, draw=none, mark size=4.2pt, mark=square, mark options={solid, myred}, forget plot]
  table[row sep=crcr]{0.04096                0.00227436\\
};
\addplot [color=myblue, line width=1.0pt, draw=none, mark size=4.2pt, mark=square, mark options={solid, myblue}, forget plot]
  table[row sep=crcr]{0.00512                0.00029081\\
};
\addplot [color=mygreen, line width=1.0pt, draw=none, mark size=4.2pt, mark=square, mark options={solid, mygreen}, forget plot]
  table[row sep=crcr]{0.00064            3.92226303e-05\\
};
\addplot [color=white!40!black, line width=1.0pt, dashed, forget plot]
  table[row sep=crcr]{0.04096                0.00227436\\
0.00064            3.5536875e-05\\
};

\end{axis}
\end{tikzpicture}        \end{minipage}
        \begin{minipage}{0.5\textwidth}
            \centering
                        \resizebox{\columnwidth}{!}{
            \begin{tabular}{ c | c  c  c }
                $\delta_t$                                  & {\color{myred}$4.10 \cdot 10^{-2}$} & {\color{myblue}$5.12 \cdot 10^{-3}$} & {\color{mygreen}$6.40 \cdot 10^{-4}$} \\
                $\delta_x, \delta_y, \delta_z$              & {\color{myred}$3.14 \cdot 10^{-1}$} & {\color{myblue}$1.57 \cdot 10^{-1}$} & {\color{mygreen}$7.90 \cdot 10^{-2}$} \\ \hline
                                                                $\| \boldsymbol{v}_f^h - \boldsymbol{v}_f \|_{\Omega_f^O \times \Omega_t, 2}$  & \color{myred}$6.22 \cdot 10^{-3}$ & \color{myblue}$7.54 \cdot 10^{-4}$ & \color{mygreen}$9.70 \cdot 10^{-5}$ \\
                $\| \boldsymbol{u}_s^h - \boldsymbol{u}_s \|_{\Omega_s^O \times \Omega_t, 2}$  & \color{myred}$3.49 \cdot 10^{-3}$ & \color{myblue}$5.99 \cdot 10^{-4}$ & \color{mygreen}$7.60 \cdot 10^{-5}$ \\
                $\| p_s^h - p_s \|_{\Omega_s^O \times \Omega_t, 2}$                            & \color{myred}$2.27 \cdot 10^{-3}$ & \color{myblue}$2.91 \cdot 10^{-4}$ & \color{mygreen}$3.90 \cdot 10^{-5}$ \\
            \end{tabular}
            }
        \end{minipage}
    \end{minipage}\\[2ex]
    \begin{minipage}{\textwidth}
        \begin{minipage}{0.5\textwidth}
                    \end{minipage}
        \begin{minipage}{0.5\textwidth}
                    \end{minipage}
    \end{minipage}
    \caption{Spatiotemporal convergence for the transient nonlinear FSI case in three dimensions:
        Temporal convergence of error for cycle $7$ (i.e., $\Omega_t = [6 T, 7 T]$)
        for fluid velocity $\boldsymbol{v}_f$, solid displacement $\boldsymbol{u}_s$
        and solid pressure $p_s$
        (and corresponding numerical approximations $\boldsymbol{v}_f^h$, $\boldsymbol{u}_s^h$, and $p_s^h$)
        for \emph{coarse} (red), \emph{medium} (blue) and \emph{fine} (green) spatial resolutions.
        Square markers highlight optimal convergence of error
        for $\delta_t / \delta_x^3 = \text{const}$ (see Table).
        Dashed line illustrates optimal rate.
        All norms represent the $L^2 (\Omega_t; L^2 (\Omega_k^0))$ norm for the $k$-domain (as appropriate for the variable).}
    \label{3D-nonlinear-tf-ts-t-conv-vf-us-lambda-pf-ps-fig}
\end{figure}
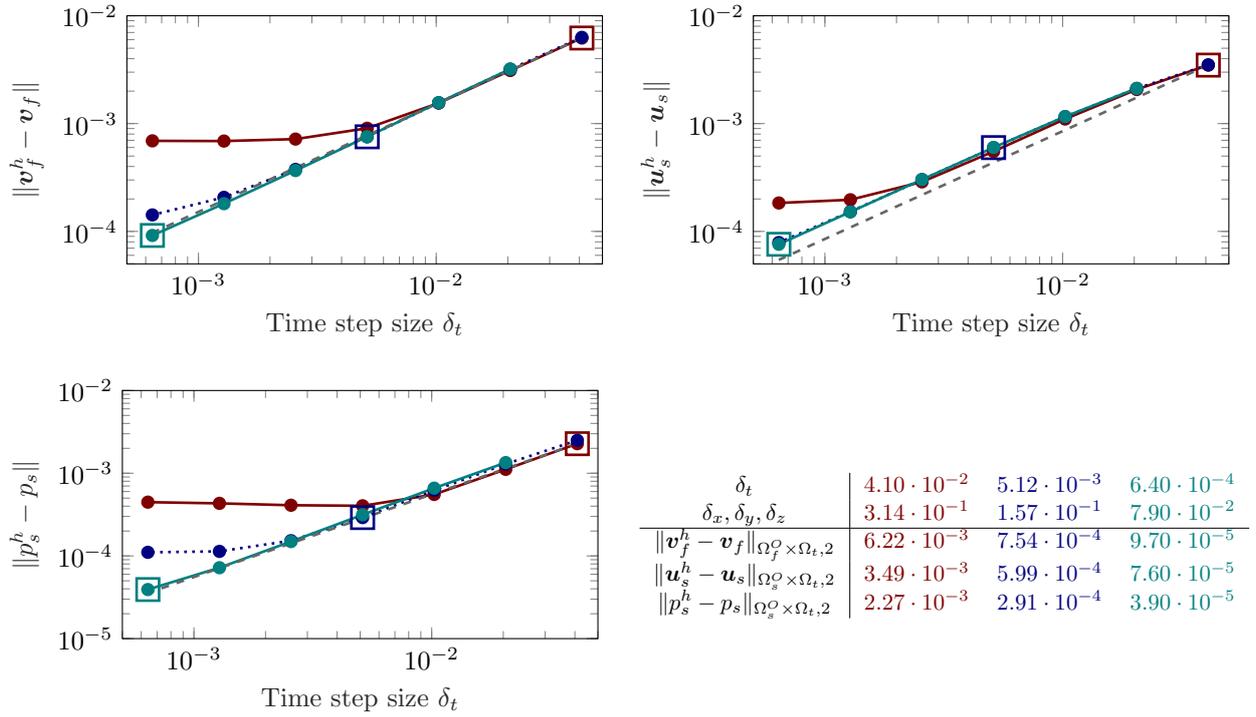
\FloatBarrier
\fi
\FloatBarrier
\clearpage

\renewcommand{\thesection}{SM\arabic{section}}
\setcounter{section}{0}\setcounter{subsection}{0}\setcounter{equation}{0}
\setcounter{figure}{0}
\setcounter{table}{0}
\setcounter{page}{1}
\renewcommand{\theequation}{SM\arabic{equation}}
\renewcommand{\thefigure}{SM\arabic{figure}}
\renewcommand{\bibnumfmt}[1]{[SM#1]}
\renewcommand{\citenumfont}[1]{SM#1}

\begin{center}
        \textbf{\large
    Supplementary Material 1:\\[1ex]
    A Class of Analytic Solutions for Verification and Convergence Analysis
    of Linear and Nonlinear Fluid-Structure Interaction Algorithms}\\[3ex]
        Andreas Hessenthaler$^{a,*}$,
    Maximilian Balmus$^{b,*}$,
    Oliver R\"ohrle$^a$,
    David Nordsletten$^{c,b}$\\[2ex]
        {\itshape \footnotesize
    ${}^a$Institute for Modelling and Simulation of Biomechanical Systems,
    University of Stuttgart,
    Pfaffenwaldring 5a, 70569 Stuttgart, Germany\\
    ${}^b$School of Biomedical Engineering and Imaging Sciences,
    King's College London, 4th FL Rayne Institute,
    St Thomas Hospital, London, SE1 7EH\\
    ${}^c$Department of Biomedical Engineering and Cardiac Surgery,
    University of Michigan,
    NCRC B20, 2800 Plymouth Rd, Ann Arbor, 48109\\
    }
    \end{center}
{\let\thefootnote\relax\footnotetext{{$^*$Authors are acknowledged as joint-first and joint-corresponding authors.}}}

\label{online-supplement-sec}

\section{Analytical derivations - continuation}
Here, we continue the set of derivations for the analytical solutions
which were left out in the main paper. The respective cases are the
2D nonlinear and the 3D linear (including all the quasi-static and
transient permutations).

\subsection{Nonlinear fluid/solid in two dimensions}
Again, we consider the nonlinear FSI problem discussed in Section~\ref{derivation-3D-sec}
and defined by Equations~\eqref{eq:strong_form_nonlinear_fluid_momentum_balance}
-\eqref{eq:strong_form_nonlinear_kinematic_constraint}),
with the mention that the definition of the first Piola-Kirchhoff stress tensor
is now adjusted for the 2D problem:
\begin{equation}
	\boldsymbol{P}_s = \frac{\mu_s}{J_s}
	\left[\boldsymbol{F} - \frac{\boldsymbol{F} : \boldsymbol{F}}{2}
	\boldsymbol{F}^{-T}\right] - J_s p_s \boldsymbol{F}^{-T}
\end{equation}
In line to Section~\ref{derivation-2D-sec}, we adapt the boundary conditions to the 2D
setting and use the same general assumption
about the behaviour of the solution. Based on these considerations, it can be shown that
the advective term, $\boldsymbol{v}_f \cdot
\nabla_{\boldsymbol{x}} \boldsymbol{v}_f$ is null and that, in effect, the fluid problem
is identical to the one in the linear case. Thus, the general solutions for
the fluid pressure is given by Equations~\eqref{eq:fluid_lin_pulse}. Similarly,
 the general quasi-static and transient fluid velocity solutions are given by
Equations~\eqref{eq:fluid_lin_pulse} and \eqref{eq:fluid_lin_general}.

To find a general solution for the solid problem, we first observe that the deformation
gradient takes the form $\boldsymbol{F} = \boldsymbol{I} +
{\partial_Y {u}_s}(\boldsymbol{e}_Y \otimes \boldsymbol{e}_X)$, where $u_s$ is the $x$ component of the deformation field.
 Based on this, we can express the first Piola-Kirchhoff
tensor in terms of ${u}_s$ and $p_s$ as follows:
\begin{equation}
	\boldsymbol{P}_s = \left[
		\begin{array}{c c}
		- \frac{\mu_s}{2} \left(\frac{\partial u_s}{\partial Y}\right)^2-p_s
		&
		\mu_s \frac{\partial u_s}{\partial Y}
		\\
		\mu_s \frac{\partial u_s}{\partial Y} + \frac{\mu_s}{2}
		\left(\frac{\partial u_s}{\partial Y}\right)^3 + p_s \frac{\partial u_s}{\partial Y}
		&
		- \frac{\mu_s}{2} \left(\frac{\partial u_s}{\partial Y}\right)^2-p_s
		\end{array}
	\right].
\end{equation}
Using this result, we expand the momentum balance equation in \eqref{eq:strong_form_nonlinear_solid_momentum_balance}
to obtain:
\begin{align}
	\left(-\frac{\partial p_s}{\partial X} + \mu_s \frac{\partial^2 u_s}{\partial Y^2}
	- \rho_s \partial_{tt} u_s
	 \right) \boldsymbol{e}_X
	+ \left[
		\frac{\partial p_s}{\partial X} \frac{\partial u_s}{\partial Y}
		- \frac{\partial p_s}{\partial Y}
		- \frac{\mu_s}{2} \frac{\partial}{\partial Y}\left(\frac{\partial u_s}{\partial Y}\right)^2
	\right] \boldsymbol{e}_Y
	= \boldsymbol{0}
	\label{eq:momentum_balance_2D_nonlin}
\end{align}
By applying $\partial_X$ to this vector, we find that $\partial^2 p_s / \partial X^2 = 0$ and
$\partial^2 p_s / (\partial X \partial Y) = 0$. Based on these results, we can integrate the $\boldsymbol{e}_y$
component of \eqref{eq:momentum_balance_2D_nonlin} with respect to $Y$ to obtain the general
solution for the pressure field:
\begin{equation}
	p_s = c_I(t) + c_{II}(t) (X + u_s) - \frac{\mu_s}{2} \left(\frac{\partial u_s}{\partial Y}\right)^2.
	\label{eq:pressure-solid-2D-nonlinear}
\end{equation}
Furthermore, assuming that $u_s = \Re\{\text{u}_s(Y) e^{i\omega t}\}$,
$\text{u}_s:[H_i,H_o]\rightarrow\mathbb{C}$
and $c_{II}(t) = \Re\{c_3e^{i\omega t}\}$, the
$X$-component of \eqref{eq:momentum_balance_2D_nonlin} can be simplified into:
\begin{equation}
	\mu_s \frac{\partial^2 \text{u}_s}{\partial Y^2} + \rho_s \omega^2 \text{u}_s - c_3 = 0.
\end{equation}
Given the fixed wall boundary condition $\text{u}_s(H_o) = 0$, the general solutions of this ODE take the form:
\begin{equation}
	\begin{aligned}
		\text{u}_s & =  \frac{c_3}{2\mu_s} (Y^2 - H_o^2) + c_4(Y - H_o), && (\rho_s = 0), \\
		\text{u}_s & = \frac{c_3}{\rho_s  \omega^2}  \left[1 - \sin(k_s Y)\csc(k_s H_o) \right]
		+ c_4 \left[ \cos (k_s Y) - \cot(k_s H_o) \sin (k_s Y) \right], && (\rho_s > 0).
		\label{eq:solid_disp_2D_nonlinear}
	\end{aligned}
\end{equation}
The last four integration constants, namely $c_1$, $c_3$ $c_4$ and $c_I$, can only be identified by
verifying that the coupling conditions are satisfied. In the following, we consider the different
permutations of quasi-static / transient behaviors and the resulting sets of closed formulations
for the constants.

\textbf{Quasi-static fluid and quasi-static solid} ($\rho_f = 0,~ \rho_s = 0$)

Let us first consider the of expansion the
traction coupling condition in Equation~\eqref{eq:strong_form_nonlinear_dynamic_constraint}:
\begin{equation}
	\left[\mu_f \left.\frac{\partial v_f}{\partial y}\right|_{y = H_i} -
	\mu_s \left.\frac{\partial u_s}{\partial Y}\right|_{Y = H_i}\right]\boldsymbol{e}_x
	+ \left[\frac{\mu_s}{2}\left. \left(\frac{\partial u_s}{\partial Y} \right)^2\right|_{Y = H_i}
	+ p_s|_{Y = H_i} - p_f|_{y=H_i}\right]\boldsymbol{e}_y = 0
\end{equation}
Taking the axial component and replacing the pressure fields with the results in
Equations~\eqref{eq:fluid_lin_pulse} and~\eqref{eq:pressure-solid-2D-nonlinear},
we obtain:
\begin{equation}
	c_I (t) + \Re\{c_3 e^{i\omega t} \} \left[x + u_s(H_i,t) \right] - \Re\{P(L - x) e^{i\omega t}\} = 0
\end{equation}
From this, it can be shown that $c_3 = -P$ and $c_I (t) = \Re\{P [L + u_s(H_i,t)] e^{i\omega t}\}$. Similarly to the
3D nonlinear case, these results are valid for all quasi-static and transient permutations due to the fact
that the general structure of the pressure solutions is independent from these factors.
Based on these results, we can expand the general nonlinear solid pressure solution in
\eqref{eq:3D_solid_nonlinear_pressure_general} to obtain a quasi-static specific version:
\begin{align}
    p_s(X,Y,t)
    =& \left[L - X + u_s(H_i,t) - u_s \right]\Re\{P e^{i\omega t} \} - \frac{\mu_s}{2}
	    \left[\frac{\partial \text{u}_s}{\partial Y} \right]^2 \nonumber \\
	=& \left[L - X +  \Re \left\{\frac{P}{2\mu_s} (Y^2 - H_i^2) e^{i\omega t} + c_4 (H_i - Y) e^{i\omega t} \right\} \right]\Re\{P e^{i\omega t} \} \nonumber \\
	& - \frac{\mu_s}{2} \left[\Re \left\{c_4e^{i\omega t} - \frac{PY}{\mu_s} e^{i\omega t} \right\}\right]^2.
	\label{eq:3D_linear_quasi_solid_pressure}
\end{align}
Furthermore, we can expand the axial component of the traction coupling condition
using Equations~\eqref{eq:fluid_lin_general} and \eqref{eq:solid_disp_2D_nonlinear} resulting in:
\begin{equation}
	\mu_f \left(-\frac{P}{2\mu_f} 2H_i  \right) -
	\mu_s \left(-\frac{P}{2\mu_s} 2H_i + c_4  \right) = 0 \quad \Rightarrow \quad c_4 = 0.
\end{equation}
Finally, using the kinematic coupling condition, it can be shown that:
\begin{equation}
	c_1 = \frac{P H_i^2}{2\mu_f} +  \frac{i\omega P}{2\mu_s} (H_o^2 - H_i^2) .
\end{equation}

\textbf{Transient fluid and quasi-static solid} ($\rho_f > 0,~\rho_s = 0$)

Similarly to the preceding case, the $x$ components of the traction and
kinematic interface conditions form a system of equations with two unknowns, $c_1$
and $c_4$, which can written as follows:
\begin{align*}
	\beta c_1 - \mu_s c_4 & = - PH_i, \\
	\alpha c_1 + i\omega(H_o - H_i)c_4 & = \frac{iP}{\rho_f \omega} + \frac{i\omega P }{2\mu_s} (H_o^2 - H_i^2),
\end{align*}
where $\alpha$ and $\beta$ retain the definitions from Table~\ref{nomenclature-tab}.
The resulting solutions are:
\begin{align}
	c_1 & = \frac{ - i\omega (H_o - H_i)H_i P + \mu_s \left[\frac{iP}{\rho_f \omega} + \frac{i\omega P }{2\mu_s} (H_o^2 - H_i^2)\right]}
	{i\omega\beta(H_o - H_i) + \mu_s \alpha }, \\
	c_4 & = \frac{\alpha H_i P + \beta \left[\frac{iP}{\rho_f \omega} + \frac{i\omega P }{2\mu_s} (H_o^2 -  H_i^2) \right] }
	{i\omega\beta(H_o - H_i) + \mu_s \alpha }.
\end{align}
The general solution for the quasi-static solid pressure (with $c_3 = - P$) given by Equation~\eqref{eq:3D_linear_quasi_solid_pressure}
applies in this case as well.

\textbf{Quasi-static fluid and transient solid} ($\rho_f = 0,~\rho_s > 0$)

Analogous to the transient fluid and quasi-static solid case, the system of equations takes the form:
\begin{align*}
\mu_s k_s\xi_1c_4 & =
\frac{\mu_s k_s \zeta_1 P}{\rho_s \omega^2}  + PH_i, \\
c_1 + i\omega \xi_2 c_4 & = \frac{PH_i^2}{2\mu_f}
- \frac{i\zeta_2P}{\rho_s \omega},
\end{align*}
where for convenience we introduced a new set of parameters:
\begin{align}
	\xi_1 & = \sin(k_s H_i)  + \cot(k_s H_o) \cos(k_s H_i), \\
	\xi_2 & = \cot(k_s H_o) \sin(k_s H_i) - \cos(k_s H_i), \\
	\zeta_1 & = \csc(k_s H_o) \cos(k_sH_i), \\
	\zeta_2 & = 1  - \sin(k_s H_i) \csc(k_s H_o).
\end{align}
Solving this system of equations gives the following closed forms for the integration constants:
\begin{align}
	c_1  &=  \frac{PH_i^2}{2\mu_f}
	- \frac{i\zeta_2 P}{\rho_s \omega}
	-  i\xi_2
	\frac{\mu_s k_s \zeta_1 P + \rho_s \omega^2PH_i }
	{\mu_s k_s\rho_s \omega\xi_1}, \\
	c_4  &= \frac{\mu_s k_s \zeta_1 P + \rho_s \omega^2PH_i }
	{\mu_s k_s\rho_s \omega^2\xi_1}.
\end{align}
The general solution for the transient solid pressure is obtained by expanding the more general form
in Equation~\eqref{eq:pressure-solid-2D-nonlinear} using $c_3 = -P$ and $c_I (t) = \Re\{P [L + u_s(H_i,t)] e^{i\omega t}\}$
and the general transient solid displacement formula in Equation~\eqref{eq:solid_disp_2D_nonlinear}:
\begin{align}
    p_s(X,Y,t)
    =& \left[L - X + u_s(H_i,t) - u_s \right] \Re\{P e^{i\omega t}\}
	    - \frac{\mu_s}{2}\left[\frac{\partial u_s}{\partial r}\right]^2 \nonumber \\
	=& \left[L - X + \Re \left\{\frac{P}{\rho_s \omega^2}
	\left(1 - \sin(k_sY) \csc(k_sH_o) - \zeta_2\right)e^{i\omega t}\right\} \right]\Re\{Pe^{i\omega t} \} \nonumber \\
	& - \Re\left\{c_4 [\cos(k_s Y) - \cot(k_s H_o) \sin(k_s Y) + \xi_2] e^{i\omega t} \right\} \Re\{P e^{i\omega t}\} \nonumber \\
	& - \frac{\mu_s}{2} \left[ \Re \left\{\frac{Pk_s}{\rho_s \omega^2} \cos(k_s Y) \csc(k_s H_o) e^{i\omega t}
	    - c_4 k_s (\sin(k_s Y) + \cot(k_sH_o) \cos(k_sY)) e^{i\omega t} \right\} \right]^2.
\end{align}
This formula also applies to the transient fluid and transient solid case provided the appropriate
$c_4$ formulation is used.

\textbf{Transient fluid and transient solid} ($\rho_f, \rho_s > 0$)

The equivalent system of equations corresponding to this case can be written as follows:
\begin{align*}
	\beta c_1 + \mu_s k_s\xi_1 c_4 & =
	\frac{\zeta_1 P }{k_s},  \\
	\alpha c_1 + i\omega\xi_2 c_4 & =
	\frac{iP}{\rho_f \omega} - \frac{i \zeta_2 P}{\rho_s \omega},
\end{align*}
and the resulting closed form solutions are:
\begin{align}
	c_1 & = \frac{iP}{ \rho_s \rho_f \omega k_s}
	\frac{\omega^2 \xi_2  \zeta_1 \rho_s \rho_f  - \mu_s  k_s^2 \xi_1  (\rho_s - \zeta_2 \rho_f)}
	{i\omega\xi_2 \beta- \alpha \mu_s k_s \xi_1}, \\
	c_4 & = \frac{P}{\rho_s \rho_f \omega k_s}
	\frac{-\omega \alpha \zeta_1 \rho_s \rho_f  + i\beta k_s (\rho_s - \zeta_2 \rho_f)}
	{i\omega\xi_2 \beta- \alpha \mu_s k_s \xi_1}.
\end{align}

\subsection{Linear fluid/solid in three dimensions}
In this section we approach the derivation of the FSI problem described in
Equations~\eqref{eq:strong_form_linear_fluid_momentum_balance} -
\eqref{eq:strong_form_linear_kinematic_constraint}, but in the 3D and
using the assumptions about the solution behaviour which were described
in Section~\ref{derivation-3D-sec}, the 3D nonlinear case.
In terms of deriving the general analytical
solutions, it should be noted that in the 3D nonlinear case,
we showed that advection term is constantly zero as a consequence
of our assumptions. Consequently, the solutions described in
Equations~\eqref{eq:3D_fluid_nonlinear_solution_periodic} and \eqref{eq:fluid_vel_3D_nonlinear}
which define the general solution for the pressure, quasi-static velocity
and transient velocity, respectively, are still applicable in this context.

We begin the derivation of the solid solution, by observing that the
the displacement gradient takes the form $\nabla\boldsymbol{u}_s
=\frac{\partial u_s}{\partial r} \boldsymbol{e}_r \otimes \boldsymbol{e}_z$ .
Based on this, the Cauchy stress tensor can be written as:
\begin{equation*}
	\boldsymbol{\sigma}_s = \mu_s u_s \left(
	\boldsymbol{e}_r \otimes \boldsymbol{e}_z +
	\boldsymbol{e}_z \otimes \boldsymbol{e}_r
	\right)  - p \boldsymbol{I}.
\end{equation*}
This can now be used to expand the solid momentum balance equation
in~(\ref{eq:strong_form_linear_solid_momentum_balance}) into the
the three cylindrical coordinates:
\begin{equation}
	 -\frac{\partial p_s}{\partial r}\boldsymbol{e}_r
	 - \frac{1}{r}
	 \frac{\partial p_s}{\partial \theta}\boldsymbol{e}_\theta
	 + \left(
	 \mu_s \frac{\partial^2 u_s}{\partial r^2}
	 +\frac{\mu_s}{r}\frac{\partial u_s}{\partial r}
	 -\frac{\partial p_s}{\partial z}
	 -\rho_s\partial_{tt} u_s
	 \right)\boldsymbol{e}_z = \boldsymbol{0}.
\end{equation}
The radial and angular components show that the pressure field is invariant
in these directions. Furthermore, by taking the axial derivative of the
$\boldsymbol{e}_z$ component, it can be shown that
$\partial^2 p_s / \partial z^2 = 0$. Base on these three observations and also taking
note of our previous assumptions (i.e. spatial and temporal components are separable, the field
is temporally periodic and described by a single harmonic), we can write
the general solid pressure solution as:
\begin{equation}
	p_s = \Re\{P_s(L-z) e^{i\omega t}\}.
\end{equation}
$P_s\in\mathbb{C}$ is the reference pressure over the length of the domain.
Substituting this into the axial component of the momentum balance
and also using the fact that $u_s = \Re\{\text{u}_s e^{i\omega t} \}$,
 we find the following non-homogeneous
PDE for the displacement:
\begin{equation}
	\mu_s \frac{\partial^2 u_s}{\partial r^2}
	+ \frac{\mu_s}{r} \frac{\partial u_s}{\partial r} + \rho_s \omega^2u_s = - P_s,
\end{equation}
which is identical to Equation~\eqref{eq:3D-nonlin-bessel-pde}, with the only difference that in
that equation non-homogenous term, i.e.~$-P_s$, was unknown.
Consequently, the
two versions of the displacement function, quasi-static and transient, that resulted in
the previous section, also apply in
this case:
\begin{equation}
	\begin{aligned}
		\text{u}_s &= \frac{P_s}{4\mu_s}(H_o^2 - r^2) + c_3 \ln \left(\frac{r}{H_o} \right),
		 && (\rho_s = 0), \\
		\text{u}_s &= -\frac{P_s}{\rho_s \omega^2} \left[1 -
		\frac{Y_0(-k_s r)}{Y_{0,s}^r}\right] + c_3 \left[J_0(-k_s r) - \gamma Y_0(-k_s r) \right], &&(\rho_s > 0).
	\end{aligned}
	\label{eq:solid_disp_3D_lin}
\end{equation}

The last two integration constants (i.e. $c_1$ and $c_3$)
are found by verifying that the coupling condition are satisfied.
In the following, we derive the closed form of these constants,
while taking into account all four quasi-static / transient permutations.

\textbf{Quasi-static fluid and quasi-static solid} ($\rho_f = \rho_s = 0$)

First, let us expand the traction interface condition into its three components:
\begin{equation}
	\left(-p_f + p_s\right)\boldsymbol{e}_r +
	\left(
	\mu_f \left.\frac{\partial v_f}{\partial r}\right|_{y=H_i} -
	\mu_s \left.\frac{\partial u_s}{\partial r}\right|_{y=H_i}
	\right)\boldsymbol{e}_z = \boldsymbol{0}.
\end{equation}
Using the radial component, it is trivial to show that $P = P_f = P_s$; a result which
applies to all quasi-static and transient permutations. In the case of the axial component,
incorporating the results of Equations~\eqref{eq:fluid_vel_3D_nonlinear} and~\ref{eq:solid_disp_3D_lin} yields:
\begin{equation}
	-\mu_f\frac{P H_i}{2\mu_f} + \mu_s\frac{P H_i}{2\mu_s} + \frac{\mu_s c_3}{H_i} =  0 \quad \Rightarrow \quad c_3 = 0.
\end{equation}
Finally, we expand the kinematic condition to obtain:
\begin{equation}
	\left(-\frac{P H_i^2}{4\mu_f} + c_1\right) - i\omega \frac{P(H_o^2 - H_i^2)}{4 \mu_s} = 0,
\end{equation}
which can be rearranged into the closed form of the last integration constant:
\begin{equation}
	c_1 = \frac{P H_i^2}{4\mu_f} + i\omega \frac{P(H_o^2 - H_i^2)}{4 \mu_s}.
\end{equation}

\textbf{Transient fluid and quasi-static solid} ($\rho_f > 0,~\rho_s = 0$)

Similarly to the previous case, the last two constants, $c_1$ and $c_3$, can be identified
by solving the system formed by the axial components of the traction and
kinematic interface conditions:
\begin{align*}
	  - i \mu_f J_{1,f}^* c_1  - \frac{\mu_s c_3}{H_i} & = - \frac{PH_i}{2} \\
	 J_{0,f}^* c_1- i\omega\ln\left(\frac{H_i}{H_o}\right) c_3 & =
	  \frac{i\omega P}{4\mu_s}(H_o^2 - H_i^2) - \frac{P}{\mu_f k_f^2} .
\end{align*}
The resulting closed forms are:
\begin{align}
	c_1 & = \frac{i \omega \ln\left(\frac{H_i}{H_o}\right) \frac{PH_i}{2}
		+ \frac{i\omega P}{4H_i}(H_o^2 - H_i^2) -
		\frac{\mu_s P}{\mu_f k_f^2 H_i}
	}{- \omega \mu_f J_{1,f}^* \ln(\frac{H_i}{H_o})
	+ J_{0,f}^*\frac{\mu_s}{H_i} } ,
	 \\
	c_3 & = \frac{ J_{0,f}^*\frac{PH_i}{2}
		+ J_{1,f}^* \frac{\mu_f \omega P}{4 \mu_s} (H_o^2 - H_i^2)
		+ J_{1,f}^* \frac{i P}{k_f^2}
	}
	{- \omega \mu_f J_{1,f}^* \ln(\frac{H_i}{H_o})
		+ J_{0,f}^*\frac{\mu_s}{H_i}}.
\end{align}

\textbf{Quasi-static fluid and transient solid} ($\rho_f = 0,~\rho_s > 0$)

Re-writing the $\boldsymbol{e}_z$ component of the traction coupling condition using
the general solutions for fluid velocity~\eqref{eq:fluid_vel_3D_nonlinear}
and solid displacement~\eqref{eq:solid_disp_3D_lin}, yields:
\begin{equation}
	-\frac{PH_i}{2} + \mu_s \left[i\frac{P\nu_1}{\rho_s \omega^2} + i\Delta_1 c_3\right] = 0.
\end{equation}
When re-arranged, we obtain the closed form of the $c_3$ constant:
\begin{equation}
	c_3 = \frac{P}{i\Delta_1} \left(\frac{H_i}{2\mu_s}
	- \frac{i \nu_1}{\rho_s\omega^2} \right).
\end{equation}
Finally, we expand the kinematic condition, which can be written as:
\begin{equation}
	-\frac{PH_i^2}{4\mu_f} + c_1 - i\omega \left[
	-(1- \nu_0)\frac{P}{\rho_s \omega^2}
	+ \Delta_0 c_3
	\right] = 0.
\end{equation}
Based on this, we find the closed form of the last unknown integration constant:
\begin{equation}
	c_1 = \frac{P H_i^2}{4 \mu_f} - i\omega P \left[
	\frac{1 - \nu_0}{\rho_s \omega^2}
	+ i \frac{\Delta_0}{\Delta_1}
	\left(\frac{H_i}{2\mu_s} -
	i \frac{\nu_1}{\rho_s \omega^2}\right)
	\right].
\end{equation}

\textbf{Transient fluid and transient solid} ($\rho_f, \rho_s > 0$)

Finally, we treat the case where both solid and fluid have transient behaviours. As in the
transient fluid~/~quasi-static solid case, $c_1$ and $c_3$ can be found by solving the system given
by the $\boldsymbol{e}_z$-components of the traction and kinematic interface conditions:
\begin{align*}
-i \mu_f J_{1,f}^* c_1  + i \mu_s\Delta_1 c_3 & = - i\nu_1\frac{\mu_sP}{\rho_s \omega^2}, \\
J_{0,f}^* c_1- i\omega\Delta_0 c_3 & = (\nu_0 - 1)  \frac{i P}{\rho_s \omega} - \frac{P}{\mu_f k_f^2} .
\end{align*}

The resulting close form solutions are:
\begin{align}
c_1 & = \frac{ - \nu_1 \Delta_0  \frac{\mu_s P}{\rho_s \omega} -
	\mu_s(1-\nu_0)\Delta_1 \frac{P}{\rho_s \omega} +
	i \mu_s \Delta_1 \frac{P}{\mu_f k_f^2}
}{- \mu_f J_{1,f}^*\omega\Delta_0 - i\mu_s J_{0,f}^* \Delta_1} ,
\\
c_3 & = \frac{i \nu_1 J_{0,f}^* \frac{\mu_s P}{\rho_s \omega^2}
	-  \mu_f(1 - \nu_0)J_{1,f}^*\frac{P}{\rho_s \omega} +   iJ_{1,f}^* \frac{P}{k_f^2} }
{- \mu_f J_{1,f}^*\omega\Delta_0 - i \mu_s J_{0,f}^* \Delta_1}.
\end{align}

\end{document}